\documentclass[11pt]{amsart}
\usepackage{graphicx}
\usepackage{amssymb}
\usepackage{amsmath}
\usepackage{amsthm}
\usepackage{mathtools}
\usepackage{tikz}
\usepackage{psfrag}
\usepackage{xcolor}

\newtheorem{thm}{Theorem}[section]
\newtheorem{lem}[thm]{Lemma}
\newtheorem{pro}[thm]{Proposition}

\theoremstyle{definition}
\newtheorem{defn}{Definition}[section]
\newtheorem{remk}{Remark}[section]

\newcommand{\N}{\mathbb N}
\newcommand{\M}{\mathbb M}

\newcommand{\R}{\mathbb R}

\makeatletter
\numberwithin{equation}{section}
\makeatother

\begin{document}



       \author{Seonghak Kim}
       \address{Institute for Mathematical Sciences\\ Renmin University of China \\  Beijing 100872, PRC}
       \email{kimseo14@gmail.com}


       \author{Baisheng Yan}

       \address{Department of Mathematics\\ Michigan State University\\ East Lansing, MI 48824, USA}
       \email{yan@math.msu.edu}



\title{Non-parabolic diffusion problems in one space dimension}


\begin{abstract}
We study some non-parabolic diffusion problems in one-space dimension, where the diffusion  flux exhibits forward and backward nature of the Perona-Malik, H\"ollig or non-Fourier type.  Classical weak solutions to such problems are constructed in a way to capture some expected and unexpected properties, including anomalous asymptotic behaviors and energy dissipation or allocation. Specific properties of solutions will  depend on the type of the diffusion flux, but the primary method of our study relies on reformulating diffusion equations involved as an inhomogeneous partial differential inclusion and on constructing solutions from the differential  inclusion by a combination of the convex integration and Baire's category methods. In doing so, we introduce the appropriate notion of subsolutions of the partial differential inclusion and their transition gauge, which plays a pivotal role in dealing with some specific features of the constructed weak solutions.
\end{abstract}
\maketitle


\section{Introduction}

In this paper, we study the initial-boundary value problem of quasilinear diffusion equations  in one space dimension:
\begin{equation}\label{main-ibP}
\left\{
\begin{array}{ll}
  u_t=(\sigma(u_x))_{x} & \mbox{in $\Omega\times(0,T)$,} \\
  u=u_0 & \mbox{on $\Omega\times\{t=0\}$,} \\
  \sigma(u_x)=0 & \mbox{on $\partial\Omega\times(0,T)$.}
\end{array}
\right.
\end{equation}
Here, $\Omega=(0,L)\subset\R$ with $L>0$, $T\in(0,\infty]$ is given, $\sigma=\sigma(s)$ is the diffusion flux function on $\R$,  $u_0=u_0(x)$ is a given initial function on $\Omega,$ and $u=u(x,t)$ is a solution to the problem.

The diffusion equation in (\ref{main-ibP}) becomes a formal \emph{gradient flow} of the  \emph{energy functional} $\mathcal{E}$ defined   by
\[
\mathcal{E}(p)=\int_0^L W(p_x)\,dx,
\]
where $W=W(s)$ is an antiderivative of the diffusion flux $\sigma.$ When the diffusion flux $\sigma$ is smooth with $\sigma'>0$ in $\R$, problem (\ref{main-ibP}) is parabolic so that  by the classical parabolic theory \cite{LSU, Ln}, it possesses a unique global classical solution that belongs to the parabolic H\"older space $C^{2+\alpha, 1+\frac{\alpha}{2}}(\bar\Omega\times [0,\tau])$ with $\alpha\in (0,1)$ for all $0<\tau<\infty$  if the initial datum $u_0$ is in $ C^{2+\alpha}(\bar\Omega)$ and  satisfies the compatibility condition (\ref{initial-comp}).

However, in many applications and models of diffusion process, the diffusion fluxes $\sigma$ are not increasing on $\R$, as have been studied in phase transition problems in thermodynamics \cite{Ho,HN}, mathematical models of stratified turbulent flows \cite{BBDU}, ecological models of  population dynamics \cite{Ok, Pa}, image processing \cite{PM}, and gradient flows associated with nonconvex energy functionals \cite{Sl}. Therefore, in such cases,  problem (\ref{main-ibP}) becomes {\em non-parabolic}, and the standard methods of parabolic theory are not applicable to study  existence and uniqueness or non-uniqueness of solutions for \emph{arbitrary} initial data $u_0\in C^{2+\alpha}(\bar\Omega)$ satisfying the compatibility condition (\ref{initial-comp}).

For a general flux function $\sigma$, the equation in (\ref{main-ibP}) is in divergence form and hence  a natural \emph{generalized} solution  to problem (\ref{main-ibP}) can be defined  in the usual variational sense.  In the following, we use $W^{1,\infty}_{loc}(\Omega\times [0,T))$ to denote the set of functions $u\in C(\bar\Omega\times[0,T))$ such that $u\in W^{1,\infty}(\Omega\times(0,\tau))$ for all $\tau\in(0,T)$.
A function $u\in W^{1,\infty}_{loc}(\Omega\times[0,T))$  is called a \emph{weak solution} to (\ref{main-ibP}) provided that equality
\begin{equation}\label{def:weak-solution}
\int_0^L \big(u(x,\tau)\zeta(x,\tau)-u_0(x)\zeta(x,0)\big)\,dx = \int_0^\tau\int_0^L \big(u\zeta_t-\sigma(u_x)\zeta_x \big)\,dxdt
\end{equation}
holds for each $\zeta\in C^\infty(\bar\Omega\times[0,T))$ and each $\tau\in[0,T).$ In case of $T=\infty$, such a solution $u$ is said to be \emph{global}.
Upon taking $\zeta\equiv 1$ in (\ref{def:weak-solution}), it is immediate  that every weak solution $u$ to (\ref{main-ibP}) fulfills the \emph{conservation of total mass} at all times:
\[
\int_0^L u(x,t)\,dx=\int_0^L u_0(x)\,dx\;\;\forall t\in[0,T).
\]

In this paper, we focus on the three main types of non-increasing  diffusion flux  functions $\sigma$: the H\"ollig type in \cite{Ho}, the Perona-Malik type in \cite{PM}, and the non-Fourier type with double-well potential $W(s)$.
Therefore, the diffusion fluxes $\sigma\in C(\R)$ under consideration have potential $\R\ni s\mapsto \int_{0}^s\sigma(l)\,dl$ that attains an absolute minimum; throughout the paper, we  let $W\in C^1(\R)$ denote the antiderivative of $\sigma$ whose absolute minimum value is $0$. The initial datum $u_0$ will always be in the class $C^{2+\alpha}(\bar\Omega)$ for a fixed $0<\alpha<1$ and satisfy  the compatibility condition
\begin{equation}\label{initial-comp}
 u_0'(0)=u_0'(L)=0.
\end{equation}

We remark that for diffusion equations of the Perona-Malik type  and certain specific  initial data $u_0$ with  $u_0'$ exhibiting backward phases, the classical solutions were constructed in the works \cite{GG1,GG2}. For such equations with general initial data, it has been known that the Neumann problem may not have a classical solution \cite{KK, Ky}.

The main purpose of this paper is to construct weak solutions to problem (\ref{main-ibP}) with diffusion fluxes $\sigma$ of the three types above that capture some expected and unexpected properties, including anomalous asymptotic behaviors and energy dissipation or allocation.
Our construction relies on the important technique by Zhang \cite{Zh1,Zh2} of  reformulating  diffusion equations involved  as an inhomogeneous partial differential inclusion; such a technique  has been recently generalized by the authors  to handle more general problems, including the high-dimensional ones, through  a combination of the convex integration and Baire's category methods; see \cite{KY, KY1,KY2,KY3}.

However, to construct weak solutions with certain specific features,  as in \cite{KY3}, we need to  introduce appropriate {\em subsolutions} of the associated  partial differential inclusion and their {\em transition gauge} and to iterate some constructions while controlling the gauge.
The controllable gauge, combined with the result on a one-dimensional parabolic problem \cite{Ki2},  will enable us to iterate the constructions and eventually lead us to extract  energy dissipation for the Perona-Malik type diffusions and energy allocation for the non-Fourier type diffusions as well as the blow-up or concentration of the spatial derivative of weak solutions.

Our main results will be introduced in Section \ref{main-res}, where we fix the hypothesis on the diffusion flux $\sigma$ of the H\"ollig, Peorna-Malik and non-Fourier types respectively, state the main results corresponding to each type, and discuss some interesting features of those results.

Other than the statements of the main results, the rest of the paper is organized as follows. In Section \ref{sec:preliminary},  we prepare  all the essential ingredients needed to prove our main results, some of which may be of independent interest; in this section  (see Subsection \ref{sec:proof-classical}), we also provide  the easy proofs for  Theorems \ref{thm:Holllig type_Classical} and \ref{thm:PM type-Classical}.   In Section \ref{sec:proof of Hollig type}, the proof of Theorem \ref{thm:Holllig type} is given to obtain global weak solutions to (\ref{main-ibP}) that are smooth after a finite time when $\sigma$ is of the H\"ollig type. For the same type of solutions to (\ref{main-ibP}) when $\sigma$ is of the Perona-Malik type, the proof of Theorems \ref{thm:PM type-Stable} and \ref{thm:PM type-Stable-Case2} is given in Section \ref{sec:proof-PM-eventual-smooth}. Section \ref{sec:proof-PM-anomalous} is devoted to the proof of Theorem \ref{thm:PM type-Unstable} on the weak solutions to (\ref{main-ibP}) with anomalous blow-up and vanishing energy when $\sigma$ is of the Perona-Malik type. In Section \ref{sec:proof-NF}, global weak solutions to (\ref{main-ibP}) with energy allocation are constructed when $\sigma$ is of the non-Fourier type to prove Theorem \ref{thm:non-Fourier type}. In the last part of the paper, Section \ref{sec:density-proof}, the key density lemma, Lemma \ref{lem:density lemma}, is proved.

\section{Statement of main results}\label{main-res}

We begin this section by fixing some notation. We denote by $\M^{m\times n}$ the space of $m\times n$ real matrices. We write $\Omega_\tau=\Omega\times(0,\tau)$ for each $\tau\in(0,\infty]$ and $\|\cdot\|_\infty=\|\cdot\|_{L^\infty(\Omega)}$. We let $\bar u_0$ denote the mean value of the initial datum $u_0$ over $\Omega$; that is, $\bar u_0=\frac{1}{L}\int_0^L u_0(x)\,dx\in\R.$ We also write
\[
M_0=\max_{\bar\Omega}u_0'\;\;\mbox{and}\;\;m_0=\min_{\bar\Omega}u_0';
\]
so $m_0\le 0\le M_0$. For a bounded open set $U\subset\R^n$ with its coordinates $(z_1,\cdots,z_n)$, $W^{1,\infty}(U)$ denotes the space of functions $f\in L^\infty(U)$ with  $f_{z_i}\in L^\infty(U)$  $(i=1,\cdots,n)$ and norm $\|f\|_{W^{1,\infty}(U)}:=\|f\|_{L^\infty(U)}+\sum_{i=1}^n \|f_{z_i}\|_{L^\infty(U)}$, and $W^{1,\infty}_0(U)$ consists of the functions $f\in W^{1,\infty}(U)$ with zero trace on $\partial U.$
For integers $k,l\ge 0$ with $2l\le k$, we denote by $C^{k,l}(\bar\Omega_\tau)$, with $0<\tau\le\infty$, the space of functions $u\in C(\bar\Omega_\tau)$ such that $D^i_x D^j_t u\in C(\bar\Omega_\tau)$ for all integers $0\le i\le k$ and  $0\le j\le l$ with $i+2j\le k$.
For other function spaces, we mainly follow the notations in the book \cite{LSU}, with one exception  that the letter $C$ is used instead of $H$ regarding suitable parabolic H\"older spaces.

\subsection{H\"ollig type equations}
In this subsection, we consider a class of flux functions $\sigma$ that includes piecewise linear ones, studied by H\"ollig \cite{Ho}, and present the results on large time behaviors of global weak or classical solutions to problem (\ref{main-ibP}).

\begin{figure}[ht]
\begin{center}
\includegraphics[scale=0.6]{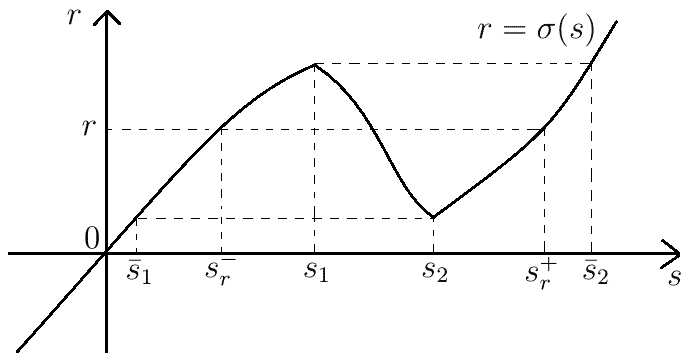}
\end{center}
\caption{Diffusion flux $\sigma(s)$ of the H\"ollig type}
\label{fig1}
\end{figure}

To be precise, we impose the following conditions on the flux function $\sigma\in C(\R)$ (see Figure \ref{fig1}): there exist three numbers $\bar s_2>s_2>s_1>0$ such that
\begin{itemize}
\item [(a)] $\sigma\in C^3((-\infty,s_1)\cup(s_2,\infty)),$
\item [(b)] $\sigma>0$ on $(s_1,s_2)$, $\sigma'>0$ in $(-\infty,s_1)\cup(s_2,\infty)$, and
\item [(c)] $\sigma(\bar s_2)=\sigma(s_1)>\sigma(s_2)>\sigma(0)=0$.
\end{itemize}
Let $\bar s_1\in(0,s_1)$ denote the unique number with $\sigma(\bar s_1)=\sigma(s_2)$.
For each $r\in(\sigma(s_2),\sigma(s_1))$, let $s^+_r\in(s_2,\bar s_2)$ and $s^-_r\in(\bar s_1,s_1)$ be the unique numbers such that $\sigma(s^\pm_r)=r$.

\subsubsection{Global classical solution}

The following is a result on the global classical solution $u$ to problem (\ref{main-ibP}) when the maximum $M_0$ of $u_0'$ on $\bar\Omega$ is less than the threshold $s_1$.

\begin{thm}[Global classical solution]\label{thm:Holllig type_Classical}
Assume
\[
M_0<s_1.
\]
Then there exists a unique global solution $u\in C^{2,1}(\bar\Omega_\infty)$ to problem (\ref{main-ibP}) satisfying the following properties:
\begin{itemize}
\item[(1)] for each $\tau>0$, $u\in C^{2+\alpha,1+\alpha/2}(\bar\Omega_\tau)$,
\item[(2)] for all $t_2>t_1\ge0,$
\[
\min_{\bar\Omega} u(\cdot,t_1)\le \min_{\bar\Omega} u(\cdot,t_2)\le \bar{u}_0\le \max_{\bar\Omega} u(\cdot,t_2) \le \max_{\bar\Omega} u(\cdot,t_1),
\]
\item[(3)] for all $t_2>t_1\ge 0,$
\[
\min_{\bar\Omega} u_x(\cdot,t_1)\le \min_{\bar\Omega} u_x(\cdot,t_2)\le0\le \max_{\bar\Omega} u_x(\cdot,t_2) \le \max_{\bar\Omega} u_x(\cdot,t_1),
\]
\item[(4)] $\|u(\cdot,t)-\bar{u}_0\|_{W^{1,\infty}(\Omega)} \le Ce^{-\gamma t}$ for all $t\ge0$, where  $C>0$ and $\gamma>0$ are some constants depending only on $u_0$, $\sigma$ and $L$.
\end{itemize}
\end{thm}

\subsubsection{Global weak solutions that are smooth after a finite time}

When $M_0$ is greater than the threshold $\bar{s}_1$, we can get infinitely many global weak solutions $u$ to problem (\ref{main-ibP}) that become identical and smooth after a finite time.

\begin{thm}[Eventual smoothing]\label{thm:Holllig type}
Assume
\[
M_0>\bar s_1
\]
so that $\bar s_1<u_0'(x_0)<\bar s_2$ for some $x_0\in\Omega.$ Let $\sigma(s_2)<r_1<r_2<\sigma(s_1)$ be any two numbers such that $s^-_{r_1}<u_0'(x_0)<s^+_{r_2}$. Then there exists a function $u^*\in C^{2,1}(\bar\Omega_\infty)$ such that for each $\epsilon>0$, there are infinitely many global weak solutions $u$ to problem (\ref{main-ibP}) satisfying the following properties:
letting
\[
\begin{split}
Q_1 & =\big\{(x,t)\in\Omega_\infty\,|\, s^-_{r_1}<u^*_x(x,t)<s^+_{r_2}\big\}, \\
Q_2 & =\big\{(x,t)\in\Omega_\infty\,|\, u^*_x(x,t)>s^+_{r_2}\big\}, \\
F^+ & =\big\{(x,t)\in\Omega_\infty\,|\, u^*_x(x,t)=s^+_{r_2}\big\}, \\
F^- & =\big\{(x,t)\in\Omega_\infty\,|\, u^*_x(x,t)=s^-_{r_1}\big\},
\end{split}
\]
we have
\begin{itemize}
\item [(1)] $u=u^*$ in $\bar\Omega_\infty\setminus Q_1$, $\|u-u^*\|_{L^\infty(Q_1)}<\epsilon$, $\|u_t-u_t^*\|_{L^\infty(Q_1)}<\epsilon$,
\item [(2)] $u_x\in[s^-_{r_1},s^-_{r_2}]\cup[s^+_{r_1},s^+_{r_2}]$ a.e. in $Q_1$, $u_x=s^+_{r_2}$ a.e. in $F^+$, $u_x=s^-_{r_1}$ a.e. in $F^-$,
\item [(3)] $\big|\big\{(x,t)\in Q_1\,|\, u_x(x,t)\in[s^\pm_{r_1},s^\pm_{r_2}]\big\}\big|>0$,
\item [(4)] $Q_1\cup Q_2\cup F^+\cup F^-$ is a bounded set whose closure is contained in $\Omega\times[0,\infty)$,
\item [(5)] $0\le \underset{(x,t)\in Q_2}{\sup}t<\underset{(x,t)\in Q_1}{\sup}t<\infty$,
\item [(6)] for each $\tau>0$, $u^*\in C^{2+\alpha,1+\alpha/2}(\bar\Omega_\tau)$,
\item [(7)] for all $t_2>t_1\ge0,$
\[
\min_{\bar\Omega} u^*(\cdot,t_1)\le \min_{\bar\Omega} u^*(\cdot,t_2)\le \bar{u}_0 \le\max_{\bar\Omega} u^*(\cdot,t_2) \le \max_{\bar\Omega} u^*(\cdot,t_1),
\]
\item [(8)] for all $t_2>t_1\ge 0,$
\[
\min_{\bar\Omega} u^*_x(\cdot,t_1)\le \min_{\bar\Omega} u^*_x(\cdot,t_2)\le0\le \max_{\bar\Omega} u^*_x(\cdot,t_2) \le \max_{\bar\Omega} u^*_x(\cdot,t_1),
\]
\item [(9)] $\|u^*(\cdot,t)-\bar{u}_0\|_{W^{1,\infty}(\Omega)} \le Ce^{-\gamma t}$ for all $t\ge0$, where $C>0$ and $\gamma>0$ are some constants depending only on $u_0$, $\sigma$, $r_1$, $r_2$ and $L$.
\end{itemize}
\end{thm}

\subsubsection{Remarks}

Depending on the maximum value $M_0$ of $u_0'$ on $\bar\Omega$, we may have several possible scenarios on the existence of classical solutions to problem (\ref{main-ibP}).
\begin{itemize}
\item ({\em Coexistence}) When $\bar{s}_1<M_0<s_1$, from Theorems \ref{thm:Holllig type_Classical} and \ref{thm:Holllig type}, we have the \emph{coexistence} of a unique global classical solution and infinitely many global weak (but eventually classical) solutions to problem (\ref{main-ibP}). This pathological feature of the existence theory seems to arise from forward and backward nature of the diffusion flux $\sigma$ of the H\"ollig type although, in this case, the classical solution may be the most natural representative.

\item ({\em Critical case}) The case $M_0=s_1$ is rather subtle. If, in addition, $\sigma\in C^3((-\infty,s_1])$ and $\sigma'(s_1)>0$, we can get a unique global classical solution $u$ to problem (\ref{main-ibP}) as in Theorem \ref{thm:Holllig type_Classical}. On the other hand, if, in addition, $\sigma\in C^3((-\infty,s_1])$ and $\sigma'(s_1)=0$, one may have to rely on the \emph{degenerate} parabolic theory, which we will not address here. However, no matter what we assume for the diffusion flux $\sigma$ at $s=s_1$, by Theorem \ref{thm:Holllig type}, we still have infinitely many global weak solutions to (\ref{main-ibP}) that are smooth after a finite time.

\item ({\em Nonexistence of classical solution}) If $M_0>s_1$, there can be no classical solution to problem (\ref{main-ibP}) at all for a large class of initial data $u_0$. For example, let us assume in addition that $\sigma\in C^\infty(\R)$ and that $\sigma'<0$ in $(s_1,s_2)$. Since $M_0>s_1$ and $u_0'=0$ on $\partial\Omega$, there is a nonempty open interval $I\subset\Omega$ on which $u_0'$ falls into the backward regime $(s_1,s_2)$. So one can readily check as in \cite{Ky} that there is no classical solution (not even a local $C^1$ solution) to (\ref{main-ibP}) unless $u_0$ were infinitely differentiable from the start. Even in the case that $u_0$ is infinitely differentiable, it is not at all obvious if there is a local or global classical solution to (\ref{main-ibP}); in this regard, see \cite{GG1,GG2}. Regardless of such possible nonexistence of a classical solution, Theorem \ref{thm:Holllig type} still guarantees the existence of global weak solutions to (\ref{main-ibP}) that are smooth after a finite time.

\item ({\em Breakdown of uniqueness}) From the first discussion above, one may expect that if $M_0\le\bar{s}_1$, there exists a unique global weak (and necessarily classical) solution to problem (\ref{main-ibP}). Regarding this, we define the number
\[
s_{cr}=\sup\left\{s\ge0\,\bigg|\,\begin{array}{l}
                           \mbox{$\forall u_0\in C^{2+\alpha}(\bar\Omega)$ with $u_0'=0$ on $\partial\Omega$} \\
                           \mbox{and $\max_{\bar\Omega}u_0'\le s,$ there exists a unique } \\
                           \mbox{global weak solution to (\ref{main-ibP})}
                         \end{array}
\right\}.
\]
From Theorem \ref{thm:Holllig type}, it is easy to see that $s_{cr}\le\bar{s}_1$. Following \cite{HN, KY2}, one can also check that $s_{cr}>0$. Nonetheless, it seems not simple to answer if  $s_{cr}=\bar{s}_1$ or not.
\end{itemize}

\subsection{Perona-Malik type equations}

Here, we consider a class of flux functions $\sigma$ that contains the \emph{Perona-Malik functions} in \cite{PM}:
\[
\sigma(s)=\frac{s}{1+s^2}\;\;\mbox{and}\;\;\sigma(s)=se^{-s^2/2}.
\]

\begin{figure}[ht]
\begin{center}
\includegraphics[scale=0.6]{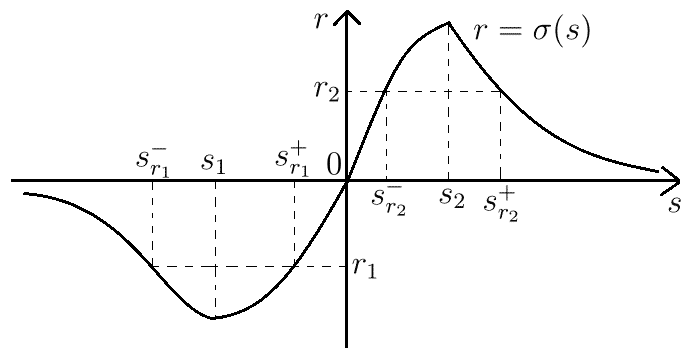}
\end{center}
\caption{Diffusion flux $\sigma(s)$ of the Perona-Malik type}
\label{fig2}
\end{figure}

Specifically, we assume the following conditions on the flux function $\sigma\in C(\R)$ (see Figure \ref{fig2}): there exist two numbers $s_1<0<s_2$ such that
\begin{itemize}
\item [(a)] $\sigma\in C^3(s_1,s_2),$
\item [(b)] $\sigma'>0$ in $(s_1,s_2)$, $\sigma(0)=0$,
\item [(c)] $\sigma$ is strictly decreasing on $(-\infty,s_1]\cup[s_2,\infty)$, and
\item [(d)] $\underset{s\to\pm\infty}{\lim}\sigma(s)=0.$
\end{itemize}
For each $r\in(0,\sigma(s_2))$, let $s^+_{r}\in(s_2,\infty)$ and $s^-_{r}\in(0,s_2)$ denote the unique numbers with $\sigma(s^\pm_{r})=r$. Similarly, for each $r\in(\sigma(s_1),0)$, let $s^+_{r}\in(s_1,0)$ and $s^-_{r}\in(-\infty,s_1)$ be the unique numbers such that $\sigma(s^\pm_{r})=r$.

\begin{figure}[ht]
\begin{center}
\includegraphics[scale=0.6]{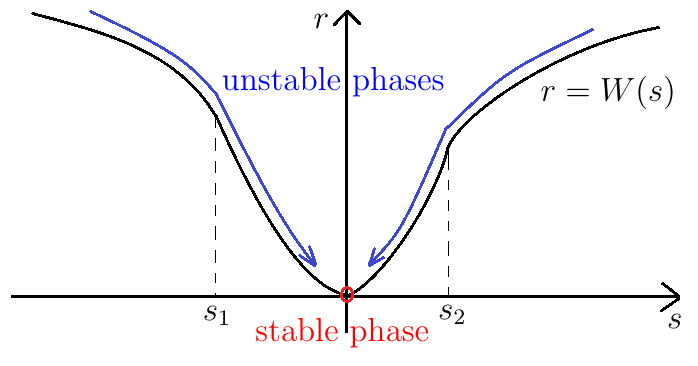}
\end{center}
\caption{Potential function $W(s)$ of the Perona-Malik type}
\label{fig2-1}
\end{figure}

\subsubsection{Global classical solution}

Since the flux function $\sigma(s)$ has a positive derivative for values $s$ near $0$, similarly as in Theorem \ref{thm:Holllig type_Classical}, we have a global classical solution $u$ to problem (\ref{main-ibP}) for initial data $u_0$ with sufficiently small $M_0$ and $m_0$.

\begin{thm}[Global classical solution]\label{thm:PM type-Classical}
Assume
\[
s_1<m_0\le0\le M_0<s_2.
\]
Then there exists a unique global solution $u\in C^{2,1}(\bar\Omega_\infty)$ to problem (\ref{main-ibP}) satisfying the following properties:
\begin{itemize}
\item[(1)] for each $\tau>0$, $u\in C^{2+\alpha,1+\alpha/2}(\bar\Omega_\tau)$,
\item[(2)] for all $t_2>t_1\ge0,$
\[
\min_{\bar\Omega} u(\cdot,t_1)\le \min_{\bar\Omega} u(\cdot,t_2)\le \bar{u}_0\le \max_{\bar\Omega} u(\cdot,t_2) \le \max_{\bar\Omega} u(\cdot,t_1),
\]
\item[(3)] for all $t_2>t_1\ge 0,$
\[
\min_{\bar\Omega} u_x(\cdot,t_1)\le \min_{\bar\Omega} u_x(\cdot,t_2)\le0\le \max_{\bar\Omega} u_x(\cdot,t_2) \le \max_{\bar\Omega} u_x(\cdot,t_1),
\]
\item[(4)] $\|u(\cdot,t)-\bar{u}_0\|_{W^{1,\infty}(\Omega)} \le Ce^{-\gamma t}$ for all $t\ge0$, where  $C>0$ and $\gamma>0$ are some constants depending only on $u_0$, $\sigma$ and $L$.
\end{itemize}
\end{thm}

\subsubsection{Global weak solutions that are smooth after a finite time}

For any \emph{nonconstant} smooth initial datum $u_0$, we have multiple global weak solutions $u$ to problem (\ref{main-ibP}) that are classical after a finite time as in the following theorems. Note that $\{0\}\subsetneq[m_0,M_0]$ when $u_0$ is nonconstant.

The first one deals with the case that none of $M_0$ and $m_0$ are equal to $0$.

\begin{thm}[Eventual smoothing: Case 1]\label{thm:PM type-Stable}
Assume
\[
m_0<0<M_0.
\]
Let $\sigma(s_1)<r_1<0<r_2<\sigma(s_2)$ be any two numbers such that
\[
s^-_{r_1}<m_0<s^+_{r_1}<0<s^-_{r_2}<M_0<s^+_{r_2},
\]
and let $r_1'\in(\sigma(s_1),r_1)$, $r_2'\in(r_2,\sigma(s_2))$ be any two numbers.
Then there exists a function $u^*\in C^{2,1}(\bar\Omega_\infty)$ such that for each $\epsilon>0$, there are infinitely many global weak solutions $u$ to problem (\ref{main-ibP}) satisfying the following properties:
letting
\[
\begin{split}
Q_1 & =\big\{(x,t)\in\Omega_\infty\,|\, s^-_{r_1}<u^*_x(x,t)<s^+_{r_1}\big\}, \\
Q_2 & =\big\{(x,t)\in\Omega_\infty\,|\, s^-_{r_2}<u^*_x(x,t)<s^+_{r_2}\big\}, \\
F_1 & =\big\{(x,t)\in\Omega_\infty\,|\, u^*_x(x,t)=s^+_{r_1}\big\}, \\
F_2 & =\big\{(x,t)\in\Omega_\infty\,|\, u^*_x(x,t)=s^-_{r_2}\big\},
\end{split}
\]
we have
\begin{itemize}
\item [(1)] $u=u^*$ in $\bar\Omega_\infty\setminus (Q_1\cup Q_2)$, $\|u-u^*\|_{L^\infty(Q_1\cup Q_2)}<\epsilon$, $\|u_t-u_t^*\|_{L^\infty(Q_1\cup Q_2)}<\epsilon$,
\item [(2)] $ \left\{ \begin{array}{l}
                         \mbox{$u_x\in[s^-_{r_1},s^-_{r_1'}]\cup[s^+_{r_1'},s^+_{r_1}]$ a.e. in $Q_1$, $u_x=s^+_{r_1}$ a.e. in $F_1$,} \\
                         \mbox{$u_x\in[s^-_{r_2},s^-_{r_2'}]\cup[s^+_{r_2'},s^+_{r_2}]$ a.e. in $Q_2$, $u_x=s^-_{r_2}$ a.e. in $F_2$,}
                       \end{array} \right.
                       $
\item [(3)] $\left\{\begin{array}{l}
                      \mbox{$\big|\big\{(x,t)\in Q_1\,|\, u_x(x,t)\in[s^-_{r_1},s^-_{r_1'}]\big\}\big|>0$,} \\
                      \mbox{$\big|\big\{(x,t)\in Q_1\,|\, u_x(x,t)\in[s^+_{r_1'},s^+_{r_1}]\big\}\big|>0$,} \\
                      \mbox{$\big|\big\{(x,t)\in Q_2\,|\, u_x(x,t)\in[s^-_{r_2},s^-_{r_2'}]\big\}\big|>0$,} \\
                      \mbox{$\big|\big\{(x,t)\in Q_2\,|\, u_x(x,t)\in[s^+_{r_2'},s^+_{r_2}]\big\}\big|>0$,}
                    \end{array}
                  \right. $
\item [(4)] $Q_1\cup Q_2\cup F_1\cup F_2$ is a bounded set whose closure is contained in $\Omega\times[0,\infty)$,
\item [(5)] for each $\tau>0$, $u^*\in C^{2+\alpha,1+\alpha/2}(\bar\Omega_\tau)$,
\item [(6)] for all $t_2>t_1\ge0,$
\[
\min_{\bar\Omega} u^*(\cdot,t_1)\le \min_{\bar\Omega} u^*(\cdot,t_2)\le\bar u_0\le \max_{\bar\Omega} u^*(\cdot,t_2) \le \max_{\bar\Omega} u^*(\cdot,t_1),
\]
\item [(7)] for all $t_2>t_1\ge 0,$
\[
\min_{\bar\Omega} u^*_x(\cdot,t_1)\le \min_{\bar\Omega} u^*_x(\cdot,t_2)\le0\le \max_{\bar\Omega} u^*_x(\cdot,t_2) \le \max_{\bar\Omega} u^*_x(\cdot,t_1),
\]
\item [(8)] $\|u^*(\cdot,t)-\bar{u}_0\|_{W^{1,\infty}(\Omega)} \le Ce^{-\gamma t}$ for all $t\ge0$, where $C>0$ and $\gamma>0$ are some constants depending only on $u_0$, $\sigma$, $r_1$, $r_1'$, $r_2$, $r_2'$ and $L$.
\end{itemize}
\end{thm}

Unlike Case 1 above, when $m_0=0$, our global weak solutions $u$ to problem (\ref{main-ibP}) have no negative spatial derivative at all times as follows.

\begin{thm}[Eventual smoothing: Case 2]\label{thm:PM type-Stable-Case2}
Assume
\[
0=m_0<M_0.
\]
Let $0<r_0<\sigma(s_2)$ be any number such that
\[
s^-_{r_0}<M_0<s^+_{r_0},
\]
and let $r_0'\in(r_0,\sigma(s_2))$ be any number.
Then there exists a function $u^*\in C^{2,1}(\bar\Omega_\infty)$ such that for each $\epsilon>0$, there are infinitely many global weak solutions $u$ to problem (\ref{main-ibP}) satisfying the following properties:
letting
\[
\begin{split}
Q_0 & =\big\{(x,t)\in\Omega_\infty\,|\, s^-_{r_0}<u^*_x(x,t)<s^+_{r_0}\big\}, \\
F_0 & =\big\{(x,t)\in\Omega_\infty\,|\, u^*_x(x,t)=s^-_{r_0}\big\},
\end{split}
\]
we have
\begin{itemize}
\item [(1)] $u=u^*$ in $\bar\Omega_\infty\setminus Q_0$, $\|u-u^*\|_{L^\infty(Q_0)}<\epsilon$, $\|u_t-u_t^*\|_{L^\infty(Q_0)}<\epsilon$,
\item [(2)] $u_x\in[s^-_{r_0},s^-_{r_0'}]\cup[s^+_{r_0'},s^+_{r_0}]$ a.e. in $Q_0$, $u_x=s^-_{r_0}$ a.e. in $F_0$,
\item [(3)] $\left\{\begin{array}{l}
                      \mbox{$\big|\big\{(x,t)\in Q_0\,|\, u_x(x,t)\in[s^-_{r_0},s^-_{r_0'}]\big\}\big|>0$,} \\
                      \mbox{$\big|\big\{(x,t)\in Q_0\,|\, u_x(x,t)\in[s^+_{r_0'},s^+_{r_0}]\big\}\big|>0$,}
                    \end{array}
                  \right. $
\item [(4)] $Q_0\cup F_0$ is a bounded set whose closure is contained in $\Omega\times[0,\infty)$,
\item [(5)] for each $\tau>0$, $u^*\in C^{2+\alpha,1+\alpha/2}(\bar\Omega_\tau)$,
\item [(6)] for all $t_2>t_1\ge0,$
\[
\min_{\bar\Omega} u^*(\cdot,t_1)\le \min_{\bar\Omega} u^*(\cdot,t_2)\le\bar u_0\le \max_{\bar\Omega} u^*(\cdot,t_2) \le \max_{\bar\Omega} u^*(\cdot,t_1),
\]
\item [(7)] for all $t_2>t_1\ge 0,$
\[
0=\min_{\bar\Omega} u^*_x(\cdot,t_1)= \min_{\bar\Omega} u^*_x(\cdot,t_2)\le \max_{\bar\Omega} u^*_x(\cdot,t_2) \le \max_{\bar\Omega} u^*_x(\cdot,t_1),
\]
\item [(8)] $\|u^*(\cdot,t)-\bar{u}_0\|_{W^{1,\infty}(\Omega)} \le Ce^{-\gamma t}$ for all $t\ge0$, where $C>0$ and $\gamma>0$ are some constants depending only on $u_0$, $\sigma$, $r_0$, $r_0'$ and $L$.
\end{itemize}
\end{thm}

We skip the details for Case 3: $M_0=0$ as those are similar to Theorem \ref{thm:PM type-Stable-Case2} in a symmetric way; in this case, our global weak solutions $u$ to problem (\ref{main-ibP}) do not have any positive spatial derivative at all times.

\subsubsection{Blow-up solutions with vanishing energy}

We now present rather surprising weak solutions $u$ to problem (\ref{main-ibP}) that exhibit abnormal behaviors of evolution. Such a solution $u$ as we can see below shares some similar properties of a global classical solution to the strictly parabolic problem in Theorem \ref{thm:uniform-parabolic}: the solution $u$ stabilizes to $\bar{u}_0$ in terms of the $L^\infty(\Omega)$-norm as $t$ increases to some number $T\in(0,\infty]$, and its energy $\mathcal{E}(u)$ vanishes in a certain sense at the same time (see Figure \ref{fig2-1}). On the other hand, the spatial derivative $u_x$ blows up in some sense as $t$ increases to $T$.

\begin{thm}[Anomalous blow-up with vanishing energy]\label{thm:PM type-Unstable}
Let $u_0$ be nonconstant; then there exist a time $T\in(0,\infty]$, an open set $G$ in $\R^2$ with $\partial\Omega\times[0,T)\subset G$,
a strictly increasing sequence $t_j\in(0,T)$ converging to $T$, and infinitely many weak solutions $u$ to problem (\ref{main-ibP}) satisfying the following properties:
\begin{itemize}
\item [(1)] for each $\tau\in(0,T)$, $u\in C^{2+\alpha,1+\alpha/2}(\overline{G\cap\Omega_\tau})$,
\item [(2)] $\|u(\cdot,t)-\bar{u}_0\|_{\infty}\to 0$ as $t\nearrow T,$
\item [(3)] $\|u_x\|_{L^\infty(\Omega\times(t_j,t_{j+1}))}\to\infty$ as $j\to\infty$,
\item [(4)] $\frac{1}{t_{j+1}-t_j}\int_{t_j}^{t_{j+1}}\mathcal{E}(u(\cdot,t))\,dt\to 0$ as $j\to\infty$.
\end{itemize}
\end{thm}

\begin{remk}
It should be noted from Theorems \ref{thm:PM type-Stable}, \ref{thm:PM type-Stable-Case2} and \ref{thm:PM type-Unstable} that for any \emph{nonconstant} initial datum $u_0\in C^{2+\alpha}(\bar\Omega)$ with $u'_0=0$ on $\partial\Omega$, there are global weak solutions to problem (\ref{main-ibP}) that are smooth after a finite time and weak solutions to (\ref{main-ibP}) whose spatial derivative blows up at some time $0<T\le\infty$. Despite of this discrepancy, solutions of each type uniformly converge to $\bar u_0$ and have a vanishing energy in a certain sense as time increases to the terminal value. Moreover, if $M_0$ and $m_0$ are close enough to $0$, it follows from Theorem \ref{thm:PM type-Classical} that (\ref{main-ibP}) also possesses a unique global classical solution; thus, in this case, there are \emph{at least} three types of solutions to the problem for the common initial datum $u_0$.
\end{remk}

\subsection{Non-Fourier type equations}

In this last subsection, we consider a class of flux functions $\sigma$ that violate the Fourier inequality: $s\sigma(s)\ge 0$  $\forall s\in\R$.

\begin{figure}[ht]
\begin{center}
\includegraphics[scale=0.6]{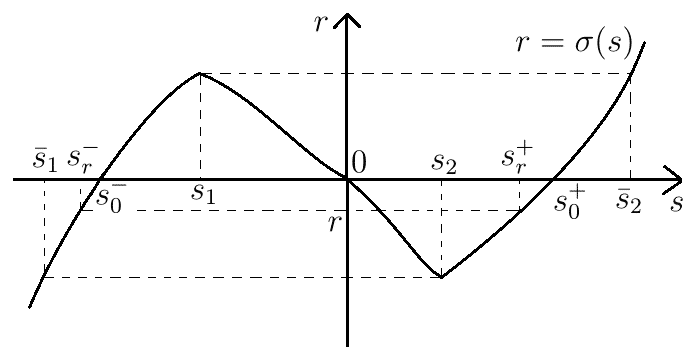}
\end{center}
\caption{Diffusion flux $\sigma(s)$ of the non-Fourier type}
\label{fig3}
\end{figure}

More precisely, we assume the following conditions on the flux function $\sigma\in C(\R)$ (see Figure \ref{fig3}): there exist four numbers $\bar s_1<s_1<0<s_2<\bar{s}_2$ such that
\begin{itemize}
\item [(a)] $\sigma\in C^3((-\infty,s_1)\cup(s_2,\infty)),$
\item [(b)] $\sigma'>0$ in $(-\infty,s_1)\cup(s_2,\infty)$,
\item [(c)]  $\sigma(\bar s_1)=\sigma(s_2)<\sigma(0)=0<\sigma(s_1)=\sigma(\bar{s}_2)$, and
\item [(d)] $\sigma$ has exactly three zeros.
\end{itemize}
For each $r\in(\sigma(s_2),\sigma(s_1))$, let $s^+_{r}\in(s_2,\bar s_2)$ and $s^-_{r}\in(\bar s_1,s_1)$ denote the unique numbers with $\sigma(s^\pm_{r})=r$. Note that the Fourier inequality is violated by $\sigma$ as follows:
\[
s\sigma(s)\left\{\begin{array}{l}
                   \mbox{$=0$ on $\{s^+_0,0,s^-_0\}$}, \\
                   \mbox{$>0$ in $(-\infty,s^-_0)\cup(s^+_0,\infty)$}, \\
                   \mbox{$<0$ in $(s^-_0,s^+_0)\setminus\{0\}$}.
                 \end{array} \right.
\]

\begin{figure}[ht]
\begin{center}
\includegraphics[scale=0.6]{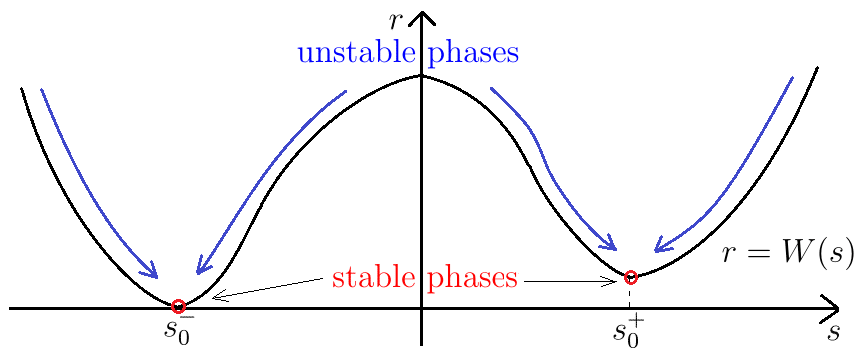}
\end{center}
\caption{Double-well potential $W(s)$ of the non-Fourier type}
\label{fig3-1}
\end{figure}

We now formulate the result on global weak solutions $u$ to problem (\ref{main-ibP}). An interesting feature of such a solution $u$ is that while it stabilizes to $\bar u_0$ in terms of the $L^\infty(\Omega)$-norm as $t\to\infty$, its spatial derivative $u_x$ concentrates on the set $\{s^+_0,s^-_0\}$ in such a way that the energy $\mathcal{E}(u)$ converges in some sense to a specific number reflecting the allocation of double-well energy  by the nonzero roots $s^\pm_0$ of $\sigma$ (see Figure \ref{fig3-1}).

\begin{thm}[Allocation  of double-well energy]\label{thm:non-Fourier type}
There exist a strictly increasing sequence $0<t_j\to\infty$ and infinitely many global weak solutions $u$ to problem (\ref{main-ibP}) satisfying the following properties:
\begin{itemize}
\item [(1)] $\|u(\cdot,t)-\bar{u}_0\|_{\infty}\to 0$ as $t\to\infty,$
\item [(2)] $\big\|\mathrm{dist}\big(u_x(\cdot,t),\{s_0^+,s_0^-\} \big)\big\|_\infty\to 0$ as $t\to\infty$ outside a null set in $(0,\infty)$,
\item [(3)] $\frac{1}{t_{j+1}-t_j}\int_{t_j}^{t_{j+1}}\mathcal{E}(u(\cdot,t)) \,dt\to L\big(\lambda_0 W(s^+_0)+(1-\lambda_0)W(s^-_0)\big)$ as $j\to\infty$, where $\lambda_0:=\frac{-s^-_0}{s^+_0-s^-_0}\in(0,1)$.
\end{itemize}
\end{thm}

Note that every constant function in $\Omega$ is a steady state solution to problem (\ref{main-ibP}), which is \emph{unstable} since its energy $LW(0)$ involves the local maximum value $W(0)>\max\{W(s^+_0),W(s^-_0)\}$. However, by the above result, even for any constant initial datum, we still have global weak solutions to (\ref{main-ibP}) with energy converging to a more \emph{stable} level $L\big(\lambda_0 W(s^+_0)+(1-\lambda_0)W(s^-_0)\big).$

\section{Preliminaries and easy proofs}\label{sec:preliminary}

In this section, we prepare some useful results that are essential ingredients for proving our main results, and some of which may be of independent interest. In particular, as an easy consequence of the result in Subsection \ref{sec:parabolic}, we provide the  proof of  Theorems \ref{thm:Holllig type_Classical} and \ref{thm:PM type-Classical} in Subsection \ref{sec:proof-classical} concerning the existence and asymptotic behaviors of a unique global classical solution $u$ to problem (\ref{main-ibP}) for suitable initial data $u_0$ when the diffusion flux $\sigma$ is either of the H\"ollig type or of the Perona-Malik type.

\subsection{Parabolic equations}\label{sec:parabolic}

Here, we consider large time behaviors of a global classical solution to problem (\ref{main-ibP}) when the problem is \emph{strictly parabolic}. The following theorem will play the role of a major building block for constructing weak solutions of interest in our main results. Since it can be proved easily with the help of \cite{Ki2}, we include its proof here.

\begin{thm}\label{thm:uniform-parabolic}
Let $\sigma\in C^3(\R)$ satisfy that $\sigma'>0\;\;\mbox{in $\R$}.$
Then there exists a unique global solution $u\in C^{2,1}(\bar\Omega_\infty)$ to problem (\ref{main-ibP})  satisfying the following properties:
\begin{itemize}
\item[(1)] for each $\tau>0$, $u\in C^{2+\alpha,1+\alpha/2}(\bar\Omega_\tau)$,
\item[(2)] for all $t_2>t_1\ge0,$
\[
\min_{\bar\Omega} u(\cdot,t_1)\le \min_{\bar\Omega} u(\cdot,t_2)\le \bar{u}_0\le \max_{\bar\Omega} u(\cdot,t_2) \le \max_{\bar\Omega} u(\cdot,t_1),
\]
\item[(3)] for all $t_2>t_1\ge 0,$
\[
\min_{\bar\Omega} u_x(\cdot,t_1)\le \min_{\bar\Omega} u_x(\cdot,t_2)\le0\le \max_{\bar\Omega} u_x(\cdot,t_2) \le \max_{\bar\Omega} u_x(\cdot,t_1),
\]
\item[(4)] $\|u(\cdot,t)-\bar{u}_0\|_{W^{1,\infty}(\Omega)} \le Ce^{-\gamma t}$ for all $t\ge0$, where  $0<\kappa<1$, $\lambda>0$,  $m>1$,
\[
\theta:=\min_{[-\|u_0'\|_\infty\frac{m}{m-1}, \|u_0'\|_\infty\frac{m}{m-1}]}\sigma',\;\;
\tilde\theta:=\max_{[-\|u_0'\|_\infty\frac{m}{m-1}, \|u_0'\|_\infty\frac{m}{m-1}]}|\sigma''|,
\]
\[
\gamma:=\frac{\kappa\theta\lambda^2e^{-\lambda L}}{\max\Big\{\frac{\|u_0'\|_\infty \tilde\theta}{(1-\kappa)\theta}+1,m \Big\}-e^{-\lambda L}},
\]
and $C>0$ is some constant depending only on $u_0$, $\sigma$, $L$, $\kappa$, $\lambda$ and $m$.
\end{itemize}
\end{thm}

\begin{proof}
To start the proof, let us fix any numbers $0<\kappa<1$, $\lambda>0$ and $m>1$. Set $s_0=\|u_0'\|_{\infty}\frac{m}{m-1}+1$, and choose a function $\tilde\sigma\in C^3(\R)$ such that
\begin{equation}\label{uniform-parabolic-1}
\left\{\begin{array}{l}
         \tilde\sigma(s)=\sigma(s)\;\;\forall s\in[-s_0,s_0], \\
         \lambda\le\tilde\sigma'(s)\le\Lambda\;\;\forall s\in\R, \\
         \mbox{$\tilde\sigma'''$ is bounded in $\R$},
       \end{array}
 \right.
\end{equation}
where $\Lambda\ge\lambda>0$ are some constants. Then it follows from \cite[Theorem 13.24]{Ln} that there exists a global solution $u\in C^{2,1}(\bar\Omega_\infty)$ to problem (\ref{main-ibP}), with $\sigma$ replaced by $\tilde\sigma$, such that $u\in C^{2+\alpha,1+\alpha/2}(\bar\Omega_\tau)$ for each $\tau>0$.

We now check that $u$ is a unique global solution to problem (\ref{main-ibP}), with original $\sigma$, satisfying (1)--(4).

Note first that (1) is satisfied as above.
Next, according to \cite[Lemma 2.2]{Ki1}, we have an improved interior regularity of $u$ as $u\in C^{3+\beta,\frac{3+\beta}{2}}(\Omega\times(0,\tau])$ for all $\tau>0$, where $\beta\in(0,1)$ is some number. In particular, $\rho:=u_x\in C^{1,0}(\bar\Omega_\infty)\cap C^{2,1}(\Omega_\infty)$ solves the initial-boundary value problem:
\begin{equation}\label{uniform-parabolic-2}
\left\{
\begin{array}{ll}
  \rho_t=(\tilde\sigma(\rho))_{xx} & \mbox{in $\Omega\times(0,\infty)$,} \\
  \rho=\rho_0 & \mbox{on $\Omega\times\{t=0\}$,} \\
  \rho=0 & \mbox{on $\partial\Omega\times(0,\infty)$},
\end{array}\right.
\end{equation}
where $\rho_0:=u_0'\in C^1(\bar\Omega)$. So \cite[Lemma 2.1]{Ki2} implies (3). Thanks to (\ref{uniform-parabolic-1}), it now follows from (3) that $u$ is a global solution to problem (\ref{main-ibP}) with original $\sigma$.
Also, from (\ref{uniform-parabolic-1}) and (\ref{uniform-parabolic-2}), it is easy to see that  \cite[Theorem 1.1]{Ki2} together with Poincar\'e's inequality yields (4). Finally, uniqueness of such a solution $u$ to (\ref{main-ibP}) and (2) are a consequence of \cite[Propositions 2.3 and 2.4]{KY2}, respectively.
\end{proof}

\subsection{Proof of  Theorems \ref{thm:Holllig type_Classical} and \ref{thm:PM type-Classical} }\label{sec:proof-classical}

We first begin with the proof of Theorem \ref{thm:Holllig type_Classical}.
Fix a number $\bar{s}\in (M_0,s_1)$. We then choose a function $\tilde\sigma\in C^3(\R)$ such that
\[
\left\{\begin{array}{l}
         \mbox{$\tilde\sigma(s)=\sigma(s)$ $\forall s\in(-\infty,\bar{s}]$,} \\
         \mbox{$\tilde\sigma'>0$ in $\R$}.
       \end{array}
 \right.
\]
Thanks to Theorem \ref{thm:uniform-parabolic}, there exists a unique global solution $u\in C^{2,1}(\bar\Omega_\infty)$ to the problem
\begin{equation*}
\left\{
\begin{array}{ll}
  u_t=(\tilde\sigma(u_x))_{x} & \mbox{in $\Omega_\infty$,} \\
  u=u_0 & \mbox{on $\Omega\times\{t=0\}$,} \\
  u_x=0 & \mbox{on $\partial\Omega\times(0,\infty)$,}
\end{array}
\right.
\end{equation*}
satisfying (1)--(3). Moreover, by the validity of (3) and the choice of $\tilde\sigma$, $u$ is also a global solution to (\ref{main-ibP}). Finally, (4) follows easily from (4) of Theorem \ref{thm:uniform-parabolic}

The proof of Theorem \ref{thm:PM type-Classical} is  similar; the only difference is that we fix two numbers $\bar{s}_1\in (s_1,m_0)$, $\bar{s}_2\in (M_0,s_2)$ and choose a function $\tilde\sigma\in C^3(\R)$ such that
\[
\left\{\begin{array}{l}
         \mbox{$\tilde\sigma(s)=\sigma(s)$ $\forall s\in[\bar{s}_1,\bar{s}_2]$,} \\
         \mbox{$\tilde\sigma'>0$ in $\R$}.
       \end{array}
 \right.
\]
Then the rest can be finished as in the proof of Theorem \ref{thm:Holllig type_Classical} above; we do not repeat the details.

\subsection{Baire's category theorem}
In this subsection, we prepare a suitable functional tool for solving a partial differential inclusion introduced in Subsection \ref{sec:phase-transition}. We begin with the definition of a related terminology.
\begin{defn}
Let $\mathfrak{X}$ and $\mathfrak{Y}$ be metric spaces. A map $f:\mathfrak{X}\to \mathfrak{Y}$ is called  \emph{Baire-one} if it is the pointwise limit of some sequence of continuous maps from $\mathfrak{X}$ into $\mathfrak{Y}$.
\end{defn}

A version of Baire's category theorem, regarding Baire-one maps, is the following.

\begin{thm}[Baire's category theorem]\label{thm:Baire-category}
Let $\mathfrak{X}$ and $\mathfrak{Y}$ be metric spaces with $\mathfrak{X}$ complete. If $f:\mathfrak{X}\to \mathfrak{Y}$ is a Baire-one map, then the set of points of discontinuity for $f$, say $\mathfrak{D}_f$, is of the first category. Therefore, the set of points in $\mathfrak{X}$ at which $f$ is continuous, that is, $\mathfrak{C}_f:=\mathfrak{X}\setminus\mathfrak{D}_f$, is dense in $\mathfrak{X}$.
\end{thm}

Our Baire-one maps to be used later are the gradient operators as follows.

\begin{pro}\label{prop:gradient-Baire-one}
Let $m$ and $n$ be any two positive integers. Let $U\subset\R^n$ be a bounded open set, and let $\mathfrak{X}\subset W^{1,\infty}(U;\R^m)$ be equipped with the $L^\infty(U;\R^m)$-metric. Then the gradient operator $\nabla:\mathfrak{X}\to L^p(U;\M^{m\times n})$ is a Baire-one map for each $p\in[1,\infty)$.
\end{pro}

Since the proof of Theorem \ref{thm:Baire-category} and Proposition \ref{prop:gradient-Baire-one} can be found in \cite[Chapter 10]{Da}, we do not reproduce it here.

\subsection{Rank-one smooth approximation}
Here, we equip with an important but technical tool for local patching to be used in the proof of the density lemma, Lemma \ref{lem:density lemma}, stated in the next subsection.  The following result  is a refinement of the $(1+1)$-dimensional version of a combination of \cite[Theorem 2.3 and Lemma 4.5]{KY1}, and its proof can be found in our recent paper \cite[Theorem 6.1]{KY3}.

\begin{lem}\label{lem:rank-1-smooth-approx}
Let $Q=(x_1,x_2)\times(t_1,t_2)\subset\R^{2}$ be an open rectangle, where $x_2>x_1$ and $t_2>t_1$ are fixed reals.
Given any $\lambda_1>0$, $\lambda_2>0$ and $\epsilon>0$, there exists a function $\omega=(\varphi,\psi)\in C^\infty_c(Q;\R^2)$ such that
\begin{itemize}
\item[(a)] $\|\omega\|_{L^\infty(Q)}<\epsilon$, $\|\varphi_t\|_{L^\infty(Q)}<\epsilon$, $\|\psi_t\|_{L^\infty(Q)}<\epsilon$,
\item[(b)] $-\lambda_1\le\varphi_x\le\lambda_2$ in $Q$,
\item[(c)] $\left\{\begin{array}{l}
              \big| |\{  (x,t)\in Q\, |\, \varphi_x(x,t) =-\lambda_1 \} |-\frac{\lambda_2}{\lambda_1+\lambda_2}|Q|\big| < \epsilon,\\
              \big| |\{  (x,t)\in Q\, |\, \varphi_x(x,t) =\lambda_2 \} |-\frac{\lambda_1}{\lambda_1+\lambda_2}|Q|\big| < \epsilon,
            \end{array}\right.
$
\item[(d)] $\psi_x=\varphi$ in $Q$, and
\item[(e)] $\int_{x_1}^{x_2}\varphi(x,t)\,dx=0$  for all $t_1<t<t_2$.
\end{itemize}
\end{lem}

\subsection{Phase transition as a differential inclusion}\label{sec:phase-transition}

In proving Theorems \ref{thm:Holllig type}, \ref{thm:PM type-Stable}, \ref{thm:PM type-Stable-Case2}, \ref{thm:PM type-Unstable} and \ref{thm:non-Fourier type}, we will keep repeatedly solving a certain partial differential inclusion in the framework of Baire's category method. To avoid an unnecessary repetition, we separate a germ of common analysis into this subsection and make use of it whenever it is needed.

\subsubsection{Initial definitions}

Let $r_1<r_2$ be any two fixed numbers, and let $g^\pm=g^\pm_{r_1,r_2}\in C([r_1,r_2])$ be functions such that
\[
g^-(r)<g^+(r)\;\;\forall r\in[r_1,r_2].
\]
Let $\tilde K^\pm=\tilde K_{g^\pm}$ denote the graph of  $g^\pm$; that is,
\[
\tilde K^\pm=\big\{(g^\pm(r),r)\in\R^2\,|\,r_1\le r\le r_2\big\}.
\]
Set $\tilde K=\tilde K_{g^+,g^-}=\tilde K^+\cup \tilde K^-$.
Let $\tilde{U}=\tilde{U}_{g^+,g^-}$ be the bounded domain in $\R^2$ given by
\[
\tilde{U}=\big\{(\lambda g^-(r)+(1-\lambda)g^+(r),r)\,|\,r_1< r< r_2, 0<\lambda<1\big\}.
\]
For each $u\in\R$, define the matrix sets
\[
\begin{split}
K(u)=K_{g^+,g^-}(u) & = \left\{\begin{pmatrix} s & c \\ u & r \end{pmatrix} \in \M^{2\times 2}\,\Big|\, (s,r)\in\tilde K, c\in\R \right\},\\
U(u)=U_{g^+,g^-}(u) & = \left\{\begin{pmatrix} s & c \\ u & r \end{pmatrix} \in \M^{2\times 2}\,\Big|\, (s,r)\in\tilde U, c\in\R \right\}.
\end{split}
\]

\subsubsection{Related differential inclusion}

We consider the \emph{inhomogeneous} partial differential inclusion
\begin{equation}\label{differential-inclusion-1}
\nabla w(x,t)\in K(u(x,t)),\;\;\mbox{a.e. $(x,t)\in Q$},
\end{equation}
where $Q\subset\R^2$ is a bounded open set, $w=(u,v):Q\to\R^2$ is a Lipschitz function, and $\nabla w=\begin{pmatrix} u_x & u_t \\ v_x & v_t\end{pmatrix}$ is the Jacobian matrix of $w$.

We say that a Lipschitz function $w:Q\to\R^2$ is a \emph{subsolution} of differential inclusion (\ref{differential-inclusion-1}) if
\begin{equation*}
\nabla w(x,t)\in K(u(x,t))\cup U(u(x,t)),\;\;\mbox{a.e. $(x,t)\in Q$},
\end{equation*}
and it is a \emph{strict subsolution} of (\ref{differential-inclusion-1})
if
\begin{equation*}
\nabla w(x,t)\in U(u(x,t)),\;\;\mbox{a.e. $(x,t)\in Q$}.
\end{equation*}

Let $w=(u,v):Q\to\R^2$ be a subsolution of differential inclusion (\ref{differential-inclusion-1}). For a.e. $(x,t)\in Q$, we can define the quantity
\[
Z_w(x,t)=Z^Q_{g^+,g^-;w}(x,t) =\frac{u_x(x,t)-g^-(v_t(x,t))}{g^+(v_t(x,t))-g^-(v_t(x,t))} \in[0,1].
\]
For each non-null subset $E$ of $Q$, we then define the \emph{transition gauge} of $w$ over $E$ by
\[
\Gamma_w^E=\Gamma^E_{g^+,g^-;w}=\frac{1}{|E|}\int_E Z_w(x,t)\,dxdt \in [0,1].
\]
This gauge measures the inclination of the diagonal $(u_x,v_t)$ of $\nabla w$ over $E$ towards the right branch $\tilde K^+$ of $\tilde K$; for instance, it is easy to see that
\[
\mbox{$\Gamma_w^E=1$ $\Longleftrightarrow$ $(u_x,v_t)\in\tilde K^+$ a.e. in $E$}
\]
and that
\[
\mbox{$\Gamma_w^E=0$ $\Longleftrightarrow$ $(u_x,v_t)\in\tilde K^-$ a.e. in $E$}.
\]
Note also that if $w$ is a strict subsolution of (\ref{differential-inclusion-1}), then $0<\Gamma_w^E<1$ for every non-null subset $E$ of $Q$.

\subsubsection{Generic setup and density lemma}

Let $b>0$ be any fixed number. For each $u\in\R$, define the matrix sets
\[
\begin{split}
K_b(u) & =K_{g^+,g^-;b}(u)= \left\{\begin{pmatrix} s & c \\ u & r \end{pmatrix} \in \M^{2\times 2}\,\Big|\, (s,r)\in\tilde K, |c|\le b \right\},\\
U_b(u) & =U_{g^+,g^-;b}(u)= \left\{\begin{pmatrix} s & c \\ u & r \end{pmatrix} \in \M^{2\times 2}\,\Big|\, (s,r)\in\tilde U, |c|< b \right\}.
\end{split}
\]
Let $Q\subset \R^2$ be a bounded open set, and let $w^*=(u^*,v^*)\in (C^{2,1}\times C^{3,1})(\bar{Q}):= C^{2,1}(\bar{Q})\times C^{3,1}(\bar{Q})$ be a Lipschitz function such that
\[
\nabla w^*\in U_b(u^*)\;\;\mbox{in $Q$};
\]
then $w^*$ is a strict subsolution of differential inclusion (\ref{differential-inclusion-1}) so that $0<\Gamma_{w^*}^E<1$ for each non-null subset $E$ of $Q$.

Let $\epsilon>0$ be any fixed number. We define a class $\mathcal{A}= \mathcal{A}_{g^+,g^-,b,w^*,Q,\epsilon}$ of strict subsolutions of (\ref{differential-inclusion-1}) by
\[
\mathcal{A}=\left\{ \begin{array}{l}
                      w=(u,v)\in \\
                      (C^{2,1}\times C^{3,1})(\bar{Q})
                    \end{array} \,\Bigg|\, \begin{array}{l}
                                             \mbox{$w$ is Lipschitz in $Q$, $w=w^*$ in $Q\setminus\bar Q^w$} \\
                                             \mbox{for some open set $Q^w\subset\subset Q$ with $|\partial Q^w|=0$,} \\
                                             \mbox{$\|u-u^*\|_{L^\infty(Q)} <\epsilon/2$, $\|u_t-u_t^*\|_{L^\infty(Q)} <\epsilon/2$,} \\
                                             \mbox{$\nabla w\in U_b(u)$ in $Q$, $|\Gamma^Q_w-\Gamma^Q_{w^*}|<\epsilon/2$}
                                           \end{array}
  \right\};
\]
then $w^*\in\mathcal{A}$ so that $\mathcal{A}\ne\emptyset$. Next, for each $\delta>0$, we define a subclass $\mathcal{A}_\delta=\mathcal{A}_{g^+,g^-,b,w^*,Q,\epsilon;\delta}$ of $\mathcal{A}$ by
\[
\mathcal{A}_\delta=\left\{ w\in\mathcal{A}\,\Big|\, \int_Q\mathrm{dist}\big(\nabla w(x,t),K_b(u(x,t))\big)\,dxdt\le\delta|Q| \right\}.
\]

We now state the pivotal density lemma that will be used throughout the rest of the paper.

\begin{lem}\label{lem:density lemma}
For each $\delta>0$, $\mathcal{A}_\delta$ is dense in $\mathcal{A}$ under the $L^\infty(Q)$-norm.
\end{lem}
The proof of this lemma is so long and complicated that we postpone it to the last section, Section \ref{sec:density-proof}.

\section{Eventual smoothing for H\"ollig type equations}\label{sec:proof of Hollig type}

This section is devoted to proving Theorem \ref{thm:Holllig type}, that is, the existence of eventually smooth global weak solutions $u$ to problem (\ref{main-ibP}) for all initial data $u_0$ with $M_0>\bar s_1$ when the diffusion flux $\sigma$ is of the H\"ollig type.

Following the notation and setup in Subsection \ref{sec:phase-transition}, we divide the proof into several steps.

\subsection*{Setup and subsolution}

\begin{figure}[ht]
\begin{center}
\includegraphics[scale=0.6]{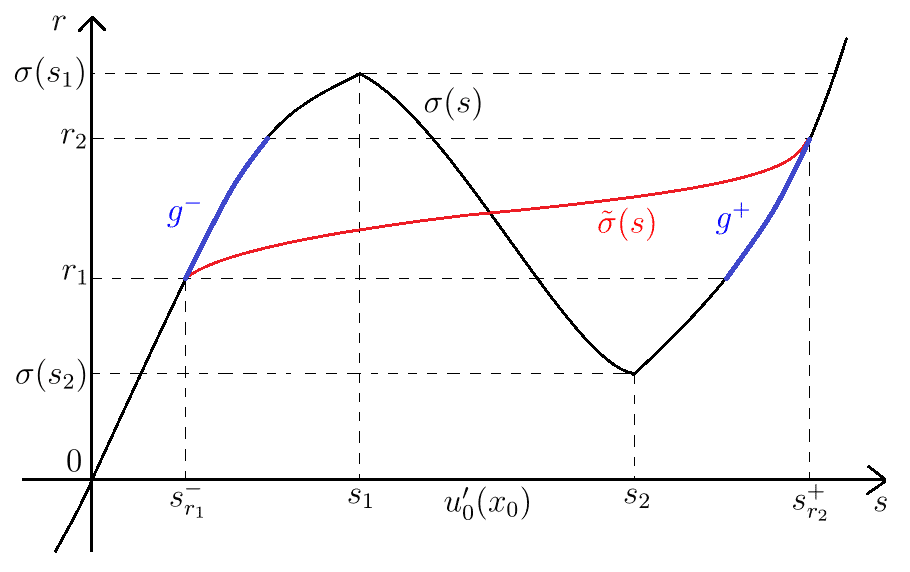}
\end{center}
\caption{Modified flux $\tilde\sigma(s)$}
\label{fig4}
\end{figure}

Let
\[
g^\pm=g^\pm_{r_1,r_2}=(\sigma|_{[s^\pm_{r_1},s^\pm_{r_2}]})^{-1}:[r_1,r_2]\to [s^\pm_{r_1},s^\pm_{r_2}];
\]
then $g^-\le g^-(r_2)=s^-_{r_2}<s^+_{r_1}=g^+(r_1)\le g^+$ on $[r_1,r_2]$. 
Next, we choose a function $\tilde\sigma=\tilde\sigma_{r_1,r_2}\in C^3(\R)$ (see Figure \ref{fig4}) such that
\begin{equation}\label{proof-thm-2.2-1}
\left\{ \begin{array}{l}
          \tilde\sigma=\sigma\;\;\mbox{on $(-\infty,s^-_{r_1}]\cup[s^+_{r_2},\infty)$}, \\
          \tilde\sigma<\sigma\;\;\mbox{on $(s^-_{r_1},s^-_{r_2}]$},\; \tilde\sigma>\sigma\;\;\mbox{on $[s^+_{r_1},s^+_{r_2})$}, \\
          \tilde\sigma'>0\;\;\mbox{in $\R$.}
        \end{array}
 \right.
\end{equation}
Thanks to Theorem \ref{thm:uniform-parabolic}, there exists a unique global solution $u^*\in C^{2,1}(\bar\Omega_\infty)$ to problem (\ref{main-ibP}), with $\sigma$ replaced by $\tilde\sigma$, satisfying (6)--(9), where the constants $\kappa$, $\lambda$ and $m$ are fixed as in (4) of Theorem \ref{thm:uniform-parabolic} to yield (9).

Now, let $Q_1$, $Q_2$ and $F^\pm$ be defined as in the statement of the theorem. It then follows from (9) that $Q_1\cup Q_2\cup F^+\cup F^-$ is bounded. Also, this set has a positive distance to the vertical boundary $\partial\Omega\times[0,\infty)$ since $u^*_x$ vanishes on $\partial\Omega\times[0,\infty)$ and is continuous on $\bar\Omega_\infty;$ thus (4) holds. Since $s^-_{r_1}<u_0'(x_0)<s^+_{r_2}$, we have $Q_1\ne \emptyset$ so that $0<\sup_{(x,t)\in Q_1}t=:t_1<\infty$ by (4). This together with (8) and (9) implies (5).

Let $Q=Q_1$ and $b=\|u^*_t\|_{L^\infty(Q)}+1.$ With the functions $g^\pm$ defined above, for each  $u\in\R$, let $K_b(u)=K_{g^+,g^-;b}(u)$ and $U_b(u)=U_{g^+,g^-;b}(u)$ be given as in Subsection \ref{sec:phase-transition}. Define an auxiliary function $v^*$ by
\[
v^*(x,t) = \int_0^t \tilde\sigma(u_x^*(x,\tau))\,d\tau + \int_0^x u_0(y)\,dy\;\;\forall (x,t)\in\Omega_\infty;
\]
then
\begin{equation}\label{proof-thm-2.2-6}
\left\{ \begin{array}{l}
          v^*_t=\tilde\sigma(u^*_x) \\
          v^*_x=u^*
        \end{array}
 \right. \mbox{in $\Omega_\infty$},
\end{equation}
and so $w^*:=(u^*,v^*)\in (C^{2,1}\times C^{3,1})(\bar\Omega_\infty)$.
From this, (\ref{proof-thm-2.2-1}) and the definition of $b$, we have
\[
\nabla w^*\in U_b(u^*)\;\;\mbox{in $Q$};
\]
thus $w^*|_{Q}$ becomes a strict subsolution of differential inclusion (\ref{differential-inclusion-1}) so that $0<\Gamma^Q_{w^*}<1$, where the gauge operator $\Gamma^Q_\cdot=\Gamma^Q_{g^+,g^-;\cdot}$ is given as in Subsection \ref{sec:phase-transition}.

\subsection*{Admissible class}

Fix any $\epsilon>0$. Following the above, let $\mathcal{A}=\mathcal{A}_{g^+,g^-,b,w^*,Q,\epsilon}$ and $\mathcal{A}_\delta=\mathcal{A}_{g_+,g_-,b,w^*,Q,\epsilon;\delta}$ $(\delta>0)$ be defined as in Subsection \ref{sec:phase-transition}. Set $t_2=t_1+1$, where $t_1>0$ is as above; then $Q\subset\Omega_{t_2-1}\subset\Omega_{t_2}$. We now define a class $\mathfrak{A}$ of \emph{admissible functions} by
\[
\mathfrak{A}= \left\{ \begin{array}{l}
                        w=(u,v)\in \\
                        (C^{2,1}\times C^{3,1})(\bar\Omega_{t_2}) \end{array} \,\Big|\, \begin{array}{l}
                        \mbox{$w=w^*$ in $\Omega_{t_2}\setminus Q$,} \\
                        w|_{Q}\in\mathcal{A}

                      \end{array}
 \right\};
\]
then $w^*\in\mathfrak{A}$ so that $\mathfrak{A}\neq\emptyset$. Next, for each $\delta>0$, let $\mathfrak{A}_\delta$ be the subclass of functions $w\in\mathfrak{A}$ such that $w|_Q\in\mathcal{A}_\delta$.

\subsection*{Baire's category method}

Let $\mathfrak{X}$ denote the closure of $\mathfrak{A}$ in the space $L^\infty(\Omega_{t_2};\R^2)$; then $(\mathfrak{X},L^\infty)$ becomes a nonempty complete metric space. From the definition of $\mathfrak{A}$, it is easily checked that
\[
\mathfrak{X}\subset w^*+W^{1,\infty}_0(\Omega_{t_2};\R^2).
\]

By Proposition \ref{prop:gradient-Baire-one}, the space-time gradient operator $\nabla:\mathfrak{X}\to L^1(\Omega_{t_2};\M^{2\times 2})$ is a Baire-one map, and thus it follows from the Baire category theorem, Theorem \ref{thm:Baire-category}, that the set of points at which $\nabla$ is continuous, say $\mathfrak{C}_\nabla$, is dense in $\mathfrak{X}$.

\subsection*{Solutions from $\mathfrak{C}_\nabla$}

Let $w=(u,v)\in\mathfrak{C}_\nabla,$ and extend it from $\Omega_{t_2}$ to $\Omega_\infty$ by setting $w= w^*$ in $\Omega_\infty\setminus\Omega_{t_2}$. We now show that $u$ is a global weak solution to problem (\ref{main-ibP}) satisfying (1)--(3).

Note first that from Lemma \ref{lem:density lemma} and the density of $\mathfrak{A}$ in $\mathfrak{X}$, $\mathfrak{A}_{1/j}$ is dense in $\mathfrak{X}$ for each $j\in\N$. So we can take a sequence $w_j=(u_j,v_j)\in\mathfrak{A}_{1/j}$ such that $\|w_j-w\|_{L^\infty(\Omega_{t_2})}\to 0$ as $j\to\infty$. By the continuity of the operator $\nabla$ at $w$, we have $\nabla w_j\to\nabla w$ in $L^1(\Omega_{t_2};\M^{2\times 2})$ and hence pointwise a.e. in $\Omega_{t_2}$ after passing to a subsequence if necessary.

We state some properties of the function $w_j=(u_j,v_j)$ inherited from its membership $w_j\in\mathfrak{A}_{1/j}$:
\begin{equation}\label{proof-thm-2.2-3}
\left\{\begin{array}{l}
         w_j=(u_j,v_j)\in (C^{2,1}\times C^{3,1})(\bar\Omega_{t_2}), \\
         \mbox{$w_j=w^*$ in $\Omega_{t_2}\setminus\bar{Q}^{w_j}$ for some open set $Q^{w_j}\subset\subset Q$ with $|\partial Q^{w_j}|=0$}, \\
         \|u_j-u^*\|_{L^\infty(Q)}<\epsilon/2,\; \|(u_j)_t-u_t^*\|_{L^\infty(Q)}<\epsilon/2,\; |\Gamma^Q_{w_j}-\Gamma^Q_{w^*}|<\epsilon/2, \\
         \mbox{$\nabla w_j\in U_b(u_j)$ in $Q$, and} \\
         \int_Q\mathrm{dist}\big(\nabla w_j(x,t),K_b(u_j(x,t))\big) \,dxdt\le\frac{1}{j}|Q|.
       \end{array}
 \right.
\end{equation}

We now check that $u$ satisfies (1)--(3). Letting $j\to\infty$ in (\ref{proof-thm-2.2-3}), we get
\[
\left\{\begin{array}{l}
         \mbox{$w=w^*$ in $\bar\Omega_{t_2}\setminus Q$,} \\
         \|u-u^*\|_{L^\infty(Q)}\le\epsilon/2,\; \|u_t-u_t^*\|_{L^\infty(Q)}\le\epsilon/2;
       \end{array}
 \right.
\]
thus (1) holds. We also have
\[
\begin{split}
\int_Q\mathrm{dist}\big(\nabla w_j & (x,t),K_b(u_j(x,t))\big) \,dxdt \\
& \to \int_Q\mathrm{dist}\big(\nabla w(x,t),K_b(u(x,t))\big) \,dxdt=0
\end{split}
\]
so that $\nabla w\in K_b(u)$ a.e. in $Q$. This inclusion implies that a.e. in $Q$,
\begin{equation}\label{proof-thm-2.2-2}
\left\{\begin{array}{l}
         v_t=\sigma(u_x), \\
         u_x\in[s^-_{r_1},s^-_{r_2}]\cup [s^+_{r_1},s^+_{r_2}].
       \end{array}
 \right.
\end{equation}
From (\ref{proof-thm-2.2-3}), we have $(u_j)_x=s^+_{r_2}$ in $F^+$ and $(u_j)_x=s^-_{r_1}$ in $F^-$. Letting $j\to\infty$, we have $u_x=s^+_{r_2}$ a.e. in $F^+$ and $u_x=s^-_{r_1}$ a.e. in $F^-$; this and (\ref{proof-thm-2.2-2}) yield (2).
Also, letting $j\to\infty$ in (\ref{proof-thm-2.2-3}), we have $|\Gamma^Q_{w}-\Gamma^Q_{w^*}|\le\epsilon/2<\epsilon$. Without loss of generality, we can assume $0<\epsilon<\min\{\Gamma^Q_{w^*},1-\Gamma^Q_{w^*}\}$ so that $0<\Gamma^Q_{w}<1$. This and (\ref{proof-thm-2.2-2}) thus imply (3).

For a later use, note from (\ref{proof-thm-2.2-3}) that  $(v_j)_t=v^*_t$ and $(u_j)_x=u^*_x\not\in(s^-_{r_1},s^+_{r_2})$ in $\Omega_{t_2}\setminus Q$, and so by (\ref{proof-thm-2.2-1}) and (\ref{proof-thm-2.2-6}),
\begin{equation}\label{proof-thm-2.2-4}
v_t=v^*_t=\tilde\sigma(u^*_x)=\sigma(u^*_x)=\sigma(u_x)\;\;\mbox{a.e. in $\Omega_{t_2}\setminus Q$}.
\end{equation}
By (\ref{proof-thm-2.2-6}) and (\ref{proof-thm-2.2-3}), we get $(v_j)_x=u_j$ in $\Omega_{t_2}$; thus we  have
\begin{equation}\label{proof-thm-2.2-5}
v_x=u\;\;\mbox{a.e. in $\Omega_{t_2}$}.
\end{equation}

Let us now check that $u$ is a global weak solution to problem (\ref{main-ibP}). To do this, let $\zeta\in C^\infty(\bar\Omega_\infty)$ and $\tau\in[0,\infty)$. Assume first that $\tau\le t_2$. From (\ref{proof-thm-2.2-6}) and (\ref{proof-thm-2.2-3}), we then have
\[
\begin{split}
& \int_0^\tau\int_0^L u_j\zeta_t\,dxdt = \int_0^\tau\int_0^L (v_j)_x\zeta_t\,dxdt \\
& = -\int_0^\tau\int_0^L (v_j)_{xt}\zeta\,dxdt + \int_0^L \big( (v_j)_x(x,\tau)\zeta(x,\tau)-(v_j)_x(x,0)\zeta(x,0) \big)\,dx \\
& = \int_0^\tau\int_0^L (v_j)_{t}\zeta_x\,dxdt - \int_0^\tau \big( (v_j)_t(L,t)\zeta(L,t)-(v_j)_t(0,t)\zeta(0,t) \big)\,dt \\
& \;\;\;\; + \int_0^L \big( (v_j)_x(x,\tau)\zeta(x,\tau)-(v_j)_x(x,0)\zeta(x,0) \big)\,dx \\
& = \int_0^\tau\int_0^L (v_j)_{t}\zeta_x\,dxdt - \int_0^\tau \big( v^*_t(L,t)\zeta(L,t)-v^*_t(0,t)\zeta(0,t) \big)\,dt \\
& \;\;\;\; + \int_0^L \big( u_j(x,\tau)\zeta(x,\tau)-u_0(x)\zeta(x,0) \big)\,dx.
\end{split}
\]
Let $j\to\infty$; then from (\ref{proof-thm-2.2-2}), (\ref{proof-thm-2.2-4}) and $v^*_t=\tilde\sigma(u^*_x)=0$ on $\partial\Omega\times[0,\infty)$, we have
\[
\begin{split}
\int_0^\tau\int_0^L u\zeta_t & \,dxdt = \int_0^\tau\int_0^L  v_{t}\zeta_x\,dxdt + \int_0^L \big( u(x,\tau)\zeta(x,\tau)-u_0(x)\zeta(x,0) \big)\,dx \\
&= \int_0^\tau\int_0^L  \sigma(u_x)\zeta_x\,dxdt + \int_0^L \big( u(x,\tau)\zeta(x,\tau)-u_0(x)\zeta(x,0) \big)\,dx;
\end{split}
\]
hence (\ref{def:weak-solution}) holds.

Next, assume that $\tau>t_2.$ Since $u=u^*$ and $u^*_x\le s^-_{r_1}$ on $\bar\Omega\times[t_2,\infty)$ and $u^*$ is a global classical solution to problem (\ref{main-ibP}) with $\sigma$ replaced by $\tilde\sigma$, it follows from (\ref{proof-thm-2.2-1}) that
\[
\begin{split}
\int_0^L\big( & u(x,\tau)\zeta(x,\tau)- u(x,t_2)\zeta(x,t_2)\big)\,dx \\
& = \int_0^L\big(u^*(x,\tau)\zeta(x,\tau)- u^*(x,t_2)\zeta(x,t_2)\big)\,dx \\
& = \int_{t_2}^\tau\int_0^L \big(u^*\zeta_t-\tilde\sigma(u^*_x)\zeta_x \big)\,dxdt  = \int_{t_2}^\tau\int_0^L \big(u\zeta_t-\sigma(u_x)\zeta_x \big)\,dxdt.
\end{split}
\]
Combining this with the previous case with $\tau=t_2$, we obtain (\ref{def:weak-solution}); therefore, $u$ is a global weak solution to (\ref{main-ibP}).

\subsection*{Infinitely many solutions} We have shown above that the first component $u$ of each function $w=(u,v)\in \mathfrak{C}_\nabla$ with common extension $u=u^*$ in $\Omega_\infty\setminus\Omega_{t_2}$ is a global weak solution to problem (\ref{main-ibP}) satisfying (1)--(3). It only remains to verify that $\mathfrak{C}_\nabla$ contains infinitely many functions and that no two different functions in $\mathfrak{C}_\nabla$ have the first component that are identical. Suppose first on the contrary that $\mathfrak{C}_\nabla$ has only finitely many functions. Then
\[
\mathfrak{C}_\nabla=\bar{\mathfrak{C}}_\nabla=\mathfrak{X} =\bar{\mathfrak{A}}=\mathfrak{A}\ni w^*=(u^*,v^*),
\]
where the closure is with respect to the $L^\infty(\Omega_{t_2};\R^2)$-norm. By the above result, we now have that $u^*$ itself is a global weak solution to (\ref{main-ibP}) satisfying (1)--(3); clearly, this is a contradiction. Thus, $\mathfrak{C}_\nabla$ contains infinitely many functions.
Next, let $w_1=(u_1,v_1),w_2=(u_2,v_2)\in\mathfrak{C}_\nabla.$ It is sufficient to show that
\[
u_1=u_2\;\;\mbox{in $\Omega_{t_2}$}\;\;\Longleftrightarrow\;\;v_1=v_2\;\;\mbox{in $\Omega_{t_2}$}.
\]
If $u_1=u_2$ in $\Omega_{t_2}$, then by (\ref{proof-thm-2.2-5}), we have $(v_1)_x=u_1=u_2=(v_2)_x$ a.e. in $\Omega_{t_2}$. Since $v_1=v_2=v^*$ on $\partial\Omega_{t_2}$, we get $v_1=v_2$ in $\Omega_{t_2}.$ The converse is also easy to check; we skip this.

The proof of Theorem \ref{thm:Holllig type} is now complete.

\section{Eventual smoothing for Perona-Malik type equations}\label{sec:proof-PM-eventual-smooth}

In this section, we prove the existence of eventually smooth global weak solutions $u$ to problem (\ref{main-ibP}) for all nonconstant initial data $u_0$ when the diffusion flux $\sigma$ is of the Perona-Malik type. We first provide the proof of Theorem \ref{thm:PM type-Stable}, which handles Case 1: $m_0<0<M_0.$
Since the proof of Theorem \ref{thm:PM type-Stable-Case2} for Case 2: $0=m_0<M_0$ is very similar to and even simpler than that of Theorem \ref{thm:PM type-Stable}, we only give a sketch for the proof at the end of the section. Likewise, for Case 3: $m_0<M_0=0$ that we do not state as a theorem, one can easily adapt the proof of Theorem \ref{thm:PM type-Stable} to deduce the expected result.

\begin{proof}[Proof of Theorem \ref{thm:PM type-Stable}]
Following the notation and setup in Subsection \ref{sec:phase-transition}, we divide the proof into several steps.

\subsection*{Setup and subsolution}

Let
\[
\begin{split}
g_1^+ & =g^+_{r_1',r_1}=(\sigma|_{[s^+_{r_1'},s^+_{r_1}]})^{-1} :[r_1',r_1]\to[s^+_{r_1'},s^+_{r_1}], \\
g_1^- & =g^-_{r_1',r_1}=(\sigma|_{[s^-_{r_1},s^-_{r_1'}]})^{-1} :[r_1',r_1]\to[s^-_{r_1},s^-_{r_1'}], \\
g_2^+ & =g^+_{r_2,r_2'}=(\sigma|_{[s^+_{r_2'},s^+_{r_2}]})^{-1} :[r_2,r_2']\to[s^+_{r_2'},s^+_{r_2}], \\
g_2^- & =g^-_{r_2,r_2'}=(\sigma|_{[s^-_{r_2},s^-_{r_2'}]})^{-1} :[r_2,r_2']\to[s^-_{r_2},s^-_{r_2'}];
\end{split}
\]
then
\[
\begin{split}
g_1^-\le g_1^-(r_1')=s^-_{r_1'}&<s^+_{r_1'}=g_1^+(r_1')\le g_1^+\;\; \mbox{on $[r_1',r_1]$}, \\
g_2^-\le g_2^-(r_2')=s^-_{r_2'}&<s^+_{r_2'}=g_2^+(r_2')\le g_2^+\;\; \mbox{on $[r_2,r_2']$}.
\end{split}
\]
Next, choose a function $\tilde\sigma=\tilde\sigma_{r_1,r_1',r_2,r_2'}\in C^3(\R)$ (see Figure \ref{fig5}) such that
\begin{equation}\label{proof-thm-2.3-1}
\left\{\begin{array}{l}
         \mbox{$\tilde\sigma=\sigma$ on $[s^+_{r_1},s^-_{r_2}]$,} \\
         \mbox{$\tilde\sigma<\min\{\sigma,r_2'\}$ on $(s^-_{r_2},M_0]$,} \\
         \mbox{$\tilde\sigma>\max\{\sigma,r_1'\}$ on $[m_0,s^+_{r_1})$,} \\
         \mbox{$\tilde\sigma'>0$ in $\R$.}
       \end{array}
 \right.
\end{equation}
Owing to Theorem \ref{thm:uniform-parabolic}, there exists a unique global solution $u^*\in C^{2,1}(\bar\Omega_\infty)$ to problem (\ref{main-ibP}), with $\sigma$ replaced by $\tilde\sigma$, satisfying (5)--(8), where the constants $\kappa$, $\lambda$ and $m$ are fixed as in (4) of Theorem \ref{thm:uniform-parabolic} to imply (8).

\begin{figure}[ht]
\begin{center}
\includegraphics[scale=0.6]{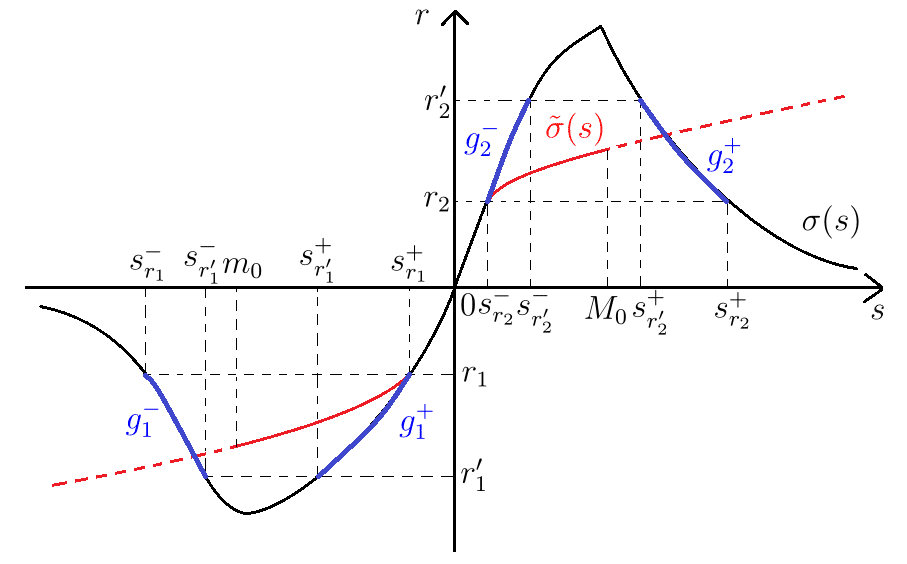}
\end{center}
\caption{Modified flux $\tilde\sigma(s)$}
\label{fig5}
\end{figure}

Now, let $Q_1$, $Q_2$, $F_1$ and $F_2$ be defined as in the statement of the theorem. From (8), it follows that $Q_1\cup Q_2\cup F_1\cup F_2$ is bounded. This set also has a positive distance to the vertical boundary $\partial\Omega\times[0,\infty)$, since $u_x^*$ vanishes on $\partial\Omega\times[0,\infty)$ and is continuous on $\bar\Omega_\infty$; thus (4) is valid.

Let $b=\|u^*_t\|_{L^\infty(Q_1\cup Q_2)}+1$. With the functions $g_i^\pm$ $(i=1,2)$ defined above, for each $u\in\R,$ let $K_{i,b}(u)=K_{g_i^+,g_i^-;b}(u)$ and $U_{i,b}(u)=U_{g_i^+,g_i^-;b}(u)$ be given as in Subsection \ref{sec:phase-transition}. Define an auxiliary function $v^*$ by
\[
v^*(x,t)=\int_0^t\tilde\sigma(u^*_x(x,\tau))\,d\tau + \int_0^x u_0(y)\,dy\;\;\forall (x,t)\in\Omega_\infty;
\]
then
\begin{equation}\label{proof-thm-2.3-6}
\left\{\begin{array}{l}
         v^*_t=\tilde\sigma(u^*_x) \\
         v^*_x=u^*
       \end{array}\mbox{in $\Omega_\infty$,}
 \right.
\end{equation}
so that $w^*:=(u^*,v^*)\in (C^{2,1}\times C^{3,1})(\bar\Omega_\infty)$. From this, (\ref{proof-thm-2.3-1}) and the definition of $b$, it follows that for $i=1,2$,
\[
\nabla w^*\in U_{i,b}(u^*)\;\;\mbox{in $Q_i$};
\]
thus $w^*|_{Q_i}$ is a strict subsolution of differential inclusion (\ref{differential-inclusion-1}) with $Q=Q_i$, and so $0<\Gamma^{Q_i}_{w^*}<1$, where the gauge operator $\Gamma^{Q_i}_\cdot=\Gamma^{Q_i}_{g_i^+,g_i^-;\cdot}$ is given as in Subsection \ref{sec:phase-transition}.

\subsection*{Admissible class}

Fix any $\epsilon>0$. Following the above notation, for $i=1,2$, let $\mathcal{A}_i=\mathcal{A}_{g_i^+,g_i^-,b,w^*,Q_i,\epsilon}$ and $\mathcal{A}_{i,\delta}=\mathcal{A}_{g_i^+,g_i^-,b,w^*,Q_i,\epsilon;\delta}$ $(\delta>0)$ be defined as in Subsection \ref{sec:phase-transition}. Set $t_1=\sup_{(x,t)\in Q_1\cup Q_2}t +1$; then $Q_1\cup Q_2\subset \Omega_{t_1-1}\subset\Omega_{t_1}.$ We now define a class $\mathfrak{A}$ of \emph{admissible functions} by
\[
\mathfrak{A}= \left\{\begin{array}{l}
                       w=(u,v)\in \\
                       (C^{2,1}\times C^{3,1})(\bar\Omega_{t_1})
                     \end{array}\,\Big|\, \begin{array}{l}
                                            \mbox{$w=w^*$ in $\Omega_{t_1}\setminus(Q_1\cup Q_2)$,} \\
                                            \mbox{$w|_{Q_i}\in\mathcal{A}_i$ for $i=1,2$}
                                          \end{array}
 \right\};
\]
then $w^*\in\mathfrak{A}$, and so $\mathfrak{A}\ne\emptyset.$ Also, for each $\delta>0$, let $\mathfrak{A}_\delta$ be the subclass of functions $w\in\mathfrak{A}$ such that $w|_{Q_i}\in\mathcal{A}_{i,\delta}$ for $i=1,2$.

\subsection*{Baire's category method}

Let $\mathfrak{X}$ denote the closure of $\mathfrak{A}$ in the space $L^\infty(\Omega_{t_1};\R^2)$; then $(\mathfrak{X},L^\infty)$ becomes a nonempty complete metric space. From the definition of $\mathfrak{A}$, it is easily checked that
\[
\mathfrak{X}\subset w^*+W^{1,\infty}_0(\Omega_{t_1};\R^2).
\]
As in the proof of Theorem \ref{thm:Holllig type}, we know that the set of points of continuity for the space-time gradient operator $\nabla:\mathfrak{X}\to L^1(\Omega_{t_1};\M^{2\times 2})$, say $\mathfrak{C}_\nabla,$ is dense in $\mathfrak{X}$.

\subsection*{Solutions from $\mathfrak{C}_\nabla$}

Let $w=(u,v)\in\mathfrak{C}_\nabla,$ and extend it from $\Omega_{t_1}$ to $\Omega_\infty$ by letting $w= w^*$ in $\Omega_\infty\setminus\Omega_{t_1}$. We now show that $u$ is a global weak solution to problem (\ref{main-ibP}) satisfying (1)--(3).

Note first that by Lemma \ref{lem:density lemma} and the density of $\mathfrak{A}$ in $\mathfrak{X}$, $\mathfrak{A}_{1/j}$ is dense in $\mathfrak{X}$ for each $j\in\N$. So we can take a sequence $w_j=(u_j,v_j)\in\mathfrak{A}_{1/j}$ such that $\|w_j-w\|_{L^\infty(\Omega_{t_1})}\to 0$ as $j\to\infty$. By the continuity of the operator $\nabla$ at $w$, we have $\nabla w_j\to\nabla w$ in $L^1(\Omega_{t_1};\M^{2\times 2})$ and hence pointwise a.e. in $\Omega_{t_1}$ after passing to a subsequence if necessary.

We now state some properties of the function $w_j=(u_j,v_j)$ inherited from its membership $w_j\in\mathfrak{A}_{1/j}$:
\begin{equation}\label{proof-thm-2.3-2}
\left\{\begin{array}{l}
         w_j=(u_j,v_j)\in (C^{2,1}\times C^{3,1})(\bar\Omega_{t_1}), \\
         \mbox{$w_j=w^*$ in $\Omega_{t_1}\setminus(\bar{Q}_1^{w_j}\cup \bar{Q}_2^{w_j})$ for some open sets $Q_1^{w_j}\subset\subset Q_1$} \\
         \mbox{\;\;\;\;\;and $Q_2^{w_j}\subset\subset Q_2$ with $|\partial Q_1^{w_j}|=|\partial Q_2^{w_j}|=0$}, \\
         \|u_j-u^*\|_{L^\infty(Q_i)}<\epsilon/2,\; \|(u_j)_t-u_t^*\|_{L^\infty(Q_i)}<\epsilon/2, \\
         \mbox{\;\;\;\, $|\Gamma^{Q_i}_{w_j}-\Gamma^{Q_i}_{w^*}|<\epsilon/2$ for $i=1,2$}, \\
         \mbox{$\nabla w_j\in U_{i,b}(u_j)$ in $Q_i$ for $i=1,2$, and} \\
         \mbox{$\int_{Q_i}\mathrm{dist}\big(\nabla w_j(x,t),K_{i,b}(u_j(x,t))\big) \,dxdt\le\frac{1}{j}|Q_i|$ for $i=1,2$}.
       \end{array}
 \right.
\end{equation}

Let us now check that $u$ satisfies (1)--(3). Letting $j\to\infty$ in (\ref{proof-thm-2.3-2}), we get
\[
\left\{\begin{array}{l}
         \mbox{$w=w^*$ in $\bar\Omega_{t_1}\setminus(Q_1\cup Q_2)$,} \\
         \|u-u^*\|_{L^\infty(Q_1\cup Q_2)}\le\epsilon/2,\; \|u_t-u_t^*\|_{L^\infty(Q_1\cup Q_2)}\le\epsilon/2;
       \end{array}
 \right.
\]
thus (1) holds. For $i=1,2$, as $j\to\infty,$ we have
\[
\begin{split}
\int_{Q_i}\mathrm{dist}\big(\nabla w_j & (x,t),K_{i,b}(u_j(x,t))\big) \,dxdt \\
& \to \int_{Q_i}\mathrm{dist}\big(\nabla w(x,t),K_{i,b}(u(x,t))\big) \,dxdt=0
\end{split}
\]
so that $\nabla w\in K_{i,b}(u)$ a.e. in $Q_i$. These inclusions imply that for $i=1,2$, we have, a.e. in $Q_i$,
\begin{equation}\label{proof-thm-2.3-3}
\left\{\begin{array}{l}
         v_t=\sigma(u_x), \\
         u_x\in[s^-_{r_i},s^-_{r_i'}]\cup [s^+_{r_i'},s^+_{r_i}].
       \end{array}
 \right.
\end{equation}
From (\ref{proof-thm-2.3-2}), we have $(u_j)_x=s^+_{r_1}$ in $F_1$ and $(u_j)_x=s^-_{r_2}$ in $F_2$. Letting $j\to\infty$, we have $u_x=s^+_{r_1}$ a.e. in $F_1$ and $u_x=s^-_{r_2}$ a.e. in $F_2$; this together with (\ref{proof-thm-2.3-3}) yields (2).
Also, letting $j\to\infty$ in (\ref{proof-thm-2.3-2}), we have $|\Gamma^{Q_i}_{w}-\Gamma^{Q_i}_{w^*}|\le\epsilon/2<\epsilon$ for $i=1,2$. Without loss of generality, we can assume $0<\epsilon<\min\{\Gamma^{Q_1}_{w^*},1-\Gamma^{Q_1}_{w^*}, \Gamma^{Q_2}_{w^*},1-\Gamma^{Q_2}_{w^*}\}$ so that $0<\Gamma^{Q_i}_{w}<1$ $(i=1,2)$. This and (\ref{proof-thm-2.3-3}) thus imply (3).

For a later use, note from (\ref{proof-thm-2.3-2}) and (7) that  $(v_j)_t=v^*_t$ and $(u_j)_x=u^*_x\in[s^+_{r_1},s^-_{r_2}]$ in $\Omega_{t_1}\setminus (Q_1\cup Q_2)$, and so by (\ref{proof-thm-2.3-1}) and (\ref{proof-thm-2.3-6}),
\begin{equation}\label{proof-thm-2.3-4}
v_t=v^*_t=\tilde\sigma(u^*_x)=\sigma(u^*_x)=\sigma(u_x)\;\;\mbox{a.e. in $\Omega_{t_1}\setminus (Q_1\cup Q_2)$}.
\end{equation}
By (\ref{proof-thm-2.3-6}) and (\ref{proof-thm-2.3-2}), we have $(v_j)_x=u_j$ in $\Omega_{t_1}$ so that
\begin{equation}\label{proof-thm-2.3-5}
v_x=u\;\;\mbox{a.e. in $\Omega_{t_1}$}.
\end{equation}

Let us now check that $u$ is a global weak solution to problem (\ref{main-ibP}). To do this, let $\zeta\in C^\infty(\bar\Omega_\infty)$ and $\tau\in[0,\infty)$. Assume first that $\tau\le t_1$. From (\ref{proof-thm-2.3-6}) and (\ref{proof-thm-2.3-2}), as in the proof of Theorem \ref{thm:Holllig type}, we have
\[
\begin{split}
\int_0^\tau\int_0^L & u_j\zeta_t\,dxdt  \\
& = \int_0^\tau\int_0^L (v_j)_{t}\zeta_x\,dxdt - \int_0^\tau \big( v^*_t(L,t)\zeta(L,t)-v^*_t(0,t)\zeta(0,t) \big)\,dt \\
& \;\;\;\; + \int_0^L \big( u_j(x,\tau)\zeta(x,\tau)-u_0(x)\zeta(x,0) \big)\,dx.
\end{split}
\]
Let $j\to\infty$; then from (\ref{proof-thm-2.3-3}), (\ref{proof-thm-2.3-4}) and $v^*_t=\tilde\sigma(u^*_x)=0$ on $\partial\Omega\times[0,\infty)$, we have
\[
\begin{split}
\int_0^\tau\int_0^L u\zeta_t & \,dxdt = \int_0^\tau\int_0^L  v_{t}\zeta_x\,dxdt + \int_0^L \big( u(x,\tau)\zeta(x,\tau)-u_0(x)\zeta(x,0) \big)\,dx \\
&= \int_0^\tau\int_0^L  \sigma(u_x)\zeta_x\,dxdt + \int_0^L \big( u(x,\tau)\zeta(x,\tau)-u_0(x)\zeta(x,0) \big)\,dx;
\end{split}
\]
hence (\ref{def:weak-solution}) holds.

Next, assume that $\tau>t_1.$ Since $u=u^*$ and $u^*_x\in [s^+_{r_1},s^-_{r_2}]$ on $\bar\Omega\times[t_1,\infty)$, it follows from (\ref{proof-thm-2.3-1}) that
\[
\begin{split}
\int_0^L\big( & u(x,\tau)\zeta(x,\tau)- u(x,t_1)\zeta(x,t_1)\big)\,dx \\
& = \int_0^L\big(u^*(x,\tau)\zeta(x,\tau)- u^*(x,t_1)\zeta(x,t_1)\big)\,dx \\
& = \int_{t_1}^\tau\int_0^L \big(u^*\zeta_t-\tilde\sigma(u^*_x)\zeta_x \big)\,dxdt  = \int_{t_1}^\tau\int_0^L \big(u\zeta_t-\sigma(u_x)\zeta_x \big)\,dxdt.
\end{split}
\]
Combining this with the previous case with $\tau=t_1$, we obtain (\ref{def:weak-solution}); therefore, $u$ is a global weak solution to (\ref{main-ibP}).

\subsection*{Infinitely many solutions} This part is identical to the corresponding one in the proof of Theorem \ref{thm:Holllig type}; we thus skip this.

Theorem \ref{thm:PM type-Stable} is now proved.
\end{proof}

To prove Theorem \ref{thm:PM type-Stable-Case2}, one can follow the above proof very closely in a way that the index $i=1$ and its related contents disappear but the index $i=2$ is replaced by the index $0$ in the statement of Theorem \ref{thm:PM type-Stable-Case2}. This is mainly due to property (3) in Theorem \ref{thm:uniform-parabolic}.

\section{Anomalous blow-up with vanishing energy for Perona-Malik type equations}\label{sec:proof-PM-anomalous}

This section aims to prove the most difficult result, Theorem \ref{thm:PM type-Unstable}, on the existence of weak solutions $u$ to problem (\ref{main-ibP}) showing anomalous blow-up of $u_x$ and energy dissipation in a certain sense as $t$ approaches some final time $0<T\le\infty$ for all nonconstant initial data $u_0$ when the diffusion flux $\sigma$ is of the Perona-Malik type.

We begin with the proof of Case 1: $0=m_0<M_0$. Then Case 2: $m_0<M_0=0$ follows in the same way by a symmetric argument; we skip the details.
The remaining one, Case 3: $m_0<0<M_0$, will be elaborated after the proof of Case 1.

\subsection{Case 1: $0=m_0<M_0$} For clarity, we divide the proof of this case into many steps.

\subsection*{Setup for iteration}
Let $0<r_0<\min\{\sigma(s_2),1\}$ be any number such that
\begin{equation}\label{proof-thm-2.7-1}
s^-_{r_0}<M_0<s^+_{r_0},
\end{equation}
and let $r_0'\in(r_0,\sigma(s_2))$ be close enough to $r_0$ so that
\[
s^+_{r_0}-s^+_{r_0'}<1.
\]
Let $0<r_1<\min\{r_0,2^{-1}\}$, and let $r_1'\in(r_1,r_0)$ be sufficiently close to $r_1$ so that
\[
s^+_{r_1}-s^+_{r_1'}<2^{-1}.
\]
Continuing this process indefinitely, we obtain two sequences $r_j,r_j'\in(0,\sigma(s_2))$ such that for each integer $j\ge0,$
\begin{equation}\label{proof-thm-2.7-10}
\left\{\begin{array}{l}
         0<r_j<\min\{r_{j-1},2^{-j}\}, \\
         r_j<r_j'<r_{j-1}, \\
         s^+_{r_j}-s^+_{r_j'}<2^{-j},
       \end{array}
 \right.
\end{equation}
where $r_{-1}:=\sigma(s_2)$, and that (\ref{proof-thm-2.7-1}) is satisfied.

\begin{figure}[ht]
\begin{center}
\includegraphics[scale=0.6]{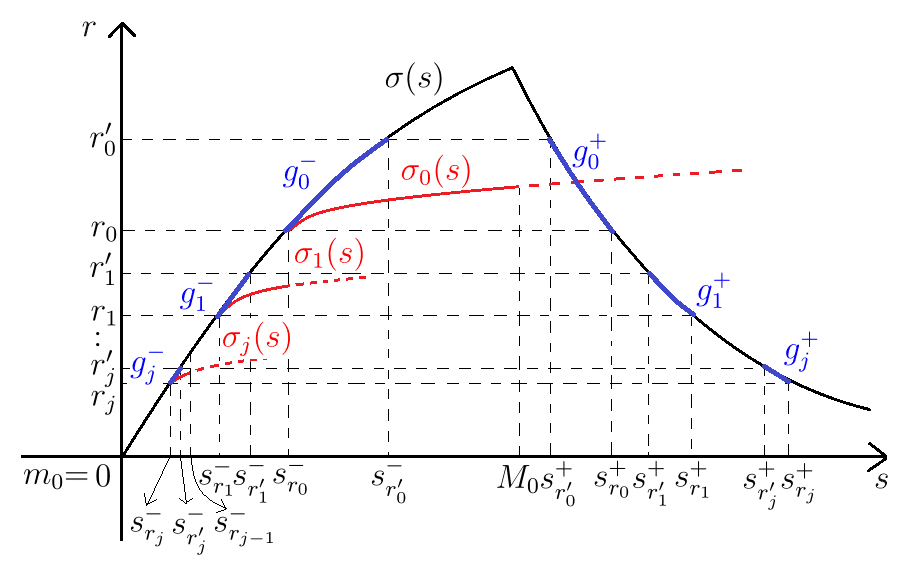}
\end{center}
\caption{Iteration scheme for Case 1: $0=m_0<M_0$}
\label{fig6}
\end{figure}

Next, for each integer $j\ge0$, choose a function $\sigma_j\in C^3(\R)$ such that
\begin{equation}\label{proof-thm-2.7-3}
\left\{\begin{array}{l}
         \mbox{$\sigma_j=\sigma$ on $[0,s^-_{r_{j}}]$,} \\
         \mbox{$\sigma_0<\min\{\sigma,r_0'\}$ on $(s^-_{r_0},M_0]$ for $j=0$,} \\
         \mbox{$\sigma_j<\min\{\sigma,r_j'\}$ on $(s^-_{r_j},s^-_{r_{j-1}}]$ for $j\ge 1$,} \\
         \mbox{$\sigma_j'>0$ in $\R$}.
       \end{array}
 \right.
\end{equation}
See Figure \ref{fig6} for an illustration of the constructions so far.

By Theorem \ref{thm:uniform-parabolic}, there exists a solution $u^{(0)*}\in C^{2,1}(\bar\Omega\times[t_0,\infty))$ to problem $(j=0)$ \begin{equation*}
\left\{
\begin{array}{ll}
  u^{(0)*}_t=(\sigma_0(u^{(0)*}_x))_{x} & \mbox{in $\Omega\times(t_0,\infty)$,} \\
  u^{(0)*}=u_0 & \mbox{on $\Omega\times\{t=t_0\}$,} \\
  u^{(0)*}_x=0 & \mbox{on $\partial\Omega\times(t_0,\infty)$,}
\end{array}
\right.
\end{equation*}
with $t_0:=0$, satisfying the properties in the statement of the theorem. By (3) and (4) of Theorem \ref{thm:uniform-parabolic}, we can choose the first time $t_1\in(t_0,\infty)$ at which $\max_{\bar\Omega}u^{(0)*}_x(\cdot,t_1)=s^-_{r_0}$. Also, from (3) of Theorem \ref{thm:uniform-parabolic}, it follows that $\min_{\bar\Omega}u^{(0)*}_x(\cdot,t)=0$ for all $t\in[t_0,t_1]$. We define $u_1=u^{(0)*}(\cdot,t_1)\in C^{2+\alpha}(\bar\Omega)$ so that $u_1'=0$ on $\partial\Omega$.

Again, by Theorem \ref{thm:uniform-parabolic}, there exists a solution $u^{(1)*}\in C^{2,1}(\bar\Omega\times[t_1,\infty))$ to problem $(j=1)$ \begin{equation*}
\left\{
\begin{array}{ll}
  u^{(1)*}_t=(\sigma_1(u^{(1)*}_x))_{x} & \mbox{in $\Omega\times(t_1,\infty)$,} \\
  u^{(1)*}=u_1 & \mbox{on $\Omega\times\{t=t_1\}$,} \\
  u^{(1)*}_x=0 & \mbox{on $\partial\Omega\times(t_1,\infty)$,}
\end{array}
\right.
\end{equation*}
satisfying the properties in the statement of the theorem. By (3) and (4) of Theorem \ref{thm:uniform-parabolic}, we can choose the first time $t_2\in(t_1,\infty)$ at which $\max_{\bar\Omega}u^{(1)*}_x(\cdot,t_2)=s^-_{r_1}$. It also follows that $\min_{\bar\Omega}u^{(1)*}_x(\cdot,t)=0$ for all $t\in[t_1,t_2]$. We define $u_2=u^{(1)*}(\cdot,t_2)\in C^{2+\alpha}(\bar\Omega)$ so that $u_2'=0$ on $\partial\Omega$.

Repeating this process indefinitely, we obtain a sequence $0=t_0<t_1<t_2<\cdots<\infty$, a sequence $u_j\in C^{2+\alpha}(\bar\Omega)$ with $u_j'=0$ on $\partial\Omega$, and a sequence $u^{(j)*}\in C^{2+\alpha,1+\alpha/2}(\bar\Omega\times[t_j,t_{j+1}])$ such that for each integer $j\ge0$,
\begin{equation}\label{proof-thm-2.7-2}
\left\{
\begin{array}{ll}
  \mbox{$u^{(j)*}_t=(\sigma_j(u^{(j)*}_x))_{x}$ in $\Omega\times(t_j,t_{j+1})$,} \\
  \mbox{$u^{(j)*}=u_j$ on $\Omega\times\{t=t_j\}$,} \\
  \mbox{$u^{(j)*}=u_{j+1}$ on $\Omega\times\{t=t_{j+1}\}$,} \\
  \mbox{$u^{(j)*}_x=0$ on $\partial\Omega\times(t_j,t_{j+1})$,} \\
  0=\underset{\bar\Omega}{\min} u_j'<\underset{\bar\Omega}{\max} u_j'=\left\{ \begin{array}{l}
                                                    \mbox{$M_0$ if $j=0$,} \\
                                                    \mbox{$s^-_{r_{j-1}}$ if $j\ge1$,}
                                                  \end{array}\right. \\
  \begin{split}
  0 & = \min_{\bar\Omega}u^{(j)*}_x(\cdot,s) =\min_{\bar\Omega}u^{(j)*}_x(\cdot,t) \\
  & < \max_{\bar\Omega}u^{(j)*}_x(\cdot,t)\le \max_{\bar\Omega}u^{(j)*}_x(\cdot,s)\;\;\mbox{for all $t_j\le s<t\le t_{j+1}$.}
  \end{split}
\end{array}
\right.
\end{equation}
We write $T=\underset{j\to\infty}{\lim} t_j\in (0,\infty]$.

\subsection*{Setup and subsolution in the $j^{th}$ step}

We fix an index $j\ge 0$ here and below and carry out analysis on the $j^{th}$ step.

Let
\[
\begin{split}
g_j^+ & = g^+_{r_j,r_j'}=(\sigma|_{[s^+_{r_j'},s^+_{r_j}]})^{-1}: [r_j,r_j']\to [s^+_{r_j'},s^+_{r_j}], \\
g_j^- & = g^-_{r_j,r_j'}=(\sigma|_{[s^-_{r_j},s^-_{r_j'}]})^{-1}: [r_j,r_j']\to [s^-_{r_j},s^-_{r_j'}];
\end{split}
\]
then $g_j^-\le g_j^-(r_j')=s^-_{r_j'}<s^+_{r_j'}=g_j^+(r_j')\le g_j^+$ on $[r_j,r_j']$ (see Figure \ref{fig6}).

Define
\[
Q_j=\big\{(x,t)\in\Omega\times(t_j,t_{j+1})\,|\,s^-_{r_j}< u^{(j)*}_x(x,t)< s^+_{r_j} \big\};
\]
then it follows from (\ref{proof-thm-2.7-2}) that $\emptyset\ne\bar{Q}_j\subset \Omega\times[t_j,t_{j+1}]$. 

Let $b_j=\|u^{(j)*}_t\|_{L^\infty(Q_j)}+1$. With the functions $g_j^\pm$ defined above, for each $u\in\R$, let $K_{b_j}(u)=K_{g_j^+,g_j^-;b_j}(u)$ and $U_{b_j}(u)=U_{g_j^+,g_j^-;b_j}(u)$ be defined as in Subsection \ref{sec:phase-transition}. Define an auxiliary function $v^{(j)*}$ by
\[
v^{(j)*}(x,t)=\int_{t_j}^t \sigma_j(u^{(j)*}_x(x,\tau))\,d\tau +\int_0^x u_j(y)\,dy\;\;\forall(x,t)\in\Omega\times(t_j,t_{j+1});
\]
then
\begin{equation}\label{proof-thm-2.7-7}
\left\{\begin{array}{l}
         v^{(j)*}_t=\sigma_j(u^{(j)*}_x) \\
         v^{(j)*}_x=u^{(j)*}
       \end{array}
 \right. \;\;\mbox{in $\Omega\times(t_j,t_{j+1})$},
\end{equation}
so that $w^{(j)*}:=(u^{(j)*},v^{(j)*})\in (C^{2,1}\times C^{3,1})(\bar\Omega\times[t_j,t_{j+1}]).$ From this, (\ref{proof-thm-2.7-3}) and the definition of $b_j$, it follows that
\[
\nabla w^{(j)*}\in U_{b_j}(u^{(j)}_*)\;\;\mbox{in $Q_j$;}
\]
thus $w^{(j)*}|_{Q_j}$ is a strict subsolution of differential inclusion (\ref{differential-inclusion-1}) with $Q=Q_j$, and so $0<\Gamma^{Q_j}_{w^{(j)*}}<1$, where the gauge operator $\Gamma^{Q_j}_{\cdot}=\Gamma^{Q_j}_{g_j^+,g_j^-;\cdot}$ is given as in Subsection \ref{sec:phase-transition}.

\subsection*{Admissible class in the $j^{th}$ step}
With $\epsilon=\min\{2^{-j},\Gamma^{Q_j}_{w^{(j)*}}\}$ and the above notation, let $\mathcal{A}_j=\mathcal{A}_{g_j^+,g_j^-,b_j,w^{(j)*},Q_j,\epsilon}$ and $\mathcal{A}_{j,\delta}=\mathcal{A}_{g_j^+,g_j^-,b_j,w^{(j)*},Q_j,\epsilon;\delta}$ $(\delta>0)$ be defined as in Subsection \ref{sec:phase-transition}. We then define a class $\mathfrak{A}_j$ of \emph{admissible functions} by
\[
\mathfrak{A}_j=\left\{\begin{array}{l}
                        w^{(j)}=(u^{(j)},v^{(j)})\in \\
                        (C^{2,1}\times C^{3,1})(\bar\Omega\times[t_j,t_{j+1}])
                      \end{array} \,\Big|\, \begin{array}{l}
                                              \mbox{$w^{(j)}=w^{(j)*}$ in} \\
                                              \mbox{$(\Omega\times (t_j,t_{j+1}))\setminus Q_j$,} \\
                                              \mbox{$w^{(j)}|_{Q_j}\in \mathcal{A}_j$}
                                            \end{array}
 \right\};
\]
then $w^{(j)*}\in\mathfrak{A}_j,$ and so $\mathfrak{A}_j\ne\emptyset$. Also, for each $\delta>0,$ let $\mathfrak{A}_{j,\delta}$ be the subclass of functions $w^{(j)}\in\mathfrak{A}_j$ such that $w^{(j)}|_{Q_j}\in \mathcal{A}_{j,\delta}$.

\subsection*{Baire's category method in the $j^{th}$ step}
Let $\mathfrak{X}_j$ denote the closure of $\mathfrak{A}_j$ in the space $L^\infty(\Omega\times(t_j,t_{j+1});\R^2)$; then $(\mathfrak{X}_j,L^\infty)$ becomes a nonempty complete metric space. From the definition of $\mathfrak{A}_j$, it is easy to check that
\[
\mathfrak{X}_j\subset w^{(j)*}+ W^{1,\infty}_0(\Omega\times(t_j,t_{j+1});\R^2).
\]
As in the proof of Theorem \ref{thm:Holllig type}, we also see that the set of points of continuity for the space-time gradient operator $\nabla_j=\nabla:\mathfrak{X}_j\to L^1(\Omega\times(t_j,t_{j+1});\M^{2\times 2})$, say $\mathfrak{C}_{\nabla_j}$, is dense in $\mathfrak{X}_j$.

\subsection*{Solutions over $[t_j,t_{j+1}]$ from $\mathfrak{C}_{\nabla_j}$}

Let $w^{(j)}=(u^{(j)},v^{(j)})\in \mathfrak{C}_{\nabla_j}.$ From Lemma \ref{lem:density lemma} and the density of $\mathfrak{A}_j$ in $\mathfrak{X}_j$, it follows that $\mathfrak{A}_{j,\frac{1}{k}}$ is dense in $\mathfrak{X}_j$ for each $k\in\N$. So we can take a sequence $w^{(j)}_{k}=(u^{(j)}_{k},v^{(j)}_{k})\in\mathfrak{A}_{j,\frac{1}{k}}$ such that $\|w^{(j)}_{k}-w^{(j)}\|_{L^\infty(\Omega\times(t_j,t_{j+1}))} \to 0$ as $k\to\infty$. By the continuity of the operator $\nabla_j$ at $w^{(j)}$, we have $\nabla w^{(j)}_k \to \nabla w^{(j)}$ in $L^1(\Omega\times(t_j,t_{j+1});\M^{2\times 2})$ and hence pointwise a.e. in $\Omega\times(t_j,t_{j+1})$ after passing to a subsequence if necessary.

We now state some properties of the function $w^{(j)}_k=(u^{(j)}_k,v^{(j)}_k)$ inherited from its membership $w^{(j)}_k\in\mathfrak{A}_{j,\frac{1}{k}}$:
\begin{equation}\label{proof-thm-2.7-4}
\left\{\begin{array}{l}
         w^{(j)}_k=(u^{(j)}_k,v^{(j)}_k)\in (C^{2,1}\times C^{3,1})(\bar\Omega\times[t_j,t_{j+1}]), \\
         \mbox{$w^{(j)}_k=w^{(j)*}$ in $(\Omega\times(t_j,t_{j+1})) \setminus\bar{Q}^{w^{(j)}_k}_j$ for some open set $Q^{w^{(j)}_k}_j\subset\subset Q_j$} \\
         \mbox{\quad\quad with $|\partial Q^{w^{(j)}_k}_j|=0$,} \\
         \|u^{(j)}_k-u^{(j)*}\|_{L^\infty(Q_j)}<2^{-(j+1)},\; \|(u^{(j)}_k)_t-u^{(j)*}_t\|_{L^\infty(Q_j)}<2^{-(j+1)},\\
         \quad\quad |\Gamma^{Q_j}_{w^{(j)}_k}- \Gamma^{Q_j}_{w^{(j)*}}|<\Gamma^{Q_j}_{w^{(j)*}}/2,\;\; \mbox{(by the choice of $\epsilon$ above)} \\
         \mbox{$\nabla w^{(j)}_k \in U_{b_j}(u^{(j)}_k)$ in $Q_j$, and} \\
         \int_{Q_j}\mathrm{dist}\big(\nabla w^{(j)}_k(x,t),K_{b_j}(u^{(j)}_k(x,t))\big)\,dxdt \le\frac{1}{k}|Q_j|.
       \end{array}
 \right.
\end{equation}

Let $k\to\infty$ in (\ref{proof-thm-2.7-4}); then we get
\begin{equation}\label{proof-thm-2.7-5}
\left\{\begin{array}{l}
         \mbox{$w^{(j)}=w^{(j)*}$ in $(\bar\Omega\times[t_j,t_{j+1}])\setminus Q_j$,} \\
         \mbox{$\nabla w^{(j)}=\nabla w^{(j)*}$ a.e. in $(\Omega\times(t_j,t_{j+1}))\setminus Q_j$,} \\
         \|u^{(j)}-u^{(j)*}\|_{L^\infty(Q_j)}\le 2^{-(j+1)}, \\
         \|u^{(j)}_t-u^{(j)*}_t\|_{L^\infty(Q_j)}\le 2^{-(j+1)}, \\
         |\Gamma^{Q_j}_{w^{(j)}}- \Gamma^{Q_j}_{w^{(j)*}}|\le\Gamma^{Q_j}_{w^{(j)*}}/2, \\
         \mbox{$\nabla w^{(j)}\in K_{b_j}(u^{(j)})$ a.e. in $Q_j$.}
       \end{array}
 \right.
\end{equation}
The last part of this also implies that a.e. in $Q_j$,
\begin{equation}\label{proof-thm-2.7-8}
\left\{\begin{array}{l}
         v^{(j)}_t=\sigma(u^{(j)}_x), \\
         u^{(j)}_x\in[s^-_{r_j},s^-_{r_j'}]\cup[s^+_{r_j'},s^+_{r_j}].
       \end{array}
 \right.
\end{equation}
On the other hand, it follows from (\ref{proof-thm-2.7-2}), (\ref{proof-thm-2.7-4}) and the definition of $Q_j$ that $(v^{(j)}_k)_t=v^{(j)*}_t$ and $(u^{(j)}_k)_x=u^{(j)*}_x\in[0,s^-_{r_j}]$ in $(\Omega\times(t_j,t_{j+1}))\setminus Q_j$. Here, let $k\to\infty$; then we have from (\ref{proof-thm-2.7-3}) and (\ref{proof-thm-2.7-7}) that a.e. in $(\Omega\times(t_j,t_{j+1}))\setminus Q_j$,
\begin{equation}\label{proof-thm-2.7-9}
v^{(j)}_t=v^{(j)*}_t=\sigma_j(u^{(j)*}_x)=\sigma(u^{(j)*}_x) =\sigma(u^{(j)}_x).
\end{equation}

Next, choose any $\zeta\in C^\infty(\bar\Omega\times[t_j,t_{j+1}])$ and $\tau\in[t_j,t_{j+1}]$. From (\ref{proof-thm-2.7-7}) and (\ref{proof-thm-2.7-4}), it follows as in the proof of Theorem \ref{thm:Holllig type} that
\[
\begin{split}
\int_{t_j}^\tau  & \int_0^L u^{(j)}_k\zeta_t\,dxdt \\
= & \int_{t_j}^\tau\int_0^L (v^{(j)}_k)_t\zeta_x\,dxdt - \int_{t_j}^\tau \big(v^{(j)*}_t(L,t)\zeta(L,t) -  v^{(j)*}_t(0,t)\zeta(0,t) \big)\,dt \\
& +\int_0^L\big(u^{(j)}_k(x,\tau)\zeta(x,\tau)- u_j(x,t_j)\zeta(x,t_j)  \big)\,dx.
\end{split}
\]
Let $k\to\infty$; then from (\ref{proof-thm-2.7-8}), (\ref{proof-thm-2.7-9}) and $v^{(j)*}_t=\sigma_j(u^{(j)*}_x)=0$ on $\partial\Omega\times[t_j,t_{j+1}]$, we have
\begin{equation}\label{proof-thm-2.7-14}
\begin{split}
\int_{t_j}^\tau  & \int_0^L u^{(j)}\zeta_t\,dxdt \\
& =  \int_{t_j}^\tau\int_0^L v^{(j)}_t\zeta_x\,dxdt +\int_0^L\big(u^{(j)}(x,\tau)\zeta(x,\tau)- u_j(x,t_j)\zeta(x,t_j)  \big)\,dx \\
& =  \int_{t_j}^\tau\int_0^L \sigma(u^{(j)}_x)\zeta_x\,dxdt +\int_0^L\big(u^{(j)}(x,\tau)\zeta(x,\tau)- u_j(x)\zeta(x,t_j)  \big)\,dx.
\end{split}
\end{equation}

\subsection*{Gauge estimate for $w^{(j)*}$}

Observe from (\ref{proof-thm-2.7-3}), (\ref{proof-thm-2.7-2}) and the definition of $Q_j$ that for all $(x,t)\in Q_j$,
\[
Z_{w^{(j)*}}(x,t) =\frac{u^{(j)*}_x(x,t)-g_j^-(v^{(j)*}_t(x,t))} {g_j^+(v^{(j)*}_t(x,t))-g_j^-(v^{(j)*}_t(x,t))}
 \le \frac{s^-_{r_{j-1}}}{s^+_{r_j'}-s^-_{r_j'}}=:e_j,
\]
where $j\ge1$ and $Z_{w^{(j)*}}=Z^{Q_j}_{g_j^+,g_j^-;w^{(j)*}}$ is as in Subsection \ref{sec:phase-transition}.
We thus have
\begin{equation}\label{proof-thm-2.7-6}
\Gamma^{Q_j}_{w^{(j)*}}=\frac{1}{|Q_j|}\int_{Q_j} Z_{w^{(j)*}}(x,t)\,dxdt \le e_j\;\;(j\ge1).
\end{equation}
Note also that $\underset{j\to\infty}{\lim}e_j=0$.

\subsection*{Blow-up estimate in the $j^{th}$ step}

From (\ref{proof-thm-2.7-6}), it follows that $0<\Gamma^{Q_j}_{w^{(j)*}}\le e_j\to 0$ as $j\to\infty;$ so there exists an index $j_{0}\ge 0$ such that $0<\Gamma^{Q_j}_{w^{(j)*}}<1/2$ for all integers $j\ge j_0$, and so $0<\Gamma^{Q_j}_{w^{(j)}}<1$ for all $j\ge j_0$ by the fifth of (\ref{proof-thm-2.7-5}). Thus, for all $j\ge j_0$, we have
\[
\left\{
\begin{array}{l}
  |Q^+_j|:=\big|\big\{(x,t)\in Q_j\,|\, u^{(j)}_x(x,t)\in[s^+_{r_j'},s^+_{r_j}]\big\}\big|>0, \\
  |Q^-_j|:=\big|\big\{(x,t)\in Q_j\,|\, u^{(j)}_x(x,t)\in[s^-_{r_j},s^-_{r_j'}]\big\}\big|>0.
\end{array}
\right.
\]
From this, we have
\begin{equation}\label{proof-thm-2.7-11}
\|u^{(j)}_x\|_{L^\infty(\Omega\times(t_j,t_{j+1}))}\ge s^+_{r'_j}\;\;\forall j\ge j_0.
\end{equation}
Note also that $\underset{j\to\infty}{\lim}s^+_{r'_j}=\infty$.

\subsection*{Energy estimate in the $j^{th}$ step}
We now consider the quantity
\[
\begin{split}
\frac{1}{t_{j+1}-t_j} & \int_{t_j}^{t_{j+1}}\mathcal{E}(u^{(j)}(\cdot,t)) \,dt = \frac{1}{t_{j+1}-t_j}\int_{t_j}^{t_{j+1}}\int_\Omega W(u^{(j)}_x(x,t)) \,dxdt \\
= & \frac{1}{t_{j+1}-t_j}\Big(\int_{Q^+_j}+\int_{Q^-_j} +\int_{(\Omega\times(t_j,t_{j+1}))\setminus Q_j}   \Big)W(u^{(j)}_x)=: I_1^+ +I_1^- +I_2.
\end{split}
\]

From the fifth of (\ref{proof-thm-2.7-5}), (\ref{proof-thm-2.7-8}), (\ref{proof-thm-2.7-6}) and the definition of $\Gamma^{Q_j}_\cdot$, we have
\[
\frac{|Q_j^+|}{|Q_j|}=\Gamma^{Q_j}_{w^{(j)}}\le \frac{3}{2}\Gamma^{Q_j}_{w^{(j)*}}\le \frac{3e_j}{2}\;\;(j\ge1).
\]
By the hypotheses on the diffusion flux $\sigma$ of the Perona-Malik type, there exists a number $c_0>0$ such that $0\le W(s)\le |s|$ for all $s\in\R$ with $|s|\ge c_0.$ Thus it follows from the definition of $Q^+_j$ and $e_j$ that for all sufficiently large $j\ge1$,
\[
I_1^+  \le \frac{1}{t_{j+1}-t_j} \int_{Q^+_j} u_x^{(j)}\,dxdt \le \frac{s^+_{r_j}|Q^+_j|}{t_{j+1}-t_j} \le \frac{3e_j s^+_{r_j}|Q_j|}{2(t_{j+1}-t_j)}\le \frac{3Ls^+_{r_j}s^-_{r_{j-1}}}{2(s^+_{r_j'}-s^-_{r_j'})}.
\]
Also, we easily have from (\ref{proof-thm-2.7-2}), the second of (\ref{proof-thm-2.7-5}) and the definition of $Q^-_j$ that
\[
I^-_1+I_2\le  \frac{W(s^-_{r_j'})|\Omega\times(t_j,t_{j+1})|}{t_{j+1}-t_j}= L W(s^-_{r_j'}).
\]
Combining these estimates, for all sufficiently large $j\ge1$, we get
\begin{equation}\label{proof-thm-2.7-12}
\frac{1}{t_{j+1}-t_j}\int_{t_j}^{t_{j+1}}\mathcal{E}(u^{(j)}(\cdot,t)) \,dt\le \frac{3Ls^+_{r_j}s^-_{r_{j-1}}}{2(s^+_{r_j'}-s^-_{r_j'})} +L W(s^-_{r_j'})=:d_j.
\end{equation}
Note also from (\ref{proof-thm-2.7-10}) that $\underset{j\to\infty}{\lim}d_j=0$.

\subsection*{Stability in the $j^{th}$ step} Recall the Poincar\'e inequality in 1-D: there exists a constant $C>0$, depending only on $\Omega=(0,L)$, such that
\[
\|f-\bar{f}\|_\infty\le C\|f'\|_\infty\;\;\forall f\in W^{1,\infty}(\Omega),
\]
where $\bar{f}:=\frac{1}{L}\int_0^L f(x)\,dx$. Using this, it follows from (\ref{proof-thm-2.7-2}) that for all $t\in[t_j,t_{j+1}]$,
\[
\|u^{(j)*}(\cdot,t)-\bar{u}_0\|_\infty \le C \|u^{(j)*}_x(\cdot,t)\|_\infty\le C s^-_{r_{j-1}}\;\;(j\ge1);
\]
thus by the first and third of (\ref{proof-thm-2.7-5}), for all $t\in[t_j,t_{j+1}]$,
\begin{equation}\label{proof-thm-2.7-13}
\|u^{(j)}(\cdot,t)-\bar{u}_0\|_\infty \le C s^-_{r_{j-1}} + 2^{-(j+1)}\;\;(j\ge1).
\end{equation}

\subsection*{Solutions by patching} In this final step, we show that the set
\[
\mathcal{S}:=\Bigg\{u=\sum_{j=0}^\infty u^{(j)}\chi_{\bar\Omega\times[t_j,t_{j+1})}\,\Big|\,
w^{(j)}=(u^{(j)},v^{(j)})\in\mathfrak{C}_{\nabla_j}\;\;\forall j\ge0
\Bigg\}
\]
consists of infinitely many weak solutions $u$ to problem (\ref{main-ibP}) satisfying (1)--(4).

Firstly, we check that $\mathcal{S}$ is an infinite set. For each $j\ge 0$, we can see as in the proof of Theorem \ref{thm:Holllig type} that $\mathfrak{C}_{\nabla_j}$ is an infinite set and that for $w^{(j)}_1=(u^{(j)}_1,v^{(j)}_1),w^{(j)}_2=(u^{(j)}_2,v^{(j)}_2)\in \mathfrak{C}_{\nabla_j}$, we have $w^{(j)}_1=w^{(j)}_2\Longleftrightarrow u^{(j)}_1=u^{(j)}_2$. Thus $\mathcal{S}$ must be an infinite set.

Secondly, we construct a suitable open set $G$ in $\R^2$ with $\partial\Omega\times[0,T)\subset G$.
By definition, we know that $\bar{Q}_j$ is a nonempty compact subset of $\Omega\times[t_j,t_{j+1}]$ for each $j\ge 0$. So  $F:=\underset{j\ge0}{\bigcup}\bar{Q}_j$ is closed in $\R\times(-\infty,T)$. Let $G=(\R\times(-\infty,T))\setminus F;$ then $G$ is open in $\R^2$, and $\partial\Omega\times[0,T)\subset G.$

Choose any $u\in\mathcal{S}$ so that
\[
u=\sum_{j=0}^\infty u^{(j)}\chi_{\bar\Omega\times[t_j,t_{j+1})}
\]
for some sequence $w^{(j)}=(u^{(j)},v^{(j)})\in\mathfrak{C}_{\nabla_j}$ $(j\ge 0)$. Finally, we verify that $u$ is a weak solution to problem (\ref{main-ibP}) fulfilling (1)--(4).

To prove (1), it is equivalent to show that $u\in C^{2+\alpha,1+\alpha/2}(\overline{G\cap\Omega_{t_j}})$ for all $j\in\N$. We proceed it  by induction on $j\in\N$. First, note from the first of (\ref{proof-thm-2.7-5}) that
\[
u=u^{(0)}=u^{(0)*}\in C^{2+\alpha,1+\alpha/2}\big (\overline{(\Omega\times(0,t_1))\setminus\bar{Q}_0}\big)= C^{2+\alpha,1+\alpha/2}(\overline{G\cap\Omega_{t_1}});
\]
so the inclusion holds for $j=1$. Next, assume that the inclusion is valid for $j=m\in\N$. Note that $(\bar{Q}_{m-1}\cup \bar{Q}_{m})\cap(\Omega\times\{t=t_m\})$ is a compact subset of $\Omega\times\{t=t_m\}$. Fix any point $x_0\in\Omega$ such that $(x_0,t_m)\not\in \bar{Q}_{m-1}\cup \bar{Q}_{m}.$ Then choose an open rectangle $B$ in $\R^2$ with center $(x_0,t_m)$, parallel to the axes, so small that its closure $\bar{B}$ does not intersect with $\bar{Q}_{m-1}\cup \bar{Q}_{m}$, $\Omega\times\{t=t_{m\pm1}\}$ and $\partial\Omega\times[t_{m-1},t_{m+1}]$. Let us write
\[
\begin{split}
B^+ & =\{(x,t)\in B\,|\,t>t_m\}, \\
B^- & =\{(x,t)\in B\,|\,t<t_m\}, \\
B^0 & =\{(x,t)\in B\,|\,t=t_m\}.
\end{split}
\]
We then have from the first of (\ref{proof-thm-2.7-5}) that $u=u^{(m-1)*}\in C^{2+\alpha,1+\alpha/2}(\bar B^-)$ and $u=u^{(m)*}\in C^{2+\alpha,1+\alpha/2}(\bar B^+)$. Also, by (\ref{proof-thm-2.7-2}), $u^{(m-1)*}$ and $u^{(m)*}$ are both equal to the function $u_m=u_m(x)$ on $\Omega\times\{t=t_m\},$ where $0\le u_m'\le s^-_{r_{m}}$ on  $\bar B^0$ by the definition of $Q_m$. From (\ref{proof-thm-2.7-3}), we have $\sigma_{m-1}=\sigma=\sigma_m$ on $[0,s^-_{r_{m}}]$; thus, with the help of (\ref{proof-thm-2.7-2}), we get $u\in C^{2+\alpha,1+\alpha/2}(\bar{B})$. On the other hand, also from the first of (\ref{proof-thm-2.7-5}), we have $u=u^{(m)*}\in C^{2+\alpha,1+\alpha/2}\big(\overline{ (\Omega\times(t_m,t_{m+1}))\setminus\bar{Q}_m}\big)$.
Using the induction hypothesis, we can now conclude that $u\in C^{2+\alpha,1+\alpha/2}(\overline{G\cap\Omega_{t_{m+1}}})$. Therefore, (1) holds.

Note that (2), (3) and (4) follow immediately from (\ref{proof-thm-2.7-13}), (\ref{proof-thm-2.7-11}) and (\ref{proof-thm-2.7-12}), respectively. It is also a simple consequence of (\ref{proof-thm-2.7-14}) by patching that $u$ is a weak solution to problem (\ref{main-ibP}).

The proof of  Case 1: $0=m_0<M_0$ is now complete.

\subsection{Case 3: $m_0<0<M_0$} 

We now focus on the proof of Case 3: $m_0<0<M_0$, which requires a careful hierarchical iteration process.

\subsection*{Setup for iteration hierarchy}

We describe an iteration scheme that shows a hierarchical pattern. In each of the steps of iteration, we encounter four possible outcomes among which the two make us stop the iteration and apply either the result of Case 1 or that of Case 2 to complete the proof. The other two possible outcomes keep us going further to the next step of iteration.

\subsubsection*{$0^{th}$ step of hierarchy}

Let $\max\{\sigma(s_1),-1\}<r_{1,0}<0<r_{2,0}<\min\{\sigma(s_2),1\}$ be any two numbers such that
\[
s^-_{r_{1,0}}<m_0<s^+_{r_{1,0}}<0<s^-_{r_{2,0}}<M_0<s^+_{r_{2,0}},
\]
and let $r'_{1,0}\in(\sigma(s_1),r_{1,0})$, $r'_{2,0}\in(r_{2,0},\sigma(s_2))$ be close enough to $r_{1,0}$, $r_{2,0}$, respectively, so that
\[
s^-_{r_{1,0}'}-s^-_{r_{1,0}}<1\;\;\mbox{and}\;\; s^+_{r_{2,0}}-s^+_{r_{2,0}'}<1.
\]
Then choose a function $\sigma_0\in C^3(\R)$ such that
\[
\left\{\begin{array}{l}
         \mbox{$\sigma_0=\sigma$ on $[s^+_{r_{1,0}},s^-_{r_{2,0}}]$,} \\
         \mbox{$\sigma_0>\max\{\sigma,r'_{1,0}\}$ on $[m_0,s^+_{r_{1,0}})$,} \\
         \mbox{$\sigma_0<\min\{\sigma,r'_{2,0}\}$ on $(s^-_{r_{2,0}},M_0]$,} \\
         \mbox{$\sigma_0'>0$ in $\R$}.
       \end{array}
 \right.
\]
By Theorem \ref{thm:uniform-parabolic}, we have a solution $u^{(0)*}\in C^{2,1}(\bar\Omega\times[t_0,\infty))$ to problem
\begin{equation*}
\left\{
\begin{array}{ll}
  u^{(0)*}_t=(\sigma_0(u^{(0)*}_x))_{x} & \mbox{in $\Omega\times(t_0,\infty)$,} \\
  u^{(0)*}=u_0 & \mbox{on $\Omega\times\{t=t_0\}$,} \\
  u^{(0)*}_x=0 & \mbox{on $\partial\Omega\times(t_0,\infty)$,}
\end{array}
\right.
\end{equation*}
with $t_0:=0$, satisfying the properties in the statement of the theorem. By (3) and (4) of Theorem \ref{thm:uniform-parabolic}, we can choose the first time $t_{1,1}\in(t_0,\infty)$ at which $\min_{\bar\Omega}u^{(0)*}_x(\cdot,t_{1,1})=s^+_{r_{1,0}}$ and the first time $t_{2,1}\in(t_0,\infty)$ at which $\max_{\bar\Omega}u^{(0)*}_x(\cdot,t_{2,1})=s^-_{r_{2,0}}$.

We now have $t_{2,1}\le t_{1,1}$ or $t_{1,1}\le t_{2,1}$. First, assume $t_{2,1}\le t_{1,1}$; then we have the dichotomy:
\[
m_{1,1}:=\max_{\bar\Omega}u^{(0)*}_x(\cdot,t_{1,1})=0\;\;\mbox{or}\;\; m_{1,1}>0.
\]
In any case, we have $m_{1,1}\le s^-_{r_{2,0}}$ by (3) of Theorem \ref{thm:uniform-parabolic}.

Second, assume $t_{1,1}\le t_{2,1}$; then we have the dichotomy:
\[
m_{2,1}:=\min_{\bar\Omega}u^{(0)*}_x(\cdot,t_{2,1})=0\;\;\mbox{or}\;\; m_{2,1}<0.
\]
Likewise, we have $m_{2,1}\ge s^+_{r_{1,0}}$.

Here and below, without any ambiguity, whenever $m_{1,1}$ appears in the context, it will be understood that we assume $t_{2,1}\le t_{1,1}$ in advance, and the same also applies to the case of $m_{2,1}$ with $t_{1,1}\le t_{2,1}$.

In the case that $m_{1,1}=0$ or that $m_{2,1}=0$, we stop the hierarchical iteration. In the case that $m_{1,1}>0$ or that $m_{2,1}<0$, we continue to the next step of iteration hierarchy as below.

\subsubsection*{$1^{st}$ step of hierarchy}
Assume  $m_{1,1}>0$ or  $m_{2,1}<0$. In the case that $m_{1,1}>0$, we set $t_1=t_{1,1}$. If $m_{2,1}<0$, then we take $t_1=t_{2,1}$. In both cases, we then let $u_1=u^{(0)*}(\cdot,t_1)\in C^{2+\alpha}(\bar\Omega)$ so that $u_1'=0$ on $\partial\Omega$.

First, assume $m_{1,1}>0$. Since $m_{1,1}\le s^-_{r_{2,0}}$, we have $0<\sigma(m_{1,1})\le r_{2,0}$.
Now, let $\max\{r_{1,0},-2^{-1}\}<r_{1,1}<0<r_{2,1}<\min\{\sigma(m_{1,1}),2^{-1}\}$ be any two numbers.

Second, assume $m_{2,1}<0$. Since $m_{2,1}\ge s^+_{r_{1,0}}$, we have $0>\sigma(m_{2,1})\ge r_{1,0}$.
Here, let $\max\{\sigma(m_{2,1}),-2^{-1}\}<r_{1,1}<0<r_{2,1}<\min\{r_{2,0},2^{-1}\}$ be any two numbers.

In both cases, let $r'_{1,1}\in(r_{1,0},r_{1,1})$, $r'_{2,1}\in(r_{2,1},r_{2,0})$ be close enough to $r_{1,1}$, $r_{2,1}$, respectively, so that
\[
s^-_{r_{1,1}'}-s^-_{r_{1,1}}<2^{-1}\;\;\mbox{and}\;\; s^+_{r_{2,1}}-s^+_{r_{2,1}'}<2^{-1}.
\]
Then choose a function $\sigma_1\in C^3(\R)$ such that
\[
\left\{\begin{array}{l}
         \mbox{$\sigma_1=\sigma$ on $[s^+_{r_{1,1}},s^-_{r_{2,1}}]$,} \\
         \mbox{$\sigma_1>\max\{\sigma,r'_{1,1}\}$ on $[s^+_{r_{1,0}},s^+_{r_{1,1}})$,} \\
         \mbox{$\sigma_1<\min\{\sigma,r'_{2,1}\}$ on $(s^-_{r_{2,1}},s^-_{r_{2,0}}]$,} \\
         \mbox{$\sigma_1'>0$ in $\R$}.
       \end{array}
 \right.
\]
By Theorem \ref{thm:uniform-parabolic}, we  have a solution $u^{(1)*}\in C^{2,1}(\bar\Omega\times[t_1,\infty))$ to problem
\begin{equation*}
\left\{
\begin{array}{ll}
  u^{(1)*}_t=(\sigma_1(u^{(1)*}_x))_{x} & \mbox{in $\Omega\times(t_1,\infty)$,} \\
  u^{(1)*}=u_1 & \mbox{on $\Omega\times\{t=t_1\}$,} \\
  u^{(1)*}_x=0 & \mbox{on $\partial\Omega\times(t_1,\infty)$,}
\end{array}
\right.
\end{equation*}
satisfying the properties in the statement of the theorem. By (3) and (4) of Theorem \ref{thm:uniform-parabolic}, we can choose the first time $t_{1,2}\in(t_1,\infty)$ at which $\min_{\bar\Omega}u^{(1)*}_x(\cdot,t_{1,2})=s^+_{r_{1,1}}$ and the first time $t_{2,2}\in(t_1,\infty)$ at which $\max_{\bar\Omega}u^{(1)*}_x(\cdot,t_{2,2}) =s^-_{r_{2,1}}$.

We now have $t_{2,2}\le t_{1,2}$ or $t_{1,2}\le t_{2,2}$. First, assume $t_{2,2}\le t_{1,2}$; then we have the dichotomy:
\[
m_{1,2}:=\max_{\bar\Omega}u^{(1)*}_x(\cdot,t_{1,2})=0\;\;\mbox{or}\;\; m_{1,2}>0.
\]
Note from (3) of Theorem \ref{thm:uniform-parabolic} that $m_{1,2}\le s^-_{r_{2,1}}$.

Second, assume $t_{1,2}\le t_{2,2}$; then we have the dichotomy:
\[
m_{2,2}:=\min_{\bar\Omega}u^{(1)*}_x(\cdot,t_{2,2})=0\;\;\mbox{or}\;\; m_{2,2}<0.
\]
Likewise, we have $m_{2,2}\ge s^+_{r_{1,1}}$.

Here and below, as in the $0^{th}$ step of hierarchy, whenever $m_{1,2}$ appears in the context, it will be understood that we assume $t_{2,2}\le t_{1,2}$ in advance, and the same also applies to the case of $m_{2,2}$ with $t_{1,2}\le t_{2,2}$.

In the case that $m_{1,2}=0$ or that $m_{2,2}=0$, we stop the hierarchical iteration. In the case that $m_{1,2}>0$ or that $m_{2,2}<0$, we continue to the next step of iteration hierarchy.

For a better understanding, we remark that, up to the $1^{st}$ step of hierarchy, we have $2^{1+1}=4$ possible  hierarchical scenarios of length 2: $(m_{1,1},m_{1,2})$, $(m_{1,1},m_{2,2})$, $(m_{2,1},m_{1,2})$ and $(m_{2,1},m_{2,2})$. In each of these cases, the second component depends on the first one; for instance, the $m_{1,2}$ in $(m_{1,1},m_{1,2})$ and $(m_{2,1},m_{1,2})$ are not necessarily equal. Also, if and only if the second component of a hierarchical scenario of length 2 is nonzero, we continue to the $2^{nd}$ step of hierarchy; otherwise, we stop here.

\subsubsection*{$j^{th}$ step of hierarchy}
Assume that we are given a hierarchical scenario $(m_{a_1,1},\cdots,m_{a_{j},j})$ of length $j\ge1$ whose each component is nonzero, where $a_1,\cdots,a_{j}\in\{1,2\}$. In the case that $m_{a_{j},j}=m_{1,j}>0$, we set $t_{j}=t_{1,j}$. If $m_{a_{j},j}=m_{2,j}<0$, then we take $t_{j}=t_{2,j}$. In both cases, we then set $u_{j}=u^{(j-1)*}(\cdot,t_{j})\in C^{2+\alpha}(\bar\Omega)$ so that $u_j'=0$ on $\partial\Omega$.

First, assume $m_{a_{j},j}=m_{1,j}>0$. Since $m_{1,j}\le s^-_{r_{2,j-1}}$, we have $0<\sigma(m_{1,j})\le r_{2,j-1}$.
Now, let $r_{1,j}$ and $r_{2,j}$ be any two numbers  with
\[
\max\{r_{1,j-1},-2^{-j}\}<r_{1,j}<0<r_{2,j} <\min\{\sigma(m_{1,j}),2^{-j}\}.
\]
Second, assume $m_{a_{j},j}=m_{2,j}<0$. Since $m_{2,j}\ge s^+_{r_{1,j-1}}$, we have $0>\sigma(m_{2,j})\ge r_{1,j-1}$.
Here, let $r_{1,j}$ and $r_{2,j}$ be any two numbers  with \[
\max\{\sigma(m_{2,j}),-2^{-j}\}<r_{1,j}<0 <r_{2,j}<\min\{r_{2,j-1},2^{-j}\}.
\]

In both cases, let $r'_{1,j}\in(r_{1,j-1},r_{1,j})$, $r'_{2,j}\in(r_{2,j},r_{2,j-1})$ be close enough to $r_{1,j}$, $r_{2,j}$, respectively, so that
\[
s^-_{r_{1,j}'}-s^-_{r_{1,j}}<2^{-j}\;\;\mbox{and}\;\; s^+_{r_{2,j}}-s^+_{r_{2,j}'}<2^{-j}.
\]
Then choose a function $\sigma_{j}\in C^3(\R)$ such that
\[
\left\{\begin{array}{l}
         \mbox{$\sigma_{j}=\sigma$ on $[s^+_{r_{1,j}},s^-_{r_{2,j}}]$,} \\
         \mbox{$\sigma_{j}>\max\{\sigma,r'_{1,j}\}$ on $[s^+_{r_{1,j-1}},s^+_{r_{1,j}})$,} \\
         \mbox{$\sigma_{j}<\min\{\sigma,r'_{2,j}\}$ on $(s^-_{r_{2,j}},s^-_{r_{2,j-1}}]$,} \\
         \mbox{$\sigma_{j}'>0$ in $\R$}.
       \end{array}
 \right.
\]
By Theorem \ref{thm:uniform-parabolic}, we  have a solution $u^{(j)*}\in C^{2,1}(\bar\Omega\times[t_{j},\infty))$ to problem
\begin{equation*}
\left\{
\begin{array}{ll}
  u^{(j)*}_t=(\sigma_{j}(u^{(j)*}_x))_{x} & \mbox{in $\Omega\times(t_{j},\infty)$,} \\
  u^{(j)*}=u_{j} & \mbox{on $\Omega\times\{t=t_{j}\}$,} \\
  u^{(j)*}_x=0 & \mbox{on $\partial\Omega\times(t_{j},\infty)$,}
\end{array}
\right.
\end{equation*}
satisfying the properties in the statement of the theorem. By (3) and (4) of Theorem \ref{thm:uniform-parabolic}, we can choose the first time $t_{1,j+1}\in(t_{j},\infty)$ at which $\min_{\bar\Omega}u^{(j)*}_x(\cdot,t_{1,j+1})=s^+_{r_{1,j}}$ and the first time $t_{2,j+1}\in(t_{j},\infty)$ at which $\max_{\bar\Omega}u^{(j)*}_x(\cdot,t_{2,j+1}) =s^-_{r_{2,j}}$.

We now have $t_{2,j+1}\le t_{1,j+1}$ or $t_{1,j+1}\le t_{2,j+1}$. First, assume $t_{2,j+1}\le t_{1,j+1}$; then we have the dichotomy:
\[
m_{1,j+1}:=\max_{\bar\Omega}u^{(j)*}_x(\cdot,t_{1,j+1})=0\;\;\mbox{or}\;\; m_{1,j+1}>0.
\]
Note $m_{1,j+1}\le s^-_{r_{2,j}}$.

Second, assume $t_{1,j+1}\le t_{2,j+1}$; then we have the dichotomy:
\[
m_{2,j+1}:=\min_{\bar\Omega}u^{(j)*}_x(\cdot,t_{2,j+1})=0\;\;\mbox{or}\;\; m_{2,j+1}<0.
\]
Note $m_{2,j+1}\ge s^+_{r_{1,j}}$.

Here and below, as in the previous steps of hierarchy, whenever $m_{1,j+1}$ appears in the context, it will be understood that we assume $t_{2,j+1}\le t_{1,j+1}$ in advance, and the same also applies to the case of $m_{2,j+1}$ with $t_{1,j+1}\le t_{2,j+1}$.

In the case that $m_{1,j+1}=0$ or that $m_{2,j+1}=0$, we stop the hierarchical iteration. In the case that $m_{1,j+1}>0$ or that $m_{2,j+1}<0$, we continue to the next step of iteration hierarchy.

\subsection*{Case of an infinite hierarchical scenario}

We first finish the proof of Case 3 for the case of an infinite hierarchical scenario $(m_{a_1,1},m_{a_2,2},\cdots)$ with $a_{j}\in\{1,2\}$ and $m_{a_j,j}\ne0$ for all $j\in\N$.

By our setup for iteration hierarchy, for each integer $j\ge0$, we have the following: with $r_{1,-1}:=\sigma(s_1)$ and $r_{2,-1}:=\sigma(s_2)$,
\[
\left\{
\begin{array}{l}
  \max\{r_{1,j-1},-2^{-j}\}<r_{1,j}<0<r_{2,j} <\min\{r_{2,j-1},2^{-j}\}, \\
  \left\{\begin{array}{l}
           \mbox{$r_{2,j}<\sigma(m_{1,j})$ if $j\ge1$ and $a_j=1$,} \\
           \mbox{$r_{1,j}>\sigma(m_{2,j})$ if $j\ge1$ and $a_j=2$,}
         \end{array}
   \right. \\
  \mbox{$r_{1,j-1}<r'_{1,j}<r_{1,j}$, $r_{2,j}<r'_{2,j}<r_{2,j-1}$,} \\
  \mbox{$s^-_{r'_{1,j}}-s^-_{r_{1,j}}<2^{-j}$, $s^+_{r_{2,j}}-s^+_{r'_{2,j}}<2^{-j}$,} \\
  \mbox{$\sigma_j\in C^3(\R)$, $\sigma'_j>0$ in $\R$,} \\
  \mbox{$\sigma_j=\sigma$ on $[s^+_{r_{1,j}},s^-_{r_{2,j}}]$,} \\
  \mbox{$\sigma_j>\max\{\sigma,r'_{1,j}\}$ on $\left\{\begin{array}{l}
                                                        \mbox{$[m_0, s^+_{r_{1,0}})$ if $j=0,$} \\
                                                        \mbox{$[s^+_{r_{1,j-1}}, s^+_{r_{1,j}})$ if $j\ge1,$}
                                                      \end{array}
   \right.$} \\
   \mbox{$\sigma_j<\min\{\sigma,r'_{2,j}\}$ on $\left\{\begin{array}{l}
                                                        \mbox{$(s^-_{r_{2,0}}, M_0]$ if $j=0,$} \\
                                                        \mbox{$(s^-_{r_{2,j}}, s^-_{r_{2,j-1}}]$ if $j\ge1,$}
                                                      \end{array}
   \right.$}
\end{array}
\right.
\]
\[
\left\{
\begin{array}{l}
  u^{(j)*}\in C^{2+\alpha,1+\alpha/2}(\bar\Omega\times[t_j,t_{j+1}]), \\
  \mbox{$u^{(j)*}_t=(\sigma_j(u^{(j)*}_x))_x$ in $\Omega\times(t_j,t_{j+1})$,} \\
  \mbox{$u^{(j)*}=u_j$ on $\Omega\times\{t=t_j\}$,} \\
  \mbox{$u^{(j)*}=u_{j+1}$ on $\Omega\times\{t=t_{j+1}\}$,}
\end{array} \right.
\]
and
\[
\left\{
\begin{array}{l}
  0<\underset{\bar\Omega}{\max}u'_j=\left\{\begin{array}{l}
                                             \mbox{$M_0$ if $j=0$,} \\
                                             \mbox{$m_{1,j}$ if $j\ge1$ and $a_j=1$, where $0<m_{1,j}\le s^-_{r_{2,j-1}}$,} \\
                                             \mbox{$s^-_{r_{2,j-1}}$ if $j\ge1$ and $a_j=2$,}
                                           \end{array}
   \right.\\
   0>\underset{\bar\Omega}{\min}u'_j=\left\{\begin{array}{l}
                                             \mbox{$m_0$ if $j=0$,} \\
                                             \mbox{$s^+_{r_{1,j-1}}$ if $j\ge1$ and $a_j=1$,} \\
                                             \mbox{$m_{2,j}$ if $j\ge1$ and $a_j=2$, where $0>m_{2,j}\ge s^+_{r_{1,j-1}}$,}
                                           \end{array}
   \right. \\
   \mbox{$\underset{\bar\Omega}{\min}u^{(j)*}_x(\cdot,s)\le \underset{\bar\Omega}{\min}u^{(j)*}_x(\cdot,t)<0$} \\
   \mbox{$\;\;\, <\underset{\bar\Omega}{\max}u^{(j)*}_x(\cdot,t) \le\underset{\bar\Omega}{\max}u^{(j)*}_x(\cdot,s)$ for all $t_j\le s<t\le t_{j+1}$.}
\end{array}
\right.
\]
Note also that $0=t_0<t_1=t_{a_1,1}<t_2=t_{a_2,2}<\cdots$. We write $T=\underset{j\to\infty}{\lim}t_j\in(0,\infty]$.

We now fix an index $j\ge 0$ here and below and setup the $j^{th}$ step in a way similar to the proof of Case 1. Due to the similarity, we only sketch some parts that are different from that of Case 1.

Let
\[
\begin{split}
g^+_{1,j} & =g^+_{r_{1,j}',r_{1,j}}=(\sigma|_{[s^+_{r_{1,j}'}, s^+_{r_{1,j}}]})^{-1}:[r_{1,j}',r_{1,j}]\to[s^+_{r_{1,j}'}, s^+_{r_{1,j}}], \\
g^-_{1,j} & =g^-_{r_{1,j}',r_{1,j}}=(\sigma|_{[s^-_{r_{1,j}}, s^-_{r_{1,j}'}]})^{-1}:[r_{1,j}',r_{1,j}]\to[s^-_{r_{1,j}}, s^-_{r_{1,j}'}], \\
g^+_{2,j} & =g^+_{r_{2,j},r_{2,j}'}=(\sigma|_{[s^+_{r_{2,j}'}, s^+_{r_{2,j}}]})^{-1}:[r_{2,j},r_{2,j}']\to[s^+_{r_{2,j}'}, s^+_{r_{2,j}}], \\
g^-_{2,j} & =g^-_{r_{2,j},r_{2,j}'}=(\sigma|_{[s^-_{r_{2,j}}, s^-_{r_{2,j}'}]})^{-1}:[r_{2,j},r_{2,j}']\to[s^-_{r_{2,j}}, s^-_{r_{2,j}'}].
\end{split}
\]

Define
\[
Q_{i,j} = \big\{(x,t)\in\Omega\times(t_j,t_{j+1})\,|\, s^-_{r_{i,j}}<u^{(j)*}_x(x,t)<s^+_{r_{i,j}} \big\}\;\;(i=1,2);
\]
then it follows from the above observations that $\bar{Q}_{1,j}$ and $\bar{Q}_{2,j}$ are nonempty subsets of $\Omega\times[t_j,t_{j+1}]$ whose intersection is $\emptyset$.

Let $b_{i,j}=\|u^{(j)*}_t\|_{L^\infty(Q_{i,j})}+1$ $(i=1,2)$. With the functions $g^\pm_{i,j}$ $(i=1,2)$ defined above, for each $u\in\R$, let $K_{b_{i,j}}(u)= K_{g^+_{i,j},g^-_{i,j};b_{i,j}}(u)$ and $U_{b_{i,j}}(u)= U_{g^+_{i,j},g^-_{i,j};b_{i,j}}(u)$ be defined as in Subsection \ref{sec:phase-transition}. Define an auxiliary function $v^{(j)*}:\Omega\times(t_j,t_{j+1})\to\R$ as in the proof of Case 1 so that $w^{(j)*}:=(u^{(j)*},v^{(j)*})\in (C^{2,1}\times C^{3,1})(\bar\Omega\times[t_j,t_{j+1}])$. With the gauge operators $\Gamma^{Q_{i,j}}_\cdot=\Gamma^{Q_{i,j}}_{g^+_{i,j},g^-_{i,j};\cdot}$ $(i=1,2)$ defined as in Subsection \ref{sec:phase-transition} and $\epsilon=\min\big\{{2^{-j},1-\Gamma^{Q_{1,j}}_{w^{(j)*}}, \Gamma^{Q_{2,j}}_{w^{(j)*}}\big\}}$, let $\mathcal{A}_{i,j}= \mathcal{A}_{g^+_{i,j},g^-_{i,j},b_{i,j},w^{(j)*},Q_{i,j},\epsilon}$ and $\mathcal{A}_{i,j,\delta}= \mathcal{A}_{g^+_{i,j},g^-_{i,j},b_{i,j},w^{(j)*},Q_{i,j},\epsilon;\delta}$ $(\delta>0)$ be defined also as in Subsection \ref{sec:phase-transition} for $i=1,2$. We then define a class $\mathfrak{A}_j$ of \emph{admissible functions} by
\[
\mathfrak{A}_j=\left\{\begin{array}{l}
                        w^{(j)}=(u^{(j)},v^{(j)})\in \\
                        (C^{2,1}\times C^{3,1})(\bar\Omega\times[t_j,t_{j+1}])
                      \end{array} \,\Big|\, \begin{array}{l}
                                              \mbox{$w^{(j)}=w^{(j)*}$ in} \\
                                              \mbox{$(\Omega\times (t_j,t_{j+1}))\setminus (Q_{1,j}\cup Q_{2,j})$,} \\
                                              \mbox{$w^{(j)}|_{Q_{i,j}}\in \mathcal{A}_{i,j}$ for $i=1,2$}
                                            \end{array}
 \right\};
\]
then $w^{(j)*}\in\mathfrak{A}_j,$ and so $\mathfrak{A}_j\ne\emptyset$. Also, for each $\delta>0,$ let $\mathfrak{A}_{j,\delta}$ be the subclass of functions $w^{(j)}\in\mathfrak{A}_j$ such that $w^{(j)}|_{Q_{i,j}}\in \mathcal{A}_{i,j,\delta}$ for $i=1,2$.

With the setup and constructions above, the rest of the proof of Case 3 for the infinite hierarchical scenario $(m_{a_j,j})_{j\in\N}$ is now a combinative repetition  of the procedures in that of Case 1 corresponding to the index $i=2$ and in that of Case 2 (not explicitly written above) to the index $i=1$. So we just sketch the remaining proof as follows without too much details.

\textbf{Baire's category method in the $j^{th}$ step.} This part can be just repeated the same as in the proof of Case 1.

\textbf{Solutions over $[t_j,t_{j+1}]$ from $\mathfrak{C}_{\nabla_j}$.} For this part, we may proceed in the same way as in the proof of Case 1 but have to take into account both the indices $i=1,2$ to deduce equation (\ref{proof-thm-2.7-14}) for any given $u^{(j)}$ with $w^{(j)}=(u^{(j)},v^{(j)})\in \mathfrak{C}_{\nabla_j}$ here and below.

\textbf{Gauge estimate for $w^{(j)*}$ over $Q_{i,j}$ $(i=1,2)$.} With $Z_{i,w^{(j)*}}=Z^{Q_{i,j}}_{g^+_{i,j},g^-_{i,j};w^{(j)*}}$ $(i=1,2)$ as in Subsection \ref{sec:phase-transition} and $j\ge 1$, we can get as in the proof of Case 1 and that of Case 2 (not explicitly written above) that
\[
\begin{split}
Z_{1,w^{(j)*}} & \ge \frac{s^+_{r_{1,j-1}}-s^-_{r_{1,j}'}}{s^+_{r_{1,j}}-s^-_{r_{1,j}}}=:e_{1,j} \;\;\mbox{in $Q_{1,j}$},\\
Z_{2,w^{(j)*}} & \le \frac{s^-_{r_{2,j-1}}}{s^+_{r_{2,j}'}-s^-_{r_{2,j}'}}=:e_{2,j}\;\;\mbox{in $Q_{2,j}$};
\end{split}
\]
thus $1>\Gamma^{Q_{1,j}}_{w^{(j)*}}\ge e_{1,j}$ and $0<\Gamma^{Q_{2,j}}_{w^{(j)*}}\le e_{2,j}$. Note also that $\underset{j\to\infty}{\lim}e_{2,j}=0$ and that
\[
e_{1,j}=\frac{s^+_{r_{1,j-1}}-s^+_{r_{1,j}}}{s^+_{r_{1,j}}-s^-_{r_{1,j}}}+1+ \frac{s^-_{r_{1,j}}-s^-_{r_{1,j}'}}{s^+_{r_{1,j}}-s^-_{r_{1,j}}} \to 1\;\;\mbox{as $j\to\infty$}.
\]

\textbf{Blow-up estimate in the $j^{th}$ step.} Using the previous gauge estimates, we can deduce as in the proof of Case 1 that
\[
\|u^{(j)}_x\|_{L^\infty(\Omega\times(t_j,t_{j+1}))}\ge \max\{-s^-_{r'_{1,j}},s^+_{r'_{2,j}}\}\;\;\mbox{for all sufficietly lrage $j\ge 1$.}
\]
Note also that $\underset{j\to\infty}{\lim }\max\{-s^-_{r'_{1,j}},s^+_{r'_{2,j}}\}=\infty$.

\textbf{Energy estimate in the $j^{th}$ step.}
Let us consider the quantity
\[
\begin{split}
& \frac{1}{t_{j+1}-t_j} \int_{t_j}^{t_{j+1}}\mathcal{E}(u^{(j)}(\cdot,t))\,dt = \frac{1}{t_{j+1}-t_j}\int_{t_j}^{t_{j+1}}\int_\Omega W(u^{(j)}_x(x,t))\,dxdt \\
& =\frac{1}{t_{j+1}-t_j} \Big(\int_{Q^+_{1,j}} +\int_{Q^-_{1,j}} +\int_{Q^+_{2,j}} +\int_{Q^-_{2,j}} +\int_{(\Omega\times(t_j,t_{j+1}))\setminus(Q_{1,j}\cup Q_{2,j})} \Big)W(u^{(j)}_x) \\
& =: I^+_{1}+I^-_{1}+I^+_{2}+I^-_{2}+I_{3},
\end{split}
\]
where
\[
\begin{split}
Q^+_{1,j} & :=\big\{(x,t)\in Q_{1,j}\,|\,u^{(j)}_x(x,t)\in[s^+_{r'_{1,j}}, s^+_{r_{1,j}}] \big\}, \\
Q^-_{1,j} & :=\big\{(x,t)\in Q_{1,j}\,|\,u^{(j)}_x(x,t)\in[s^-_{r_{1,j}}, s^-_{r'_{1,j}}] \big\}, \\
Q^+_{2,j} & :=\big\{(x,t)\in Q_{2,j}\,|\,u^{(j)}_x(x,t)\in[s^+_{r'_{2,j}}, s^+_{r_{2,j}}] \big\}, \\
Q^-_{2,j} & :=\big\{(x,t)\in Q_{2,j}\,|\,u^{(j)}_x(x,t)\in[s^-_{r_{2,j}}, s^-_{r'_{2,j}}] \big\}.
\end{split}
\]

As in the proof of Case 1, for all sufficiently large $j\ge1$, we have
\[
I^+_2\le \frac{3Ls^+_{r_{2,j}}s^-_{r_{2,j-1}}}{2(s^+_{r'_{2,j}} -s^-_{r'_{2,j}})}\;\;\mbox{and}\;\;
I^+_1+I^-_2+I_3\le L\max\big\{W(s^+_{r'_{1,j}}),W(s^-_{r'_{2,j}}) \big\}.
\]
On the other hand, observe
\[
\frac{|Q^-_{1,j}|}{|Q_{1,j}|}=1-\Gamma^{Q_{1,j}}_{w^{(j)}}\le \frac{3(1-e_{1,j})}{2},
\]
so that
\[
\begin{split}
I^-_1 & \le\frac{-s^-_{r_{1,j}}|Q^-_{1,j}|}{t_{j+1}-t_j}\le \frac{-3s^-_{r_{1,j}}(1-e_{1,j})|Q_{1,j}|}{2(t_{j+1}-t_j)} \\
& \le \frac{-3Ls^-_{r_{1,j}}(1-e_{1,j})}{2} = \frac{3L s^-_{r_{1,j}}}{2} \Big(\frac{s^+_{r_{1,j-1}}-s^+_{r_{1,j}}}{s^+_{r_{1,j}}-s^-_{r_{1,j}}}+ \frac{s^-_{r_{1,j}}-s^-_{r_{1,j}'}}{s^+_{r_{1,j}}-s^-_{r_{1,j}}}\Big);
\end{split}
\]
thus
\[
\begin{split}
\frac{1}{t_{j+1}-t_j} \int_{t_j}^{t_{j+1}}\mathcal{E}(u^{(j)}(\cdot,t))
\,dt & \le  \frac{3Ls^+_{r_{2,j}}s^-_{r_{2,j-1}}}{2(s^+_{r'_{2,j}} -s^-_{r'_{2,j}})} + L\max\big\{W(s^+_{r'_{1,j}}),W(s^-_{r'_{2,j}}) \big\} \\
& + \frac{3L s^-_{r_{1,j}}}{2} \Big(\frac{s^+_{r_{1,j-1}}-s^+_{r_{1,j}}}{s^+_{r_{1,j}}-s^-_{r_{1,j}}}+ \frac{s^-_{r_{1,j}}-s^-_{r_{1,j}'}}{s^+_{r_{1,j}}-s^-_{r_{1,j}}}\Big)=: d_j
\end{split}
\]
for all sufficiently large $j\ge 1$. Here, it is easy to see that $\underset{j\to\infty}{\lim}d_j=0$.

\textbf{Stability in the $j^{th}$ step.} As in the proof of Case 1, we easily have
\[
\|u^{(j)}(\cdot,t)-\bar{u}_0\|_\infty\le C\max\{-s^+_{r_{1,j-1}},s^-_{r_{2,j-1}}\}+2^{-(j+1)}\;\;(j\ge1),
\]
where $C>0$ is some constant depending only on $\Omega=(0,L).$

\textbf{Solutions by patching.} Let $\mathcal{S}$ be the set defined as in the proof of Case 1. With $G:=(\R\times(-\infty,T))\setminus\underset{j\ge0}{\bigcup}(\bar{Q}_{1,j}\cup \bar{Q}_{2,j})$, it follows from the above justifications that the set $\mathcal{S}$ consists of infinitely many weak solutions $u$ to problem (\ref{main-ibP}) fulfilling (1)--(4).

The proof of Case 3 for the infinite hierarchical scenario $(m_{a_j,j})_{j\in\N}$ is now complete.

\subsection*{Case of a finite hierarchical scenario}

Finally, we finish the proof of Case 3 for the case of a finite hierarchical scenario $(m_{a_1,1},\cdots, m_{a_k,k})$ of length $k\ge 1$ such that $a_j\in\{1,2\}$ $\forall 1\le j\le k$, $m_{a_j,j}\ne0$ $\forall 1\le j\le k-1$, and $m_{a_{k},k}=0.$

From the $0^{th}$ step to the $(k-1)^{th}$ step of hierarchy, we simply repeat the same process in the above case of an infinite hierarchical scenario to produce infinitely many weak solutions $\tilde{u}$  to problem (\ref{main-ibP}) up to time $t=t_{k}=t_{a_{k},k}$. At time $t=t_{k}$, we have the two possibilities on the datum $\tilde{u}(\cdot,t_k)=u_{k}\in C^{2+\alpha}(\bar\Omega)$ with $u_k'=0$ on $\partial\Omega$:
\[
\left\{\begin{array}{l}
         \mbox{$\underset{\bar\Omega}{\min} u_{k}' <0=m_{1,k}= \underset{\bar\Omega}{\max} u_{k}'$ if $a_{k}=1$,} \\
         \mbox{$\underset{\bar\Omega}{\min} u_{k}'=m_{2,k}=0 <\underset{\bar\Omega}{\max} u_{k}'$ if $a_{k}=2$.}
       \end{array}
 \right.
\]
In the case that $a_{k}=1,$ we can apply the result of Case 2 to get infinitely many weak solutions $\hat{u}$ to (\ref{main-ibP}) with the initial datum $u_{k}$ at $t=t_{k}$ and some final time $T\in(t_k,\infty]$; then we patch these $\hat{u}$ to the above functions $\tilde{u}$ at time $t=t_{k}$ to obtain the desired weak solutions $u$ to (\ref{main-ibP}) satisfying (1)--(4).
If $a_{k}=2,$ we instead apply the result of Case 1 and do the same job to obtain infinitely many weak solutions $u$ to (\ref{main-ibP}) fulfilling (1)--(4).

The proof of Case 3 is now complete, and therefore, Theorem \ref{thm:PM type-Unstable} is finally proved.

\section{Allocation of double-well energy for non-Fourier type equations}\label{sec:proof-NF}

In this section, we deal with the proof of the last result of the paper, Theorem \ref{thm:non-Fourier type}, on the existence of global weak solutions $u$ to problem (\ref{main-ibP}) showing energy allocation on the double-well as $t\to\infty$ for all initial data $u_0$ when the diffusion flux $\sigma$ is of the non-Fourier type.

We divide the proof into several steps.

\subsection*{Setup for iteration and subsolution}

\begin{figure}[ht]
\begin{center}
\includegraphics[scale=0.6]{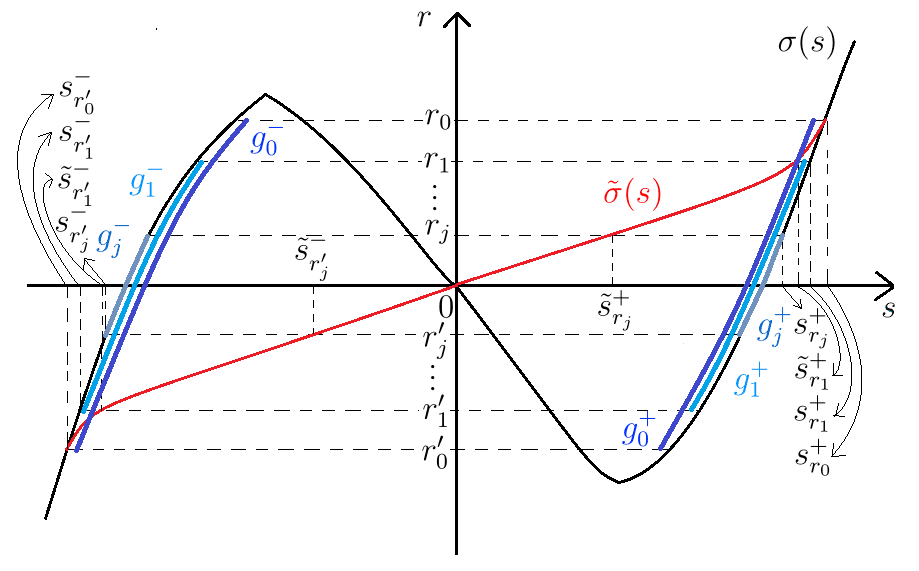}
\end{center}
\caption{Modified flux $\tilde\sigma(s)$ and iteration scheme}
\label{fig7}
\end{figure}

Let $\{r_j\}_{j\ge0}$ be a strictly decreasing sequence of positive numbers with $r_0<\min\{\sigma(s_1),-\sigma(s_2)\}$ and $\underset{j\to\infty}{\lim}r_j=0$. We write $r_j'=-r_j$ for all integers $j\ge0$. Choose a function $\tilde\sigma=\tilde\sigma_{r_0}\in C^3(\R)$ (see Figure \ref{fig7}) such that
\begin{equation}\label{proof-NF-1}
\left\{
\begin{array}{l}
  \mbox{$\tilde\sigma=\sigma$ on $(-\infty,s^-_{r_0'}]\cup[s^+_{r_0},\infty)$,} \\
  \mbox{$\tilde\sigma<\sigma$ on $(s^-_{r_0'},s^-_0]$, $\tilde\sigma>\sigma$ on $[s^+_0,s^+_{r_0})$, $\tilde\sigma(0)=0$,} \\
  \mbox{$\tilde\sigma'>0$ in $\R$.}
\end{array}
\right.
\end{equation}
By Theorem \ref{thm:uniform-parabolic}, we have a unique global solution $u^*\in C^{2,1}(\bar\Omega_\infty)$ to problem (\ref{main-ibP}), with $\sigma$ replaced by $\tilde\sigma$, satisfying the properties in the statement of the theorem. We then define an auxiliary function $v^*$ by
\[
v^*(x,t)=\int_0^t\tilde\sigma(u^*_x(x,\tau))\,d\tau +\int_0^x u_0(y)\,dy\;\;\forall(x,t)\in\Omega_\infty;
\]
then
\begin{equation}\label{proof-NF-2}
\left\{
\begin{array}{l}
  v^*_t=\tilde\sigma(u^*_x) \\
  v^*_x=u^*
\end{array}
\right.\;\;\mbox{in $\Omega_\infty$},
\end{equation}
and so $w^*:=(u^*,v^*)\in (C^{2,1}\times C^{3,1})(\bar\Omega_\infty).$

For each $j\ge 1$, let $\tilde{s}^+_{r_j}\in(0,s^+_{r_j})$ and $\tilde{s}^-_{r_j'} \in(s^-_{r_j'},0)$ denote the unique numbers with $\tilde\sigma(\tilde{s}^+_{r_j})=r_j$ and $\tilde\sigma(\tilde{s}^-_{r_j'})=r_j'$, respectively.
By Theorem \ref{thm:uniform-parabolic}, for each $j\ge1$, we can \emph{inductively} choose  a time $t_j\in(t_{j-1},\infty)$ with $t_j>j$ such that
\begin{equation}\label{proof-NF-3}
\tilde{s}^-_{r_j'}<\min_{\bar\Omega}u^*_x(\cdot,t_j)\le 0\le \max_{\bar\Omega}u^*_x(\cdot,t_j)<\tilde{s}^+_{r_j},
\end{equation}
where $t_0:=0$.

Define
\[
Q_0=\big\{(x,t)\in\Omega\times(t_0,t_{1})\,|\,s^-_{r_0'} <u^*_x(x,t)<s^+_{r_0}\big\}.
\]
Also, for each $j\ge1$, define
\[
Q_j=\big\{(x,t)\in\Omega\times(t_j,t_{j+1})\,|\, \tilde{s}^-_{r_j'}<u^*_x(x,t) <\tilde{s}^+_{r_j}\big\};
\]
then from Theorem \ref{thm:uniform-parabolic} and (\ref{proof-NF-3}), we have $Q_j=\Omega\times(t_j,t_{j+1})$. For each $j\ge1$, we also let $u_j=u^*(\cdot,t_j)\in C^{2+\alpha}(\bar\Omega)$ so that $u_j'=0$ on $\partial\Omega$.

Next, for each $j\ge0$, let
\[
g_j^\pm  =g_{r_j',r_j;\pm}=(\sigma|_{[s^\pm_{r_j'},s^\pm_{r_j}]})^{-1}: [r_j',r_j]\to[s^\pm_{r_j'},s^\pm_{r_j}];
\]
then $g^-_j\le g^-_j(r_j)=s^-_{r_j}<0<s^+_{r_j'}=g^+_j(r_j')\le g^+_j$ on $[r_j',r_j]$ (see Figure \ref{fig7}).

Fix an integer $j\ge0$ here and below. Let $b_j=\|u^*_t\|_{L^\infty(Q_j)}+1$. With the functions $g^\pm_j$ above, for each $u\in\R$, let $K_{b_j}(u)=K_{g^+_j,g^-_j;b_j}(u)$ and $U_{b_j}(u)=U_{g^+_j,g^-_j;b_j}(u)$ be defined as in Subsection \ref{sec:phase-transition}.
From (\ref{proof-NF-1}), (\ref{proof-NF-2}) and the definition of $b_j$, it follows that
\[
\nabla w^*\in U_{b_j}(u^*)\;\;\mbox{in $Q_j$};
\]
thus $w^*|_{Q_j}$ becomes a strict subsolution of differential inclusion (\ref{differential-inclusion-1}) with $Q=Q_j$ so that $0<\Gamma^{Q_j}_{w^*}<1$, where the gauge operator $\Gamma^{Q_j}_\cdot=\Gamma^{Q_j}_{g^+_j,g^-_j;\cdot}$ is given as in Subsection \ref{sec:phase-transition}.

\subsection*{Admissible class in the $j^{th}$ step}
With $\epsilon=\min\{2^{-j},j^{-1}\Gamma^{Q_j}_{w^*}\}$ and the above notation, let $\mathcal{A}_j=\mathcal{A}_{g_j^+,g_j^-,b_j,w^*,Q_j,\epsilon}$ and $\mathcal{A}_{j,\delta}=\mathcal{A}_{g_j^+,g_j^-,b_j,w^*,Q_j,\epsilon;\delta}$ $(\delta>0)$ be defined as in Subsection \ref{sec:phase-transition}. We then define a class $\mathfrak{A}_j$ of \emph{admissible functions} by the following: for $j=0$, we take
\[
\mathfrak{A}_0=\left\{\begin{array}{l}
                        w^{(0)}=(u^{(0)},v^{(0)})\in \\
                        (C^{2,1}\times C^{3,1})(\bar\Omega\times[t_0,t_{1}])
                      \end{array} \,\Big|\, \begin{array}{l}
                                              \mbox{$w^{(0)}=w^*$ in} \\
                                              \mbox{$(\Omega\times (t_0,t_{1}))\setminus Q_0$,} \\
                                              \mbox{$w^{(j)}|_{Q_0}\in \mathcal{A}_0$}
                                            \end{array}
 \right\},
\]
and for $j\ge1$, we simply let $\mathfrak{A}_j=\mathcal{A}_j$;
then $w^*|_{\Omega\times(t_0,t_1)}\in\mathfrak{A}_0$ and $w^*|_{Q_j}\in\mathfrak{A}_j$ $\forall j\ge1$. Also, for each $\delta>0,$ let $\mathfrak{A}_{j,\delta}$ be the subclass of functions $w^{(j)}\in\mathfrak{A}_j$ with $w^{(j)}|_{Q_j}\in \mathcal{A}_{j,\delta}$.

\subsection*{Baire's category method in the $j^{th}$ step}
Let $\mathfrak{X}_j$ denote the closure of $\mathfrak{A}_j$ in the space $L^\infty(\Omega\times(t_j,t_{j+1});\R^2)$; then $(\mathfrak{X}_j,L^\infty)$ becomes a nonempty complete metric space. From the definition of $\mathfrak{A}_j$, it is easy to check that
\[
\mathfrak{X}_j\subset w^*+ W^{1,\infty}_0(\Omega\times(t_j,t_{j+1});\R^2).
\]
As in the proof of Theorem \ref{thm:Holllig type}, we also see that the set of points of continuity for the space-time gradient operator $\nabla_j=\nabla:\mathfrak{X}_j\to L^1(\Omega\times(t_j,t_{j+1});\M^{2\times 2})$, say $\mathfrak{C}_{\nabla_j}$, is dense in $\mathfrak{X}_j$.

\subsection*{Solutions over $[t_j,t_{j+1}]$ from $\mathfrak{C}_{\nabla_j}$}

Let $w^{(j)}=(u^{(j)},v^{(j)})\in\mathfrak{C}_{\nabla_j}$. Repeating the same procedure as in the proof of Theorem \ref{thm:PM type-Unstable} for Case 1 under the Baire category framework above, we obtain the following:
\[ (j=0)
\left\{
\begin{array}{l}
  \mbox{$w^{(0)}=w^*$ in $(\bar\Omega\times[t_0,t_1])\setminus Q_0$}, \\
  \mbox{$\nabla w^{(0)}=\nabla w^*$ a.e. in $(\Omega\times(t_0,t_1))\setminus Q_0$}, \\
  \mbox{$v^{(0)}_t=\sigma(u^{(0)}_x)$ a.e. in $(\Omega\times(t_0,t_1))\setminus Q_0$},
\end{array}
\right.
\]
\begin{equation}\label{proof-NF-4}
\left\{
\begin{array}{l}
  \|u^{(j)}-u^*\|_{L^\infty(Q_j)}\le 2^{-(j+1)}, \\
  \|u^{(j)}_t-u^*_t\|_{L^\infty(Q_j)}\le 2^{-(j+1)}, \\
  |\Gamma^{Q_j}_{w^{(j)}}-\Gamma^{Q_j}_{w^*}|\le \frac{\Gamma^{Q_j}_{w^*}}{2j}, \\
  \mbox{$v^{(j)}_t=\sigma(u^{(j)}_x)$ a.e. in $Q_j$,} \\
  \mbox{$u^{(j)}_x\in[s^-_{r_j'},s^-_{r_j}]\cup[s^+_{r_j'}, s^+_{r_j}]$  a.e. in $Q_j$}.
\end{array}
\right.
\end{equation}
From the last of this, we see
\begin{equation}\label{proof-NF-4-2}
\begin{split}
\big\|\mathrm{dist}\big(u^{(j)}_x & ,\{s^+_0,s^-_0\}\big) \big\|_{L^\infty(\Omega\times(t_j,t_{j+1}))}\\
& \le \max\big\{|s^+_{r_j}-s^+_0|,|s^+_{r_j'}-s^+_0|, |s^-_{r_j}-s^-_0|,|s^-_{r_j'}-s^-_0|\big\}\;\;(j\ge1).
\end{split}
\end{equation}
As in the proof of Theorem \ref{thm:PM type-Unstable} for Case 1, we can also deduce that for any $\zeta\in C^\infty(\bar\Omega\times[t_j,t_{j+1}])$ and any $\tau\in[t_j,t_{j+1}]$,
\begin{equation}\label{proof-NF-5}
\begin{split}
\int_{t_j}^\tau  & \int_0^L u^{(j)}\zeta_t\,dxdt \\
& =  \int_{t_j}^\tau\int_0^L \sigma(u^{(j)}_x)\zeta_x\,dxdt +\int_0^L\big(u^{(j)}(x,\tau)\zeta(x,\tau)- u_j(x)\zeta(x,t_j)  \big)\,dx.
\end{split}
\end{equation}

\subsection*{Gauge estimate for $w^*$ over $[t_j,t_{j+1}]$}

Observe from (\ref{proof-NF-1}), (\ref{proof-NF-3}) and the definition of $Q_j$ that at any point $(x,t)\in Q_j,$
\[
\frac{\tilde d_{0,j}}{d_{1,j}} \le Z_{w^*}=\frac{u^*_x-g^-_j(v^*_t)}{g^+_j(v^*_t)-g^-_j(v^*_t)} \le\frac{\tilde d_{1,j}}{d_{0,j}},
\]
where $j\ge1$, $Z_{w^*}=Z^{Q_j}_{g^+_j,g^-_j;w^*}$ is as in Subsection \ref{sec:phase-transition}, and
\[
\left\{
\begin{array}{l}
  d_{0,j}:=\underset{[r_j',r_j]}{\min}(g^+_j-g^-_j),\;\;
  d_{1,j}:=\underset{[r_j',r_j]}{\max}(g^+_j-g^-_j), \\
  \tilde{d}_{0,j}:=\underset{[r_j',r_j]}{\min}(\tilde g_j-g^-_j), \;\;
  \tilde{d}_{1,j}:=\underset{[r_j',r_j]}{\max}(\tilde g_j-g^-_j),
\end{array}
\right.
\]
with $\tilde g_j:=(\tilde\sigma|_{[\tilde{s}^-_{r_j'}, \tilde{s}^+_{r_j}]})^{-1}:[r_j',r_j]\to[\tilde{s}^-_{r_j'}, \tilde{s}^+_{r_j}]$. We thus have
\begin{equation}\label{proof-NF-6}
\frac{\tilde d_{0,j}}{d_{1,j}} \le \Gamma^{Q_j}_{w^*} \le\frac{\tilde d_{1,j}}{d_{0,j}}\;\;(j\ge1).
\end{equation}
On the other hand, it is easy to check that
\[
\lim_{j\to\infty}\frac{\tilde d_{0,j}}{d_{1,j}} = \lim_{j\to\infty} \frac{\tilde d_{1,j}}{d_{0,j}}=\frac{-s^-_0}{s^+_0-s^-_0}=\lambda_0\in(0,1).
\]
To see this, choose two numbers $r_{1,j},\tilde{r}_{0,j}\in[r_j',r_j]$ such that
\[
g^+_j(r_{1,j})-g^-_j(r_{1,j})=d_{1,j}\;\;\mbox{and}\;\; \tilde g_j(\tilde{r}_{0,j})-g^-_j(\tilde{r}_{0,j})=\tilde{d}_{0,j}.
\]
As $j\to\infty$, we have $r_{1,j},\tilde{r}_{0,j}\to 0$, and so
\[
d_{1,j}=g^+_j(r_{1,j})-g^-_j(r_{1,j})\longrightarrow s^+_0-s^-_0,
\]
\[
\tilde{d}_{0,j}=\tilde g_j(\tilde{r}_{0,j})-g^-_j(\tilde{r}_{0,j}) \longrightarrow -s^-_0.
\]
Thus $\underset{j\to\infty}{\lim}\frac{\tilde d_{0,j}}{d_{1,j}}= \lambda_0$, and the other limit can be shown in the same way.

\subsection*{Energy estimate in the $j^{th}$ step}

We now consider the quantity
\[
\begin{split}
\frac{1}{t_{j+1}-t_j} \int_{t_j}^{t_{j+1}} \mathcal{E}(u^{(j)} (\cdot,t))  & \,dt= \frac{1}{t_{j+1}-t_j}\int_{t_j}^{t_{j+1}}\int_\Omega W(u^{(j)}_x(x,t))\,dxdt \\
& = \frac{1}{t_{j+1}-t_j}\bigg(\int_{Q^+_j}+\int_{Q^-_j} \bigg)W(u^{(j)}_x)=: I^+_j+ I^-_j,
\end{split}
\]
where $j\ge 1$, and $Q_j=\Omega\times(t_j,t_{j+1})$ is divided into the two sets $Q^\pm_j$ up to measure zero by the last of (\ref{proof-NF-4}):
\[
\begin{split}
Q^\pm_j=\big\{(x,t)\in Q_j\,|\, u^{(j)}_x(x,t)\in [s^\pm_{r_j'},s^\pm_{r_j}]\big\}.
\end{split}
\]

From the third of (\ref{proof-NF-4}), (\ref{proof-NF-6}) and the definition of $\Gamma^{Q_j}_\cdot$, we have
\[
\begin{split}
\Big(1-\frac{1}{2j}\Big)\frac{\tilde{d}_{0,j}}{d_{1,j}} &  \le\Big(1- \frac{1}{2j}\Big) \Gamma^{Q_j}_{w^*}\le \frac{|Q^+_j|}{|Q_j|}=\Gamma^{Q_j}_{w^{(j)}} \\
&\le\Big(1+\frac{1}{2j}\Big) \Gamma^{Q_j}_{w^*} \le\Big(1+\frac{1}{2j}\Big)\frac{\tilde{d}_{1,j}}{d_{0,j}};
\end{split}
\]
hence
\[
\Big(1-\frac{1}{2j}\Big)\frac{\tilde{d}_{0,j}}{d_{1,j}}|Q_j| \le |Q^+_j| \le\Big(1+\frac{1}{2j}\Big)\frac{\tilde{d}_{1,j}}{d_{0,j}} |Q_j|.
\]
Using this, we get
\[
\begin{split}
L\Big(1-\frac{1}{2j}\Big)\frac{\tilde{d}_{0,j}}{d_{1,j}} W(s^+_0) &  \le \frac{1}{t_{j+1}-t_j}W(s^+_0) |Q^+_j|\le I^+_j \\
& \le \frac{1}{t_{j+1}-t_j}\max\big\{W(s^+_{r_j}), W(s^+_{r_j'})\big\}|Q^+_j| \\
& \le L\Big(1+\frac{1}{2j}\Big)\frac{\tilde{d}_{1,j}}{d_{0,j}} \max\big\{W(s^+_{r_j}), W(s^+_{r_j'})\big\}.
\end{split}
\]
Note here that the far left and right terms are approaching $L\lambda_0 W(s^+_0)$ as $j\to\infty$; that is, $\underset{j\to\infty}{\lim}I^+_j=L\lambda_0 W(s^+_0)$. In a similar way, one can check that $\underset{j\to\infty}{\lim}I^-_j=L(1-\lambda_0)W(s^-_0)$; we omit the details. Therefore, as $j\to\infty$,
\begin{equation}\label{proof-NF-7}
\frac{1}{t_{j+1}-t_j} \int_{t_j}^{t_{j+1}}\mathcal{E}(u^{(j)}(\cdot,t)) \,dt \longrightarrow L\big(\lambda_0 W(s^+_0)+ (1-\lambda_0)W(s^-_0)\big).
\end{equation}

\subsection*{Solutions by patching} At this final stage, we check that the set
\[
\mathcal{S}:=\Bigg\{u=\sum_{j=0}^\infty u^{(j)}\chi_{\bar\Omega\times[t_j,t_{j+1})}\,\Big|\,
w^{(j)}=(u^{(j)},v^{(j)})\in\mathfrak{C}_{\nabla_j}\;\;\forall j\ge0
\Bigg\}
\]
consists of infinitely many global weak solutions $u$ to problem (\ref{main-ibP}) satisfying (1)--(3).

As in the proof of Theorem \ref{thm:PM type-Unstable} for Case 1, we easily see that $\mathcal{S}$ is an infinite set.

Let $u\in\mathcal{S}$. Then (1) follows from the first of (\ref{proof-NF-4}) and (4) of Theorem \ref{thm:uniform-parabolic}. Moreover, (2) and (3) are an immediate consequence of (\ref{proof-NF-4-2}) and (\ref{proof-NF-7}), respectively.

The proof of Theorem \ref{thm:non-Fourier type} is finally complete.

\section{Proof of density lemma}\label{sec:density-proof}

In this final section, we provide a long proof of the density lemma, Lemma \ref{lem:density lemma}.

To start the proof, fix any $\delta>0$, and let $w=(u,v)\in\mathcal{A}$; that is,
\begin{equation}\label{proof-density-1}
\left\{
\begin{array}{l}
  \mbox{$w\in (C^{2,1}\times C^{3,1})(\bar{Q})$, $w$ is Lipschitz in $Q$,} \\
  \mbox{$w=w^*$ in $Q\setminus\bar{Q}^w$ for some open set $Q^w\subset\subset Q$ with $|\partial Q^w|=0$,} \\
  \mbox{$\|u-u^*\|_{L^\infty(Q)}<\epsilon/2$, $\|u_t-u_t^*\|_{L^\infty(Q)}<\epsilon/2$, $|\Gamma^Q_w-\Gamma^Q_{w^*}|<\epsilon/2$, } \\
  \mbox{$\nabla w\in U_b(u)$ in $Q$}.
\end{array}
\right.
\end{equation}
Let $\eta>0$. Our goal is to construct a function $w_\eta=(u_\eta,v_\eta)\in\mathcal{A}_\delta$ with $\|w_\eta-w\|_{L^\infty(Q)}<\eta$, that is, a function $w_\eta\in (C^{2,1}\times C^{3,1})(\bar{Q})$ such that
\begin{equation}\label{proof-density-goal}
\left\{
\begin{array}{l}
  \mbox{$w_\eta$ is Lipschitz in $Q$, $\|w_\eta-w\|_{L^\infty(Q)}<\eta$,} \\
  \mbox{$w_\eta=w^*$ in $Q\setminus\bar{Q}^{w_\eta}$ for some open set $Q^{w_\eta}\subset\subset Q$ with $|\partial Q^{w_\eta}|=0$,} \\
  \mbox{$\|u_\eta-u^*\|_{L^\infty(Q)}<\epsilon/2$, $\|(u_\eta)_t-u_t^*\|_{L^\infty(Q)}<\epsilon/2$, $|\Gamma^Q_{w_\eta}-\Gamma^Q_{w^*}|<\epsilon/2$, } \\
  \mbox{$\nabla w_\eta\in U_b(u_\eta)$ in $Q$, $\int_{Q}\mathrm{dist}\big(\nabla w_\eta(x,t),K_b(u_\eta(x,t)) \big)\,dxdt\le\delta|Q|$}.
\end{array}
\right.
\end{equation}

As the construction is long and complicated, we divide it into many steps.

\subsection*{Step 1} We first choose an open set $G\subset\subset Q\setminus\partial Q^w$ with $|\partial G|=0$ such that
\begin{equation}\label{proof-density-2}
\int_{Q\setminus\bar{G}} \mathrm{dist}\big(\nabla w(x,t),K_b(u(x,t))\big)\,dxdt \le\frac{\delta}{k}|Q|,
\end{equation}
where $k\in\N$ is to be specified later.

By (\ref{proof-density-1}), we have
\[
(u_x,v_t)\in \tilde{U}\;\;\mbox{and}\;\;|u_t|<b\;\;\mbox{on $\bar{G}$};
\]
so
\begin{equation}\label{proof-density-3}
d':=\min_{\bar{G}}\mathrm{dist}\big((u_x,v_t),\partial\tilde{U}\big)>0 \;\;\mbox{and}\;\;b':=b-\max_{\bar{G}}|u_t|>0.
\end{equation}
Also, by (\ref{proof-density-1}), we have
\[
d'':=(\epsilon/2-|\Gamma^Q_w-\Gamma^Q_{w^*}|)/2>0.
\]
By the uniform continuity of $g^\pm$ on $[r_1,r_2]$, we can choose a $\kappa>0$ such that
\begin{equation}\label{proof-density-4}
\mbox{$|g^\pm(r)-g^\pm(r')|\le d''/l$ whenever $r,r'\in[r_1,r_2]$ and $|r-r'|\le\kappa,$}
\end{equation}
where $l\in\N$ will be chosen later.

Now, we choose finitely many disjoint open squares $Q_1,\cdots,Q_N\subset G$, parallel to the axes, such that
\begin{equation}\label{proof-density-5}
\int_{G\setminus(\cup_{i=1}^N\bar{Q}_i)} \mathrm{dist}\big(\nabla w(x,t),K_b(u(x,t))\big)\,dxdt \le\frac{\delta}{k}|Q|.
\end{equation}
Dividing these squares $Q_1,\cdots,Q_N$ into disjoint open sub-squares up to measure zero if necessary, we can have that
\begin{equation}\label{proof-density-6}
\big|\big(u_x(x_1,t_1),v_t(x_1,t_1)\big)-\big(u_x(x_2,t_2),v_t(x_2,t_2) \big)\big|\le\min\Big\{\frac{\delta}{4k},\frac{d'}{4}, \frac{\kappa}{2}\Big\}
\end{equation}
and
\begin{equation}\label{proof-density-7}
|u_t(x_1,t_1)-u_t(x_2,t_2)|\le\frac{b'}{4}
\end{equation}
whenever $(x_1,t_1),(x_2,t_2)\in\bar{Q}_i$ and $1\le i\le N$.

\subsection*{Step 2}
Consider the two continuous functions $f^\pm:[0,\infty)\to[0,\infty)$ defined by
\[
f^\pm(\beta)=\min_{r,r'\in[r_1,r_2]}\big|\big(g^\pm(r)\mp\beta,r\big)- \big(g^\pm(r'),r'\big)\big| \;\;\forall \beta\ge0;
\]
then $f^\pm$ are strictly increasing on $[0,\infty)$ with $f^\pm(0)=0$ and $\underset{\beta\to\infty}{\lim}f^\pm(\beta)=\infty$. So we can choose two unique numbers $\beta^\pm_{k}>0$ such that $f^\pm(\beta^\pm_{k})=\frac{\delta}{k-1},$ where we let $k\in\N$ satisfy
\begin{equation}\label{proof-density-k-extra-1}
k\ge2.
\end{equation}
We also choose $k\in\N$ so large that
\begin{equation}\label{proof-density-k-extra-2}
g^-(r)+\beta^-_{k}<g^+(r)-\beta^+_{k}\;\;\forall r\in[r_1,r_2].
\end{equation}
We then define two disjoint bounded domains $D^\pm_{k}\subset\tilde{U}$ by
\[
\begin{split}
D^+_{k} & =\big\{(s,r)\in\R^2\,|\, r_1<r<r_2, g^+(r)-\beta^+_{k}<s<g^+(r) \big\}, \\
D^-_{k} & =\big\{(s,r)\in\R^2\,|\, r_1<r<r_2, g^-(r)<s<g^-(r)+\beta^-_{k} \big\}.
\end{split}
\]

Fix an index $i\in\mathcal{I}:=\{1,\cdots,N\}$, and let us denote by $(x_i,t_i)$ the center of the square $Q_i$. We also write
\[
(s_i,\gamma_i)=\big(u_x(x_i,t_i),v_t(x_i,t_i)\big)\;\;\mbox{and}\;\; c_i=u_t(x_i,t_i);
\]
then $(s_i,\gamma_i)\in\tilde{U}$ and $|c_i|\le b-b'.$

Next, we split the index set $\mathcal{I}$ into the three sets
\[
\mathcal{I}^{\pm}_1 :=\big\{i\in\mathcal{I}\,|\, \mathrm{dist}\big((s_i,\gamma_i),\tilde{K}^\pm\big)\le\delta/k\big\} \;\;\mbox{and}\;\;
\mathcal{I}_2 :=\mathcal{I}\setminus(\mathcal{I}^+_1\cup \mathcal{I}^-_1).
\]
It is then easy to check that $(s_i,\gamma_i)\in D^\pm_k$ for all $i\in\mathcal{I}^{\pm}_1$; thus $\mathcal{I}^{+}_1\cap \mathcal{I}^{-}_1=\emptyset$.

Let $i\in\mathcal{I}^+_1$. Then there exists a $\bar{\gamma_i}\in [r_1,r_2]$ such that
\begin{equation}\label{proof-density-8}
\big|(s_i,\gamma_i)-\big(g^+(\bar{\gamma_i}),\bar{\gamma_i}\big)\big|\le \frac{\delta}{k}.
\end{equation}
Let $(x,t)\in Q_i$; then by (\ref{proof-density-1}), (\ref{proof-density-6}) and (\ref{proof-density-8}), we have
\[
\left|\begin{pmatrix} u_x(x,t) & u_t(x,t) \\ v_x(x,t) & v_t(x,t)  \end{pmatrix} - \begin{pmatrix} g^+(\bar{\gamma_i}) & u_t(x,t) \\ u(x,t) & \bar{\gamma_i}  \end{pmatrix}  \right| \le\frac{5\delta}{4k},
\]
where
\[
\begin{pmatrix} g^+(\bar{\gamma_i}) & u_t(x,t) \\ u(x,t) & \bar{\gamma_i}  \end{pmatrix} \in K_b(u(x,t)).
\]
So
\[
\int_{Q_i}\mathrm{dist}\big(\nabla w(x,t),K_b(u(x,t)) \big)\,dxdt \le\frac{5\delta}{4k}|Q_i|,
\]
and the same inequality can be shown to hold in a similar manner when $i\in\mathcal{I}^-_1.$
We thus have
\begin{equation}\label{proof-density-9}
\int_{\cup_{i\in\mathcal{I}^+_1\cup \mathcal{I}^-_1}Q_i}\mathrm{dist}\big(\nabla w(x,t),K_b(u(x,t)) \big)\,dxdt \le\frac{5\delta}{4k}|Q|.
\end{equation}

\subsection*{Step 3} Fix an index $i\in\mathcal{I}_2$. Since
 $\mathrm{dist}\big((s_i,\gamma_i),\tilde{K}^\pm\big)>\delta/k$, we can choose two positive numbers $\lambda_{i}^\pm$ so that
\begin{equation}\label{proof-density-10}
\left\{
\begin{array}{l}
  \mathrm{dist}\big((s_i+\lambda_{i}^+,\gamma_i),\tilde{K}^+\big)=\delta/k, \;(s_i+\lambda_{i}^+,\gamma_i)\in D^+_{k}, \\
  \mathrm{dist}\big((s_i-\lambda_{i}^-,\gamma_i),\tilde{K}^-\big)=\delta/k, \;(s_i-\lambda_{i}^-,\gamma_i)\in D^-_{k}.
\end{array}
\right.
\end{equation}

For a given $\epsilon_i>0$ to be specified later, we can apply Lemma \ref{lem:rank-1-smooth-approx} to obtain a function $\omega_i=(\varphi_i,\psi_i)\in C^\infty_c(Q_i;\R^2)$ such that
\begin{itemize}
\item[(a)] $\|\omega_i\|_{L^\infty(Q_i)}<\epsilon_i$, $\|(\varphi_i)_t\|_{L^\infty(Q_i)}<\epsilon_i$, $\|(\psi_i)_t\|_{L^\infty(Q_i)}<\epsilon_i$,
\item[(b)] $-\lambda_{i}^-\le(\varphi_i)_x\le\lambda_{i}^+$ in $Q_i$,
\item[(c)] $\left\{\begin{array}{l}
              \big| |\{  (x,t)\in Q_i\, |\, (\varphi_i)_x(x,t) =\lambda_{i}^+ \} |-\frac{\lambda_{i}^-}{\lambda_{i}^++\lambda_{i}^-}|Q_i|\big| < \epsilon_i,\\
              \big| |\{  (x,t)\in Q_i\, |\, (\varphi_i)_x(x,t) =-\lambda_{i}^- \} |-\frac{\lambda_{i}^+}{\lambda_{i}^++\lambda_{i}^-}|Q_i|\big| < \epsilon_i,
              \end{array}\right.
$
\item[(d)] $(\psi_i)_x=\varphi_i$ in $Q_i$, and
\item[(e)] $\int_{x_{i,1}}^{x_{i,2}}\varphi_i(x,t)\,dx=0$  for all $t_{i,1}<t<t_{i,2}$,
\end{itemize}
where $Q_i=(x_{i,1},x_{i,2})\times(t_{i,1},t_{i,2})$. We now define
\begin{equation}\label{proof-density-11}
w_\eta=(u_\eta,v_\eta)=w+\sum_{i\in\mathcal{I}_2}\omega_i\chi_{Q_i}\;\;\mbox{in $Q$}.
\end{equation}

\subsection*{Step 4} To finish the proof, we show that upon choosing suitable numbers $l,k\in\N$ and $\epsilon_i>0$ $(i\in\mathcal{I}_2)$, the function $w_\eta=(u_\eta,v_\eta)$ defined in (\ref{proof-density-11}) satisfies the desired properties in (\ref{proof-density-goal}). Since this step consists of many arguments to verify, we separate it into several substeps.

\subsubsection*{\textbf{Substep 4-1}} We begin with simpler parts to prove.

From (\ref{proof-density-1}) and (\ref{proof-density-11}), it follows that
\begin{equation}\label{goal-1}
\mbox{$w_\eta\in (C^{2,1}\times C^{3,1})(\bar Q)$  is Lipschitz in $Q$.}
\end{equation}

Set $Q^{w_\eta}=G\cup Q^w$; then $Q^{w_\eta}\subset\subset Q$ is an  open set with $|\partial Q^{w_\eta}|=0$. From (\ref{proof-density-1}) and (\ref{proof-density-11}), we also have
\begin{equation}\label{goal-2}
w_\eta=w=w^*\;\;\mbox{in $Q\setminus\bar Q^{w_\eta}$}.
\end{equation}

Next, we choose
\begin{equation}\label{proof-density-15}
\mbox{$\epsilon_i<\min\Big\{\eta,\; \frac{\epsilon}{2} -\|u-u^*\|_{L^\infty(Q)},\, \frac{\epsilon}{2}-\|u_t-u^*_t\|_{L^\infty(Q)}\Big \}$\; for all $i\in\mathcal{I}_2$;}
\end{equation}
then  from (a) and (\ref{proof-density-11}), we have
\begin{equation}\label{goal-3}
\begin{cases}
\|w_\eta-w\|_{L^\infty(Q)}<\eta,\\
\|u_\eta-u^*\|_{L^\infty(Q)}<{\epsilon}/{2},\\
\|(u_\eta)_t-u^*_t\|_{L^\infty(Q)}<{\epsilon}/{2}.
\end{cases}
\end{equation}

\subsubsection*{\textbf{Substep 4-2}} In this substep, we show that
\begin{equation}\label{goal-4}
\nabla{w}_\eta(x,t)\in U_b(u_\eta(x,t))\;\;\forall (x,t)\in Q.
\end{equation}

Note from (\ref{proof-density-7}), (\ref{proof-density-11}) and (a) that
\[
\begin{split}
|(u_\eta)_t(x,t)| & =|u_t(x,t)+(\varphi_i)_t(x,t)|\le |u_t(x,t)-c_i|+|c_i|+|(\varphi_i)_t(x,t)| \\
& \le \frac{b'}{4}+b-b'+ \epsilon_i< b-\frac{b'}{2}\;\;\forall (x,t)\in Q_i,\,\forall i\in\mathcal{I}_2,
\end{split}
\]
where we let
\begin{equation}\label{proof-density-18}
\mbox{$\epsilon_i< b'/4$\;\;$\forall i\in\mathcal{I}_2$.}
\end{equation}
This  implies that
\begin{equation}\label{proof-density-19}
|(u_\eta)_t(x,t)|< b\;\;\forall (x,t)\in Q.
\end{equation}

From (\ref{proof-density-1}), we have $v_x=u$ in $Q$. Thus by (\ref{proof-density-11}) and (d), we have
\[
\begin{split}
(v_\eta)_x (x,t) & =v_x(x,t) + (\psi_i)_x(x,t) \\
& =u(x,t) + \varphi_i(x,t)=u_\eta(x,t)\;\;\forall (x,t)\in Q_i,
\end{split}
\]
where $i\in\mathcal{I}_2$; thus
\begin{equation}\label{proof-density-20}
(v_\eta)_x (x,t)  =u_\eta(x,t)\;\;\forall (x,t)\in Q.
\end{equation}

Let $i\in\mathcal{I}_2$.
From (\ref{proof-density-3}) and (\ref{proof-density-10}), we see that
\begin{equation}\label{proof-density-21}
\mathrm{dist}\big( [(s_i-\lambda_{i}^-,\gamma_i),(s_i+\lambda_{i}^+,\gamma_i)] , \partial \tilde U \big)\ge \min\Big\{\frac{\delta}{k},d'\Big\}.
\end{equation}
Let $(x,t)\in Q_i.$ Then by (a) and (\ref{proof-density-6}),
\[
\begin{split}
\big|\big((u_\eta)_x (x,t), &  (v_\eta)_t(x,t)\big)  - \big(s_i+(\varphi_i)_x(x,t),\gamma_i\big)\big| \\
\le &   \big|\big(u_x(x,t), v_t(x,t)\big) - (s_i,\gamma_i)\big| + |(\psi_i)_t(x,t) | \\
\le & \min\Big\{\frac{\delta}{4k},\frac{d'}{4}\Big\} + \epsilon_i <\min\Big\{\frac{\delta}{2k},\frac{d'}{2}\Big\},
\end{split}
\]
where we let
\begin{equation}\label{proof-density-22}
\epsilon_i<  \min\Big\{\frac{\delta}{4k},\frac{d'}{4}\Big\}.
\end{equation}
By (b), we have $-\lambda_{i}^-\le(\varphi_i)_x(x,t)\le\lambda_{i}^+$, and so it follows from  (\ref{proof-density-21}) and the previous inequality that
\[
\big((u_\eta)_x(x,t), (v_\eta)_t(x,t)\big)\in \tilde U.
\]
Adopting (\ref{proof-density-22}) for all $i\in\mathcal{I}_2$, this inclusion holds for all $(x,t)\in Q.$

This inclusion together with (\ref{proof-density-19}) and (\ref{proof-density-20}) implies inclusion (\ref{goal-4}).

\subsubsection*{\textbf{Substep 4-3}} Here, we prove that
\begin{equation}\label{goal-5}
|\Gamma^Q_{w_\eta}-\Gamma^Q_{w^*}|<\frac{\epsilon}{2}.
\end{equation}

By the definition of the gauge operator $\Gamma^Q_\cdot$ and (\ref{proof-density-11}), we have
\[
\begin{split}
 \Gamma^Q_{w_\eta}-\Gamma^Q_w = &   \frac{1}{|Q|} \int_{Q} \big(Z_{w_\eta}(x,t)-Z_{w}(x,t)\big) \,dxdt \\
= &  \frac{1}{|Q|} \sum_{i\in\mathcal{I}_2} \int_{Q_i}  \big(Z_{w+\omega_i}(x,t)-Z_{w}(x,t)\big)\,dxdt \\
= & \frac{1}{|Q|} \sum_{i\in\mathcal{I}_2} \bigg(\int_{Q_{i,1}}  + \int_{Q^+_{i,2}}  + \int_{Q^-_{i,2}} \bigg)(Z_{w+\omega_i}-Z_{w}),
\end{split}
\]
where $\omega_i$ is defined to be zero outside its compact support, and
\[
\begin{split}
Q_{i,1} & :=\big\{(x,t)\in Q_i\,|\, (\varphi_i)_x(x,t)\not\in\{\lambda_{i}^+,-\lambda_{i}^-\} \big \}, \\
Q^+_{i,2} & :=\big\{(x,t)\in Q_i\,|\, (\varphi_i)_x(x,t)=\lambda_{i}^+  \big\}, \\
Q^-_{i,2} & :=\big\{(x,t)\in Q_i\,|\, (\varphi_i)_x(x,t)=-\lambda_{i}^- \big \}.
\end{split}
\]

The trivial part to estimate here follows from (c) as
\begin{equation}\label{goal-5-estimate-1}
\bigg|\int_{Q_{i,1}} (Z_{w+\omega_i}-Z_{w})  \bigg|\le 2|Q_{i,1}|<4\epsilon_i\;\;(i\in\mathcal{I}_2).
\end{equation}


Next, let $i\in\mathcal{I}_2$ and $(x,t)\in Q^-_{i,2}$. Then from (\ref{proof-density-6}), we have at the point $(x,t)$ here and below that
\begin{equation}\label{proof-density-25}
\begin{split}
\big|(u+\varphi_i)_x- & g^-\big((v+\psi_i)_t\big)\big|=\big|u_x  -\lambda_{i}^- -g^-\big(v_t+(\psi_i)_t\big)\big| \\
\le & |u_x -s_i| + |s_i-\lambda_{i}^- -g^-(\gamma_i)|  + \big|g^-(\gamma_i) -g^-\big(v_t+(\psi_i)_t\big)\big| \\
\le & \frac{\delta}{4k}+ |s_i-\lambda_{i}^- -g^-(\gamma_i)|  + \big|g^-(\gamma_i) -g^-\big(v_t+(\psi_i)_t\big)\big|.
\end{split}
\end{equation}
From (\ref{proof-density-10}), it follows that $\big|(s_i-\lambda_{i}^-,\gamma_i)-\big(g^-(\bar\gamma_i), \bar\gamma_i\big)\big|=\delta/k$ for some  $\bar\gamma_i\in[r_1,r_2]$.
Thus, we have from (\ref{proof-density-4}) that
\begin{equation}\label{proof-density-26}
|s_i-\lambda_{i}^--g^-(\gamma_i)|\le |s_i-\lambda_{i}^--g^-(\bar\gamma_i)|+ |g^-(\bar\gamma_i)-g^-(\gamma_i)| \le\frac{\delta}{k}+\frac{d''}{l},
\end{equation}
where we let $k\in\N$ satisfy
\begin{equation}\label{proof-density-27}
\frac{\delta}{k}\le\kappa.
\end{equation}
Likewise, with (\ref{proof-density-27}), we also have
\begin{equation}\label{proof-density-28}
|s_i+\lambda_{i}^+-g^+(\gamma_i)| \le\frac{\delta}{k}+\frac{d''}{l}.
\end{equation}
Also, from (\ref{proof-density-6}) and (a), we have
\[
\big|\gamma_i- \big (v_t+(\psi_i)_t\big)\big| \le |\gamma_i-v_t|+|(\psi_i)_t|\le \frac{\kappa}{2}+\epsilon_i\le\kappa,
\]
where we let
\begin{equation}\label{proof-density-29}
\epsilon_i\le\frac{\kappa}{2}.
\end{equation}
With the help of (\ref{proof-density-4}), this implies that
\begin{equation}\label{proof-density-30}
\big|g^-(\gamma_i) -g^-\big(v_t+(\psi_i)_t\big)\big|\le\frac{d''}{l}.
\end{equation}
Thus (\ref{proof-density-25}), (\ref{proof-density-28}) and (\ref{proof-density-30}) yield that
\[
\big|(u+\varphi_i)_x- g^-\big((v+\psi_i)_t\big)\big|\le\frac{5\delta}{4k}+\frac{2d''}{l} \le\frac{3d''}{l}
\]
if we choose $k\in\N$ so large that
\begin{equation}\label{proof-density-31}
\frac{5\delta}{4k}\le\frac{d''}{l}.
\end{equation}
Since $g^+\big((v+\psi_i)_t\big)-g^-\big((v+\psi_i)_t\big)\ge \underset{[r_1,r_2]}{\min}(g^+-g^-)=:d_{0}>0$, we thus have
\[
Z_{w+\omega_i}(x,t) \le \frac{3d''}{ld_{0}}\quad \forall\, (x,t)\in  Q_{i,2}^-;
\]
hence
\begin{equation}\label{goal-5-estimate-2}
\bigg|\int_{Q^-_{i,2}} Z_{w+\omega_i}\bigg| \le \frac{3d''}{ld_{0}}|Q^-_{i,2}|.
\end{equation}

Let $i\in\mathcal{I}_2$ and $Q_{i,2}=Q_{i,2}^+\cup Q_{i,2}^-$. We now estimate the quantity
\[
\bigg| \int_{Q_{i,2}^+} Z_{w+\omega_i} - \int_{Q_{i,2}} Z_{w}\bigg|.
\]


First, we consider
\[
\begin{split}
\bigg| & \int_{Q_{i,2}^+} Z_{w+\omega_i}(x,t) \,dxdt - \frac{\lambda_{i}^-}{\lambda_{i}^++\lambda_{i}^-}|Q_i|\bigg| \\
&\le \bigg|\int_{Q_{i,2}^+} \bigg( \frac{u_x+(\varphi_i)_x - g^-\big(v_t+(\psi_i)_t\big)}{g^+\big(v_t+(\psi_i)_t\big) -g^-\big(v_t+(\psi_i)_t\big)} \\
& \;\;\;\;\quad\quad\quad\;- \frac{s_i+\lambda_{i}^+ -g^-(\gamma_i)}{g^+\big(v_t+(\psi_i)_t\big) -g^-\big(v_t+(\psi_i)_t\big)}\bigg) \,dxdt\bigg|\\
& \;\;\;\;+ \bigg|\int_{Q_{i,2}^+} \bigg( \frac{s_i+\lambda_{i}^+ -g^-(\gamma_i)}{g^+\big(v_t+(\psi_i)_t\big) -g^-\big(v_t+(\psi_i)_t\big)} - \frac{s_i+\lambda_{i}^+ -g^-(\gamma_i)}{g^+(\gamma_i)-g^-(\gamma_i)}\bigg) \,dxdt\bigg|\\
& \;\;\;\;+ \bigg| \frac{s_i+\lambda_{i}^+-g^-(\gamma_i)}{g^+(\gamma_i)-g^-(\gamma_i)}|Q^+_{i,2}|- \frac{\lambda_{i}^-}{\lambda_{i}^++\lambda_{i}^-}|Q_i| \bigg|=:I_1+I_2+I_3.
\end{split}
\]
Here, let $\epsilon_i>0$ satisfy (\ref{proof-density-29}); then as in (\ref{proof-density-30}), we have
\begin{equation}\label{proof-density-33}
\big|g^\pm\big(v_t+(\psi_i)_t\big)-g^\pm(\gamma_i)\big| \le\frac{d''}{l}\;\;\mbox{in $Q^+_{i,2}$}.
\end{equation}
Using this, (\ref{proof-density-6}) and the fact that $(\varphi_i)_x=\lambda_{i}^+$ in $Q^+_{i,2}$, we deduce that
\[
I_1\le \frac{1}{d_0}\Big(\frac{\delta}{4k}+\frac{d''}{l}\Big)|Q^+_{i,2}|\le \frac{6d''}{5ld_0}|Q^+_{i,2}|,
\]
where we let $k\in\N$ fulfill (\ref{proof-density-31}). Having the common denominator in the integrand, we  have from (\ref{proof-density-33}) that
\[
I_2\le \big(s_i+\lambda_{i}^+-g^-(\gamma_i)\big)\frac{2d''}{ld_0^2} |Q^+_{i,2}| \le \frac{2d''d_1}{ld_0^2}|Q^+_{i,2}|,
\]
where $d_1:=\underset{[r_1,r_2]}{\max}(g^+-g^-)$.
Note from (\ref{proof-density-28}) and (c) that
\[
\begin{split}
I_3 & \le \bigg| \frac{s_i+\lambda_{i}^+-g^-(\gamma_i)}{g^+(\gamma_i)-g^-(\gamma_i)} |Q^+_{i,2}|- |Q^+_{i,2}|\bigg|+\bigg||Q^+_{i,2}|- \frac{\lambda_{i}^-}{\lambda_{i}^++\lambda_{i}^-}|Q_i| \bigg| \\
& \le \Big(\frac{\delta}{k}+\frac{d''}{l}\Big) \frac{1}{d_0}|Q^+_{i,2}|+\epsilon_i \le \frac{9d''}{5ld_0}|Q^+_{i,2}|+\epsilon_i,
\end{split}
\]
where we let $k\in\N$ also satisfy (\ref{proof-density-27}) and (\ref{proof-density-31}). Combining the estimates on $I_1,$ $I_2$ and $I_3$, we obtain
\begin{equation}\label{proof-density-34}
\bigg|\int_{Q_{i,2}^+}  Z_{w+\omega_i} (x,t)\,dxdt -  \frac{\lambda_{i}^-}{\lambda_{i}^++\lambda_{i}^-}|Q_i|\bigg|
\le \bigg(\frac{3}{d_0} + \frac{2d_1}{d_0^{2}}\bigg)\frac{d''}{l} |Q^+_{i,2}|+\epsilon_i.
\end{equation}

Second, we handle
\[
\begin{split}
\bigg|\int_{Q_{i,2}} & Z_w(x,t)\,dxdt - \frac{\lambda_{i}^-}{\lambda_{i}^++\lambda_{i}^-}|Q_i|\bigg| \\
\le &  \bigg|\int_{Q_{i,2}}  \bigg (\frac{u_x-g^-(v_t)}{g^+(v_t)-g^-(v_t)} - \frac{\lambda_{i}^-}{g^+(v_t)-g^-(v_t)}\bigg ) \,dxdt\bigg| \\
& + \bigg|\int_{Q_{i,2}}   \frac{\lambda_{i}^-}{g^+(v_t)-g^-(v_t)}\,dxdt -\frac{\lambda_{i}^-}{g^+(\gamma_i)-g^-(\gamma_i)}|Q_{i,2}|\bigg|\\
& + \bigg|\frac{\lambda_{i}^-}{g^+(\gamma_i)-g^-(\gamma_i)}|Q_{i,2}| - \frac{\lambda_{i}^-}{\lambda_{i}^++\lambda_{i}^-}|Q_i|\bigg|=:J_1+J_2+J_3.
\end{split}
\]
To estimate $J_1$, we first note from (\ref{proof-density-4}), (\ref{proof-density-6}) and (\ref{proof-density-26}) that in $Q_{i,2}$,
\[
\begin{split}
|u_x-g^-(v_t)-\lambda_{i}^-| & \le |u_x-s_i|+|g^-(\gamma_i)-g^-(v_t)| +| s_i-\lambda_{i}^- - g^-(\gamma_i)|\\
& \le \frac{\delta}{4k}+\frac{d''}{l}+\Big(\frac{\delta}{k}+\frac{d''}{l}\Big)\le\ \frac{3d''}{l},
\end{split}
\]
where we let $k\in\N$ satisfy (\ref{proof-density-27}) and (\ref{proof-density-31}). From this, we get
\[
J_1\le \frac{3d''}{ld_0}|Q_{i,2}|.
\]
From (\ref{proof-density-4}) and (\ref{proof-density-6}), we have
\[
J_2\le \frac{2\lambda_{i}^-d''}{ld_0^2}|Q_{i,2}|\le \frac{2d''d_1}{ld_0^2}|Q_{i,2}|.
\]
To estimate $J_3$, we observe from (\ref{proof-density-26}) and (\ref{proof-density-28}) that
\[
\begin{split}
\big|\lambda_{i}^++\lambda_{i}^--\big(g^+(\gamma_i)-g^-(\gamma_i) \big)\big| & =\big|s_i+\lambda_{i}^+ -g^+(\gamma_i)-\big(s_i-\lambda_{i}^- -g^-(\gamma_i)\big)\big| \\
&\le \frac{2\delta}{k}+ \frac{2d''}{l}\le \frac{18d''}{5l},
\end{split}
\]
with $k\in\N$ satisfying (\ref{proof-density-27}) and (\ref{proof-density-31}).
Note from this and (c) that
\[
\begin{split}
J_3  \le & \bigg|\frac{\lambda_{i}^-}{g^+(\gamma_i)- g^-(\gamma_i)}|Q_{i,2}| - \frac{\lambda_{i}^-}{\lambda_{i}^++\lambda_{i}^-}|Q_{i,2}|\bigg|\\
& + \bigg|\frac{\lambda_{i}^-}{\lambda_{i}^++\lambda_{i}^-}|Q_{i,2}| - \frac{\lambda_{i}^-}{\lambda_{i}^++\lambda_{i}^-}|Q_i|\bigg| \\
\le & \frac{18d''d_1}{5ld_0^2}|Q_{i,2}|+2\epsilon_i.
\end{split}
\]
By the estimates on $J_1$, $J_2$ and $J_3$, we now have
\[
\bigg|\int_{Q_{i,2}}  Z_w(x,t) \,dxdt - \frac{\lambda_{i}^-}{\lambda_{i}^++\lambda_{i}^-}|Q_i|\bigg|
 \le \bigg(\frac{3}{d_0}+ \frac{28d_1}{5d_0^2}\bigg)\frac{d''}{l}|Q_{i,2}|+2\epsilon_i.
\]

Combining this estimate with (\ref{proof-density-34}), we have
\begin{equation}\label{goal-5-estimate-3}
\begin{split}
\bigg| \int_{Q_{i,2}^+} &Z_{w+\omega_i}  - \int_{Q_{i,2}}Z_{w}\bigg|  \\
& \quad\quad\quad\le \bigg(\frac{6}{d_0}+ \frac{38d_1}{5d_0^2}\bigg)\frac{d''}{l} |Q_{i,2}|+3\epsilon_i.
\end{split}
\end{equation}

Thanks to the estimates (\ref{goal-5-estimate-1}), (\ref{goal-5-estimate-2}) and (\ref{goal-5-estimate-3}), it follows that for all $i\in\mathcal{I}_2,$
\[
\begin{split}
\bigg|\bigg ( \int_{Q_{i,1}}+&\int_{Q^+_{i,2}}+\int_{Q^-_{i,2}}\bigg) (Z_{w+\omega_i} -Z_w) \bigg| \\
& \le \bigg(\frac{9}{d_0}+ \frac{38d_1}{5d_0^2}\bigg)\frac{d''}{l}|Q_{i,2}|+7\epsilon_i\le  \frac{17d''d_1}{ld_0^2}|Q_{i}|+7\epsilon_i.
\end{split}
\]
Summing this over the indices $i\in\mathcal{I}_2,$ we then have
\[
|\Gamma^Q_{w_\eta}-\Gamma^Q_w|\le  \frac{17d''d_1}{ld_0^2} + \frac{7}{|Q|}\sum_{i\in\mathcal{I}_2}\epsilon_i< d''
\]
if $l\in\N$ is taken so large that
\begin{equation}\label{proof-density-36}
\frac{17d_1}{ld_0^2}<\frac{1}{2},
\end{equation}
$k\in\N$ satisfies (\ref{proof-density-27}) and (\ref{proof-density-31}), and the numbers $\epsilon_i>0$ $(i\in\mathcal{I}_2)$ fulfill (\ref{proof-density-29}) and
\begin{equation}\label{proof-density-36-1}
\epsilon_i<\frac{d''}{14N}|Q|.
\end{equation}
Under such choices of numbers $l,$ $k$ and $\epsilon_i$ $(i\in\mathcal{I}_2)$, we finally have from (\ref{proof-density-1}) and the definition of $d''$ that
\[
\begin{split}
|\Gamma^Q_{w_\eta}-\Gamma^Q_{w^*}| & \le |\Gamma^Q_{w_\eta}-\Gamma^Q_{w}|+|\Gamma^Q_{w}-\Gamma^Q_{w^*}| \\
& <d'' +|\Gamma^Q_{w}-\Gamma^Q_{w^*}|=\frac{\epsilon}{4}+ \frac{|\Gamma^Q_{w}-\Gamma^Q_{w^*}|}{2} <\frac{\epsilon}{2};
\end{split}
\]
hence, our goal (\ref{goal-5}) of this substep is indeed achieved.

\subsubsection*{\textbf{Substep 4-4}} We now prove
\begin{equation}\label{goal-6}
\int_{Q} \mathrm{dist}\big(\nabla w_\eta(x,t) ,K_b(u_\eta(x,t))\big)\,dxdt\le \delta|Q|.
\end{equation}

Let $i\in\mathcal{I}_2$.
By (c), we have $|Q_{i,1}|<2\epsilon_i$. Let $(x,t)\in Q_{i,1}$ and $r\in[r_1,r_2]$; then by (\ref{goal-4}), we have, at the point $(x,t)$,
\[
\begin{split}
  \bigg|    \begin{pmatrix} (u_\eta)_x  & (u_\eta)_t  \\ (v_\eta)_x  & (v_\eta)_t  \end{pmatrix}  &- \begin{pmatrix} g^+(r) & (u_\eta)_t \\ u_\eta & r \end{pmatrix}\bigg| \\
& =\big|\big((u_\eta)_x, (v_\eta)_t \big)-\big(g^{+}(r), r\big)\big| \\
&  \le \mathrm{diam}(\tilde{U}),
\end{split}
\]
where $\begin{pmatrix} g^+(r) & (u_\eta)_t \\ u_\eta & r \end{pmatrix}\in K_b(u_\eta)$.
Thus, we have
\begin{equation*}
\int_{Q_{i,1}} \mathrm{dist}\big(\nabla w_\eta(x,t),K_b(u_\eta(x,t))\big)\,dxdt \le 2\epsilon_i \mathrm{diam}(\tilde{U})\le\frac{\delta}{Nk}|Q|,
\end{equation*}
where we let
\begin{equation}\label{proof-density-38}
\epsilon_i\le \big(\mathrm{diam}(\tilde{U})\big)^{-1}\frac{\delta}{2Nk}|Q|.
\end{equation}
Having this choice for all $i\in\mathcal{I}_2$, we get
\begin{equation}\label{proof-density-39}
\sum_{i\in\mathcal{I}_2}\int_{Q_{i,1}} \mathrm{dist}\big(\nabla w_\eta(x,t),K_b(u_\eta(x,t))\big)\,dxdt \le  \frac{\delta}{k}|Q|.
\end{equation}

Let $i\in\mathcal{I}_2$ and $(x,t)\in Q_{i,2}$; then $(\varphi_i)_x(x,t)\in\{\lambda_{i}^+,-\lambda_{i}^-\}$. Suppose $(\varphi_i)_x(x,t)=-\lambda_{i}^-.$ By (\ref{proof-density-10}), we can  choose a number $\bar{r}_i\in[r_1,r_2]$ such that
\[
\big| (s_i-\lambda_{i}^-,\gamma_i) - \big(g_-(\bar{r}_i),\bar{r}_i\big) \big|=\frac{\delta}{k}.
\]
Then it follows from (a), (\ref{proof-density-6}) and the previous equality that  at the point $(x,t)$,
\[
\begin{split}
  \bigg|  & \begin{pmatrix} (u_\eta)_x  & (u_\eta)_t \\ (v_\eta)_x  & (v_\eta)_t  \end{pmatrix}
 - \begin{pmatrix} g^-(\bar{r}_i) & (u_\eta)_t \\ u_\eta  & \bar{r}_i \end{pmatrix}\bigg| \\
&  =\big|\big((u_\eta)_x, (v_\eta)_t\big)-\big(g^-(\bar{r}_i), \bar{r}_i\big)\big| \\
& \le |(u_x, v_t) - (s_i,\gamma_i)|  + \big| (s_i-\lambda_{i}^-,\gamma_i) - \big(g^-(\bar{r}_i), \bar{r}_i\big) \big| + |(\psi_i)_t|  \\
&  \le \frac{\delta}{4k} + \frac{\delta}{k} + \epsilon_i \le \frac{3\delta}{2k},
\end{split}
\]
where we let $\epsilon_i$ satisfy (\ref{proof-density-22}). The same can be shown in case of $(\varphi_i)_x(x,t)=\lambda_{i}^+$; we omit the detail. Thus we get
\[
\int_{Q_{i,2}} \mathrm{dist}\big(\nabla w_\eta(x,t),K(u_\eta(x,t))\big)\,dxdt \le \frac{3\delta}{2k} |Q_{i,2}|,
\]
and so
\begin{equation}\label{proof-density-40}
\sum_{i\in\mathcal{I}_2}\int_{Q_{i,2}} \mathrm{dist}\big(\nabla w_\eta(x,t),K(u_\eta(x,t))\big)\,dxdt \le \frac{3\delta}{2k} |Q|.
\end{equation}

Gathering all of (\ref{proof-density-2}), (\ref{proof-density-5}), (\ref{proof-density-9}), (\ref{proof-density-39}) and  (\ref{proof-density-40}), we obtain
\begin{equation*}
\begin{split}
\int_{Q} \mathrm{dist}&\big(\nabla w_\eta(x,t),K_b(u_\eta(x,t))\big)\,dxdt \\
= & \int_{Q\setminus\bar{G}} \mathrm{dist}\big(\nabla w(x,t),K_b(u(x,t))\big)\,dxdt \\
& + \int_{G\setminus(\cup_{i=1}^N\bar{Q}_i)} \mathrm{dist}\big(\nabla w(x,t),K_b(u(x,t))\big)\,dxdt \\
& + \int_{\cup_{i\in\mathcal{I}^+_1\cup \mathcal{I}^-_1}Q_i} \mathrm{dist}\big(\nabla w(x,t),K_b(u(x,t))\big)\,dxdt \\
& + \sum_{i\in\mathcal{I}_2}\int_{Q_{i,1}} \mathrm{dist}\big(\nabla w_\eta(x,t),K_b(u_\eta(x,t))\big)\,dxdt \\
& + \sum_{i\in\mathcal{I}_2}\int_{Q_{i,2}} \mathrm{dist}\big(\nabla w_\eta(x,t),K_b(u_\eta(x,t))\big)\,dxdt \\
\le & \Big(\frac{\delta}{k}+\frac{\delta}{k} + \frac{5\delta}{4k}+\frac{\delta}{k} +\frac{3\delta}{2k}\Big)|Q| = \frac{23\delta}{4k}|Q|<\delta|Q|,
\end{split}
\end{equation*}
where we let $k\ge6$; hence, (\ref{goal-6}) is satisfied.

\subsubsection*{\textbf{Substep 4-5}}
In conclusion, let us take an $l\in\N$ satisfying (\ref{proof-density-36}), choose a $k\in\N$ to fulfill $k\ge 6$ (so (\ref{proof-density-k-extra-1}) holds), (\ref{proof-density-k-extra-2}), (\ref{proof-density-27}) and (\ref{proof-density-31}), and select numbers $\epsilon_i>0$ $(i\in\mathcal{I}_2)$ such that (\ref{proof-density-15}), (\ref{proof-density-18}), (\ref{proof-density-22}), (\ref{proof-density-29}), (\ref{proof-density-36-1})  and (\ref{proof-density-38}) are satisfied. Then we have (\ref{goal-1}), (\ref{goal-2}), (\ref{goal-3}), (\ref{goal-4}), (\ref{goal-5}) and (\ref{goal-6}); that is, the function $w_\eta=(u_\eta,v_\eta)$ indeed satisfies all the desired properties in (\ref{proof-density-goal}).

The proof of Lemma \ref{lem:density lemma}  is now complete.

\end{document}